\documentclass[twoside,11pt]{article}

\usepackage{amsmath}
\usepackage{amssymb}
\usepackage{amsthm}
\usepackage{amsfonts}
\usepackage{amscd}

\begin{document}

\unitlength0.5cm

\newcommand{\Co}{{\mathbb C}}

\newcommand{\lV}{{{}^*V}}
\newcommand{\lU}{{{}^*U}}
\newcommand{\lW}{{{}^*W}}

\newcommand{\G}{\ensuremath{{\cal G}}}
\newcommand{\M}{\ensuremath{{\cal M}}}
\newcommand{\N}{\ensuremath{{\cal N}}}
\newcommand{\A}{\ensuremath{{\cal A}}}
\newcommand{\dG}{\ensuremath{\hat{{\cal G}}}}
\newcommand{\D}{\ensuremath{{\cal D}}}
\newcommand{\B}{{\cal B}}   
\newcommand{\hG}{\hat\G}
\renewcommand{\H}{{\cal H}}
\newcommand{\C}{{\cal C}}
\newcommand{\F}{{\cal F}}
\newcommand{\R}{{\cal R}}

\newcommand{\GG}[1][]{{\rm \boldsymbol{\Gamma}}^{#1}}
\newcommand{\GGt}[1][]{{\rm \tilde{\boldsymbol{\Gamma}}}^{#1}}
\newcommand{\JJ}[1][]{{\rm \bf J}^{#1}}

\newcommand{\LLt}[1][]{{\rm \bf \tilde{L}}^{#1}}
\newcommand{\RRt}{{\rm \bf \tilde{R}}}
\newcommand{\TTt}[1][]{{\rm \bf \tilde{T}}^{#1}}
\newcommand{\Lt}{\tilde{L}}
\newcommand{\FF}[1][]{{\rm \bf F}^{#1}}
\newcommand{\XX}{{\rm \bf X}}
\newcommand{\RRp}{{\rm \bf R'}}

\newcommand{\Fb}{\bar{F}}
%
\newcommand{\Fi}[2][]{\phi^{#1}_{#2}}
\newcommand{\Fii}[2][]{(\phi_{#2}^{-1})^{#1}}
\newcommand{\vi}{\varphi}
\newcommand{\Fib}[2][]{\bar{\phi}^{#1}_{#2}}
\newcommand{\Pb}{\bar{\Psi}}
\newcommand{\la}{\lambda}
\newcommand{\La}{\Lambda}
\newcommand{\ep}{\epsilon}
\newcommand{\om}{\omega}
\newcommand{\Om}{\Omega}
\newcommand{\de}{\delta}
\newcommand{\laR}{{\la}_R}
\newcommand{\rhoL}{{\rho}_L}
\newcommand{\Psh}{\hat{\Psi}}
\newcommand{\Phh}{\hat{\bar{\Psi}}}
\newcommand{\gam}{\gamma}
\newcommand{\rhot}{{\tilde{\rho}}}
\newcommand{\lat}{{\tilde{\la}}}
\newcommand{\cop}{\Delta}
\newcommand{\copd}{\Delta^D}
\newcommand{\copf}{\Delta_f}

\newcommand{\me}{m_{(1)}}
\newcommand{\mo}{m_{(0)}}
\newcommand{\mme}{m_{(-1)}}
\newcommand{\mmz}{m_{(-2)}}
\newcommand{\mz}{m_{(2)}}

\newcommand{\vie}{\vi_{(1)}}
\newcommand{\viz}{\vi_{(2)}}
\newcommand{\vid}{\vi_{(3)}}

\newcommand{\psie}{\psi_{(1)}}
\newcommand{\psiz}{\psi_{(2)}}
\newcommand{\psid}{\psi_{(3)}}
\newcommand{\Si}{S^{-1}}
\newcommand{\tp}{\otimes}

\newcommand{\pfeil}{\longrightarrow}
\newcommand{\e}{{\bf 1}}
\newcommand{\eG}{{\e}_{\G}}
\newcommand{\eA}{{\e}_{\A}}
\newcommand{\eB}{{\e}_{\B}}
\newcommand{\eC}{{\e}_{\C}}
\newcommand{\edG}{{\e}_{\dG}}
\newcommand{\eM}{{\e}_{\M}}
\newcommand{\eN}{{\e}_{\N}}

\newcommand{\lpa}{\langle}
\newcommand{\rpa}{\rangle}
%
\newcommand{\wir}{\triangleright}

\newcommand{\gotH}{{\mathfrak H}} 
\newcommand{\goth}{{\mathfrak h}} 
\newcommand{\gotK}{{\mathfrak K}} 
\newcommand{\gotL}{{\mathfrak L}} 
%
\newcommand{\bo}{\boxtimes}
\newcommand{\so}{\odot}
\newcommand{\dol}{\gtrdot}
\newcommand{\dor}{\lessdot}

\newcommand{\id}{{\rm id}}
\newcommand{\idM}{{\id}_{\M}}
\newcommand{\idA}{{\id}_{\A}}
\newcommand{\idMR}{{\id}_{\MR}}
\newcommand{\idML}{{\id}_{\ML}}
\newcommand{\idG}{{\id}_{\G}}
\newcommand{\iddG}{{\id}_{\dG}}

\numberwithin{equation}{section}

\newcommand{\ncm}{\newcommand}
\ncm{\rncm}{\renewcommand}


\newtheorem{theorem}{Theorem }[section]
\newtheorem{lemma}[theorem]{Lemma }
\newtheorem{proposition}[theorem]{Proposition }
\newtheorem{corollary}[theorem]{Corollary }
\newtheorem{definition}[theorem]{Definition }

\ncm{\Theorem}[2]{\begin{theorem} \label{Thm #1} {\it #2} \end{theorem}}
\ncm{\thm}[1]{Theorem \ref{Thm #1}}

\ncm{\Definition}[2]{\begin{definition} \label{Def #1}{\rm #2}\end{definition}}
\ncm{\defi}[1]{Definition \ref{Def #1}}

\ncm{\Lemma}[2]{\begin{lemma} \label{Lem #1} {\it #2} \end{lemma}}
\ncm{\lem}[1]{Lemma \ref{Lem #1}}

\ncm{\Proposition}[2]{\begin{proposition}\label{Prop #1}{\it #2}
                      \end{proposition}}
\ncm{\prop}[1]{Proposition \ref{Prop #1}}

\ncm{\Corollary}[2]{\begin{corollary}\label{Cor #1} {\it #2} \end{corollary}}
\ncm{\cor}[1]{Corollary \ref{Cor #1}}

\ncm{\Thm}{\Theorem}
\ncm{\Def}{\Definition}
\ncm{\Lem}{\Lemma}
\ncm{\Prop}{\Proposition}
\ncm{\Cor}{\Corollary}

\ncm{\setc}[1]{\setcounter{equation}{#1}}
\rncm{\sec}{\setc{0}\section}
\ncm{\no}[1]{(\ref{#1})}
\ncm{\Eq}[1]{Eq.\ \no{#1}}

\ncm{\bsn}{\bigskip\noindent}

%
\def\Ga{\boldsymbol{\Gamma}}
\ncm{\RR}{{\bf R}}
\ncm{\LL}{{\bf L}}
\ncm{\TT}{{\bf T}}
\ncm{\DD}{{\bf D}}
\newcommand{\CC}{{\mathbb C}}
\newcommand{\ZZ}{{\mathbb Z}}

\ncm{\bea}{\begin{eqnarray}}
\ncm{\eea}{\end{eqnarray}}
\ncm{\beanon}{\begin{eqnarray*}}
\ncm{\eeanon}{\end{eqnarray*}}
\ncm{\ba}{\begin{array}}
\ncm{\ea}{\end{array}}

%
\def\<cros{\ltimes}
\def\crosl{\<cros\!\!_\l\,\,}
\def\croslt{{\<cros\!_{\tilde{\l}}}\,}
\def\croslC{\<cros\!\!_{\l_\C}\,\,}
\def\croslB{\<cros\!\!_{\l_\B}\,\,}
\def\cros{\rtimes}
\def\>cros{\cros}
\def\crosr{\cros\!\!_\rho\,\,}
\def\crosrA{\cros\!\!_{\rho_\A}\,\,}
\def\crosrB{\cros\!\!_{\rho_\B}\,\,}
\def\crosrt{{\cros\!_{\tilde{\rho}}}\,}
\def\arr{\rightharpoonup}
\def\arl{\leftharpoonup}

\def\re{{\,\hbox{$\textstyle\triangleright$}\,}}
\def\li{{\,\hbox{$\textstyle\triangleleft$}\,}}

\def\reli{\bowtie}
\def\reliD{\reli\!\!_D\,\,}
\def\relihD{\reli\!\!_{\hat D}\,\,}
\def\relid{\reli\!\!_\delta\,\,}
\def\relidop{\reli\!\!_{\delta_{op}}\,\,}
\def\relidp{\reli\!\!_{\delta'}\,\,}
\def\relir{\reli\!\!_\rho\,\,}
\def\relidr{\reli\!\!_{\delta_r}\,\,}
\def\relidl{\reli\!\!_{\delta_l}\,\,}
\def\relilr{{}_\la\!\!\reli\!\!_\rho\,\,}
\def\doubcross{\reli\!\!_{double}\,\,}

\def\DeltaM{\Delta_\M}
\def\Re{\ -\!\!\!>}
\def\Li{<\!\!\!-\ }
%
%
\def\Rep{\mbox{Rep}\,}
\def\Ker{\mbox{Ker}\,}
\def\End{\mbox{End}\,}
\def\Ad{\mbox{Ad}\,}
\def\Hom{\mbox{Hom}}

\def\bra{\langle}
\def\ket{\rangle}
\def\o{\otimes}
\ncm{\one}{{\bf 1}}
%
%
\def\e{\epsilon}
\def\x{\times}
\def\al{\alpha}
\def\be{\beta}
\def\r{\rho}
\def\s{\sigma}
\rncm{\l}{\lambda}
\ncm{\0}{_{(0)}}
\ncm{\1}{_{(1)}}
\ncm{\2}{_{(2)}}
\ncm{\3}{_{(3)}}
\ncm{\4}{_{(4)}}
\ncm{\5}{_{(5)}}
\ncm{\6}{_{(6)}}
\ncm{\7}{_{(7)}}
\ncm{\8}{_{(8)}}
\ncm{\9}{_{(9)}}
\def\m1{_{(-1)}}
\def\n2{_{(-2)}}

\def\lri{$\lambda\rho$-intertwiner}

%

\title{ \bf Diagonal Crossed Products by\\
Duals of Quasi--Quantum Groups}

\renewcommand{\thefootnote}{\fnsymbol{footnote}}
\author{{\sc Frank Hau{\ss}er\footnotemark[1]}\hspace{1ex} and {\sc Florian
    Nill\footnotemark[2]} 
\\
{\small Freie Universit\"at Berlin,
Institut f\"ur theoretische Physik,}\\
{\small Arnimallee 14, D-14195 Berlin, Germany}}

\maketitle

\footnotetext[1]{
supported by DFG, SFB 288 {\em Differentialgeometrie und
  Quantenphysik},\quad
e-mail: hausser@physik.fu-berlin.de}
\footnotetext[2]{
supported by DFG, SFB 288 {\em Differentialgeometrie und
  Quantenphysik},\quad
e-mail:  nill@physik.fu-berlin.de}

\renewcommand{\thefootnote}{\arabic{footnote}}
\setcounter{footnote}{0}

\begin{abstract}
A two-sided coaction $\delta:\M\to\G\o\M\o\G$ of a Hopf algebra
 $(\G,\Delta ,\e,S)$
on an associative algebra $\M$ is an algebra map of the form
$\delta=(\l\o\id_\M)\circ\rho=(\id_\M\o\rho)\circ\l$, where $(\l,\rho)$
is a commuting pair of left and right $\G$-coactions on $\M$,
respectively.
Denoting the associated commuting right and left actions of the dual
Hopf algebra $\hat\G$ on $\M$ by
$\li$ and $\re$, respectively, we define the {\em diagonal crossed
  product} $\M\reli\hat\G$ to be the algebra generated by $\M$ and
$\hat\G$ with relations given by
$$
 \varphi\, m = (\varphi\1\re m\li
\hat S^{-1}(\varphi\3))\, \varphi\2,\ m\in\M,\,\varphi\in\hat\G\ .
$$
We give a natural generalization of this construction to the case
where $\G$ is a quasi--Hopf algebra in the sense of Drinfeld and, more
generally, also in the sense of Mack and Schomerus (i.e., where the
coproduct $\Delta$ is non-unital).
In these cases our diagonal crossed product will still be an
associative algebra structure on $\M\o\hat\G$ extending
$\M\equiv\M\o\hat\one$, even though
the analogue of an ordinary crossed product
$\M\cros \hat\G$ in general is not well defined as an associative
algebra.\\
Applications of our formalism include the field algebra constructions
with quasi-quantum group symmetry given by Mack and Schomerus [MS,S]
as well as the  
formulation of Hopf spin chains or lattice current algebras based on
truncated quantum groups at roots of unity.\\
In the case $\M=\G$  and $\l=\rho=\Delta$ we obtain an explicit
definition of
the quantum double $\D(\G)$ for quasi--Hopf algebras $\G$, which
before had been described in the form of an implicit Tannaka-Krein
reconstruction procedure by S. Majid [Ma2]. We prove that
$\D(\G)$ is itself a (weak) quasi-bialgebra and that any
diagonal crossed product $\M\reli\hat\G$ naturally admits a
two-sided $\D(\G)$-coaction.
In particular, the above mentioned lattice models
always admit the quantum double $\D(\G)$ as a localized cosymmetry,
generalizing results of Nill and Szlach\'anyi [NSz1].
A complete proof that $\D(\G)$ is even a (weak) quasi-triangular
quasi-Hopf algebra will be given in a separate paper [HN1].
\end{abstract}

\newpage

\renewcommand{\baselinestretch}{0.6}
\footnotesize
\tableofcontents
\newpage
\renewcommand{\baselinestretch}{1}
\normalsize

%
\sec{Introduction and summary of results}

During the last decade quantum groups have become the most fashionable
candidates describing the symmetry in low dimensional quantum field
theories (QFT)\footnote{see [BaWi], [BL], [G], [DPR], [FGV], [FrKe],
  [MS1,2], [MoRe], [M\"u], [PSa1], [ReSm], [Sz,V]}
or lattice models\footnote{see [AFFS], [AFSV], [AFS], [ByS], [Fa],
  [FG], [KaS], [NSz1,2], [P], [PSa2], [SzV]}.
Here, in an axiomatic approach by a ``symmetry algebra'' $\G$ one means a
$*$-algebra acting on the Hilbert space of physical states $\H$,
such that
\begin{itemize}
\item
observables and space-time translations commute with $\G$,
\item
charge creating fields fall into multiplets transforming covariantly
under the action of $\G$,
\item
equivalence classes of irreducible representations of $\G$ are in
one-to-one correspondence with the Doplicher-Haag-Roberts (DHR)
superselection sectors of the observable algebra $\A$, such that also
the fusion rules of $\Rep_{\!\mbox{\sc\footnotesize DHR}}\,\A$ 
and $\Rep\G$ coincide.
\end{itemize}
It is well known that the results of Doplicher-Roberts [DR] 
characterizing $\G$ as a compact group (or the associated group 
algebra) break down in low dimensions due to the appearance of 
braid statistics. It was soon realized that at least for 
rational theories (i.e. with a finite number of sectors) 
quantum groups are also ruled out, unless all sectors have 
integer statistical dimensions (see e.g. [FrKe] for a review or 
[N2] for a specific discussion of q-dimensions in finite 
quantum groups).

Based on the theory of quasi-Hopf algebras introduced by 
Drinfel'd [Dr2], G. Mack and V. Schomerus [MS] have proposed the 
notion of {\em weak quasi-Hopf algebras} $\G$ as appropriate 
symmetry candidates, where ``weak'' means that the tensor 
product of two ``physical'' representations of $\G$ may also 
contain ``unphysical'' subrepresentations (i.e. of q-dimension 
$\le 0$), which have to be discarded. Examples are semisimple 
quotients of q-deformations of classical groups at $q=$ roots 
of unity.

In this way non-integer dimensions could successfully be 
incorporated. The price to pay was that now commutation 
relations of $\G$-covariant charged fields involve operator 
valued $R$-matrices and, more drastically, the operator product 
expansion for $\G$-covariant multiplets of charged fields 
involves non-scalar coefficients with values in $\G$.
Thus, the analogue of the ``would-be'' DHR-field algebra 
$\F$ is no longer algebraically closed. Instead, Mack and 
Schomerus have proposed a new ``covariant product'' for charged 
fields, which does not lead outside of $\F$, but which is no 
longer associative. In [S] Schomerus has analyzed this scenario 
somewhat more systematically in the framework of DHR-theory, 
showing that a weak quasi-Hopf algebra $\G$ and a field 
``algebra'' $\F$ may always be constructed such that the 
combined algebra $\F\vee\G$ is associative and satisfies all 
desired properties, except that $\F\subset\F\vee\G$ is only a 
linear subspace but not a subalgebra.
Technically, the reason for this lies in the fact that the dual 
$\hat\G$ of a quasi-Hopf algebra is not an associative algebra.

\bigskip
One should also remark at this point that the above 
reconstruction of $\G$ from the category of DHR-endomorphisms is not unique.
Also, in a more mathematical framework a general Tannaka-Krein like 
reconstruction theorem for quasi-Hopf algebras has been 
obtained by S. Majid [M1] and for weak quasi--Hopf algebras by [H\"a].

\bigskip
To study quantum symmetries on the lattice in an axiomatic approach,
K.\ Szlach\'anyi and P.\ Vecsernyes [SzV]  have proposed an ``amplified'' 
version of the DHR-theory, which also applies to locally finite 
dimensional lattice models. This setting has been further 
developed by [NSz1,2], where based on the example of Hopf spin 
chains the authors proposed the notion of a {\em universal 
localized cosymmetry} $\rho:\A\to\A\o\G$, incorporating all 
sectors $\rho_I$ of $\A$ via $\rho_I=(\id_\A\o\pi_I)\circ\rho,\ 
\pi_I\in\Rep\G$.
In the specific example studied by [NSz1,2] $\G$ was given by a 
quantum double and the cosymmetry $\rho$ was given by a {\em 
coaction} of $\G$ on $\A$. Related results have later been 
obtained for lattice current algebras [AFFS], the later 
actually being a special case of the Hopf spin chains of [NSz1,2] 
(see [N1] and Sect. 11.3). The analogue of a DHR--field algebra for
these models is now given by the standard crossed product $\F \equiv
\A\>cros\dG$ [NSz2], where $\dG$ is the Hopf algebra dual to $\G$.

Now the methods and results of these works were still 
restricted to ordinary Hopf algebras and therefore to integer 
dimensions. To formulate lattice current algebras at roots of 
unity one may of course identify them with the boundary part of 
lattice Chern-Simons algebras [AGS,AS] defined on a disk. 
Nevertheless, it remains unclear whether and how for $q=$root 
of unity the structural 
results of [AFFS] survive the truncation to the semi-simple 
(``physical'') quotients. Similarly, the generalizations of the 
model, the methods and the results of [NSz] to weak quasi 
quantum groups are by no means obvious. In particular one would 
like to know whether and in what sense in such models universal 
localized cosymmetries $\rho:\A\to\A\o\G$ still provide 
coactions and whether $\G$ would still be (an analogue of) a 
quantum double of a quasi-Hopf algebra, possibly in the sense 
recently described by Majid [M2].

\bigskip
In this work we present a theory of left, right and two-sided coactions of a 
(weak) quasi-Hopf algebra $\G$ on an associative algebra $\M$. 
Based on these structures we then provide a new construction of 
what we call the {\em diagonal crossed product} $\M\reli\hG$, 
which we will show to be the appropriate mathematical structure 
underlying all constructions discussed above. In particular, 
$\M\reli\hG$ will always be an associative algebra extending 
$\M\equiv\M\reli\hat\one$. On the other hand, the linear 
subspace $\one_\M\reli\hG$ will in general not be a subalgebra 
of $\M\reli\hG$, unless $\G$ is an ordinary (i.e. 
coassociative) Hopf algebra.

The basic idea for this construction comes from generalizing the
relations defining the quantum double.
To this end we start from an algebra $\M$ equipped with a
(quasi-)commuting pair of  
right and left $\G$-coactions, $\rho:\M\to\M\o\G$ and 
$\l:\M\to\G\o\M$ and denote 
$\delta_\ell:=(\l\o\id)\circ\rho$ and 
$\delta_r:=(\id\o\rho)\circ\l$ as the associated equivalent {\em two-sided 
coactions}. In the simplest case of $\G$ being an 
ordinary Hopf algebra and $(\l,\rho)$  being strictly commuting
(i.e. $\delta_\ell =\delta_r$) this amounts to providing a 
commuting pair of left and right Hopf module actions 
$\re:\hG\o\M\to\M$ (dual to $\rho$) and $\li:\M\o\hG\to\M$ 
(dual to $\l$) of the dual Hopf algebra $\hG$ on $\M$. In this 
case our diagonal crossed product $\M\reli\hG$ is defined to be 
generated by $\M$ 
and $\hG$ as unital subalgebras with commutation relations 
given by ($\hat S:\hG\to\hG$ being the antipode)
\begin{equation}
\label{0.1}
 \varphi\, m = (\varphi\1\re m\li
\hat S^{-1}(\varphi\3))\, \varphi\2,\ m\in\M,\,\varphi\in\hat\G\ .
\end{equation}
Note that for $\M=\G$ and $\rho=\l=\Delta$ the coproduct on $\G$,
 these are the defining relations of the quantum double $\D(\G)$ [Dr1],
 and therefore $\G\reli\hG=\D(\G)$.
Introducing the ``generating matrix''
$$
\Ga:=\sum_\mu e_\mu\o e^\mu\in\G\o\hG\subset\G\o(\M\reli\hG),
$$
where $e_\mu\in\G$ is a basis with dual basis $e^\mu\in\hG$, 
\Eq{0.1} is equivalent to
\begin{equation}
\label{0.2}
\Ga\l(m)=\rho^{op}(m)\Ga,\quad\forall m\in\M.
\end{equation}
Moreover, in this case $\hG\subset\M\reli\hG$ being a unital
subalgebra is equivalent to
\begin{eqnarray}
(\e\o\id)(\Ga)&=&\one\label{0.3}
\\
\Ga^{13}\Ga^{23}&=&(\Delta\o\id)(\Ga)\label{0.4}
\end{eqnarray}
where \no{0.4} is an identity in $\G\o\G\o(\M\reli\hG)$, the 
indices denoting the embeddings of tensor factors. 
We call $\Ga$ the {\em universal normal and coherent $\l\rho$-intertwiner}
in $\G\o(\M\reli\hG)$, where normality is the property \no{0.3} and
coherence is the property \no{0.4}.
Again, for $\M=\G$ and $\M\reli\G=\D(\G)$ Eqs. \no{0.2}-\no{0.4} are 
precisely the defining relations for the  generating 
matrix $\DD\equiv\Ga_{\D(\G)}$ of the quantum double (see e.g. 
[N1, Lem.5.2]).

Inspired by the techniques of [AGS, AS] we show in the main 
body of this work how to 
generalize the notion of coherent $\l\rho$-intertwiners to the case 
of (weak) quasi-Hopf algebras $\G$, such that analogues of
the Eqs. \no{0.2}-\no{0.4} still serve 
as the defining relations of an associative algebra extending 
$\M\equiv\M\reli\hat\one$. We also 
show that diagonal crossed products may equivalently be 
modeled on the linear spaces $\M\o\hG$ or $\hG\o\M$ (or -- in the 
weak case -- certain subspaces thereof).

\bigskip
The basic model for this generalization is again given by $\M=\G$ with 
its natural two-sided $\G$-coactions 
$\delta_\ell:=(\Delta\o\id)\circ\Delta$ and
$\delta_r:=(\id\o\Delta)\circ\Delta$.
In this case our construction provides a definition of the 
quantum double $\D(\G)$ 
for (weak) quasi-Hopf algebras $\G$. In fact, we 
show that $\Rep\D(\G)$ coincides with 
what has been called the ``double of the category'' $\Rep\G$ in 
[M2]. Hence our definition provides a concrete realization of 
the abstract Tannaka-Krein like reconstruction of the quantum 
double given by [M2].
We also give a proof that $\D(\G)$ is a (weak) quasi-bialgebra. 
In [HN1] we will show, that $\D(\G)$ is in fact a 
(weak) quasi-triangular quasi-Hopf algebra, and there we will 
also visualize many of our (otherwise almost untraceable) algebraic 
identities in terms of graphical proofs.

 The field algebra construction of [MS,S] may also be described 
as a diagonal crossed product $\M\reli\hG$ within our formalism
by putting $\M=\A\o\G$, where $\A$ is the observable algebra. 
In this case the right $\G$-coaction $\rho:\M\to\M\o\G$ is a 
localized cosymmetry acting only on $\A$, whereas the left 
$\G$-coaction $\l:\M\to\G\o\M$ only acts on $\G$, where it is 
given by the coproduct $\Delta$, see Sect. 11.4 for a rough sketch.
A more detailed account of this within an appropriate 
von-Neumann algebraic framework will be given elsewhere.

The application of our formalism to $\G$-spin quantum chains is 
given by putting in the previous example also $\A=\G$ and 
$\rho=\Delta$, in which case 
$\M\reli\hG\cong\G\>cros\hG\<cros\G$ becomes 
a {\it two--sided crossed product}. We take this construction as
building block 
of a quantum chain living on two neighboring sites (carrying the 
copies of $\G$) joined by a link (carrying the copy of $\hG$).
We show how this construction iterates to provide a local net 
of associative algebras $\A(I)$ for any lattice interval $I$ 
bounded by sites. Generalizing the methods of [NSz] we also 
construct localized coactions of the quantum double 
$\D(\G)$ on such (weak) quasi-Hopf spin chains. Periodic 
boundary conditions for these models are again described as a 
diagonal crossed product of the open chain by a copy of $\hG$ 
sitting on the link joining the end points.
In this way we arrive at a formulation of lattice current 
algebras at roots of unity 
by adjusting the transformation rules of [N1] to the 
quasi-coassociative setting. We expect to have more detailed 
results also on these models in the near future.

\bigskip
We now describe the plan of this paper.
In Part I we prepare our language and develop our basic ideas 
within a strictly coassociative setting, putting emphasis on a
pedagogical presentation of the ``read thread''.
In Sect.\ 2 we start with reviewing the notions of left and right coactions and
crossed products. In Sect. 3 we introduce two-sided coactions 
$\delta:\M\to\G\o\M\o\G$ and
identify them with pairs $(\l,\rho)$ of commuting left and right
coactions. We then define the diagonal crossed product
$\M\,_\l\!\!\relir\hG$ and show that it may be identified with the
subalgebra of $\hG\,\,_\l\!\!\<cros\M\crosr\hG$ generated by $\M$ and the
``diagonal'' $\hat\Delta_{op}(\hG)\subset\hG\o\hG$. For
$\M=\G$ this gives the quantum double $\D(\G)$. Moreover any
diagonal crossed product $\M\reli\hG$ canonically admits a two-sided
$\D(\G)$-coaction.

In Sect.\ 4 we give examples and applications of our
formalism. Section 4.1 discusses the relation of our diagonal crossed
products with Majid's notion of {\it double crossed products}
[M3,4]. In Section 4.2 we introduce
 {\em two-sided crossed products}
$\A\>cros\hG\<cros\B$ as special examples obtained for the case
$\M=\A\o\B$, where $\rho$ acts trivially on $\B$ and $\l$ acts trivially on
$\A$.
By putting $\A=\B=\G$ and iterating, these become the building structures
of Hopf spin chains described in Section 4.3.

In Sect.\ 5 we reformulate our constructions using the generating
matrix formalism. This leads to the description of diagonal crossed
products as being generated by $\M$ and the ``matrix entries'' of a
normal coherent $\l\rho$-intertwiner $\Ga\in\G\o(\M\reli\hG)$
subjected to the relations \no{0.2}-\no{0.4}.

This is the appropriate language to be generalized to the
quasi-coassociative setting in Part~II. Sect.\ 6 starts with a short review
of Drinfeld's axioms for quasi-Hopf algebras $\G$ and some of their
properties. In Sects.\ 7 and 8 we provide a generalization of the
notions of left, right and two-sided $\G$-coactions on an associative
algebra $\M$. We also show that up to twist equivalence two-sided
$\G$-coactions are in one-to-one correspondence with {\em
  quasi-commuting} pairs of (left and right) $\G$-coactions
$(\l,\rho)$.

In Sect.\ 9 we give a representation theoretic interpretation of these
concepts by showing that they give rise to left, right and two-sided
``actions'' of $\Rep\G$ on $\Rep\M$, respectively. By this we mean that
representations $\gamma$ of $\M$ may be ``tensored'' with
representations $\pi$ of $\G$ to yield again representations of $\M$,
such that an analogue of McLane's natural associativity constraints for
monoidal categories is satisfied.
This formalism allows to motivate and explain many otherwise almost
untraceable algebraic identities in terms of commuting diagrams. 
Yet, the most technical parts of these proofs are deferred to an Appendix.

In Sect.\ 10 we proceed to the generalized construction of diagonal
crossed products associated with any two-sided $\G$-coaction
$\delta:\M\to\G\o\M\o\G$. Sect.\ 10.1 gives the associativity proof and
states the equivalence of a ``left'' and a ``right'' convention for
this construction. In Sect.\ 10.2 we describe these two conventions by
the associated generating matrices $\LL$ and $\RR$, called {\em left
  (right) diagonal $\delta$-implementers}.
In Sect 10.3 we generalize the notion of coherent
$\l\rho$-intertwiners $\GG$ to the quasi-coassociative setting and show that
they are always in one-to-one correspondence with the above coherent left (or
right) diagonal $\delta$-implementers. This provides a proof of our
main result in Sect.\ 10.4, showing that diagonal crossed products are
uniquely (up to equivalence) described by a quasi-coassociative
version of the relations \no{0.2}-\no{0.4}, see Theorem II on page
\pageref{seite}.
In Sect.\ 11 we provide examples and applications by discussing the
quantum double $\D(\G)$, the two-sided crossed products
$\A\cros\hG\<cros\B$ and their appearance in quasi-Hopf spin chains
and lattice current algebras, and finally a reformulation of the
Mack-Schomerus scenario within our formalism.

Finally, in Part III we generalize all results to weak quasi-Hopf
algebras $\G$.  

\bigskip
In conclusion we point out that most of our algebraic constructions
are based on representation categorical concepts and may therefore
also be visualized by graphical proofs. We will put more emphasis on
this technique in [HN1] when providing further results on the quantum
double $\D(\G)$.

We also remark that without mentioning explicitly at every instance
the (weak, quasi) Hopf algebras $\G$ are always supposed to be finite
dimensional. Although we believe that many aspects of our formalism
would also carry over to an infinite dimensional setting, we don't
consider it worthwhile to discuss this complication at present.
\newpage
More importantly, in applications to quantum physics one should extend our
formalism to incorporate $C^*$- or von-Neumann algebraic structures,
which we will come back to in the near future
when discussing our examples in more detail.

\part{The Coassociative Setting}

To strip off all technicalities from the main ideas, in this first
part we restrict ourselves to strictly coassociative Hopf algebras $\G$.
Throughout by an algebra we will mean an associative unital algebra over $\CC$
and unless stated differently all algebra morphisms are supposed to be
unit preserving.
Given an algebra $\M$, two algebra extensions $\M_1\supset\M$ and
$\M_2\supset\M$ are called {\em equivalent}, if there exists an
isomorphism $\M_1\cong\M_2$ restricting to the identity on $\M$.
The main result we are going for is given by

\vspace{0.25cm}
  \begin{sloppypar} \noindent{\bf Theorem I}
{\it 
Let $(\G,\Delta,\e,S)$ be a finite dimensional Hopf algebra and let
$(\l,\rho)$ be a commuting pair
of (left and right) $\G$--coactions on an associative algebra $\M$.
\begin{enumerate}
\item Then there exists a unital associative algebra extension $\M_1
  \supset \M$ together with a linear map $\Gamma : \dG\pfeil \M_1$
  satisfying the following universal property: \\
  $\M_1$ is algebraically generated by $\M$ and $\Gamma(\dG)$ and for
  any algebra map $\gamma :\M\pfeil\A$ into some target algebra $\A$
  the relation 
  \begin{equation}
    \label{1.00}
   \gam_T (\Gamma(\vi)) = (\vi\tp\id)(\TT),\quad \vi \in \dG
  \end{equation}
provides a one--to--one correspondence between algebra maps
$\gam_T : \M_1 \pfeil \A$ extending $\gam$ and elements $\TT \in
\G\tp\A$ satisfying $(\e\o\id_\A)(\TT)=\one_\A$ and
\begin{align}
  \label{1.01}
\TT \, \l_\A(m) &= \rho_\A^{op}(m) \,\TT, \quad \forall m \in \M \\
  \label{1.02}
\TT^{13}\TT^{23} &= (\cop \tp\id_\A)(\TT),
\end{align}
where $\l_\A (m) := (\id_\G\tp\gam)(\l(m))$ and
$\rho_\A (m) := (\gam\o\id_\G)(\rho(m))$.
\item If $\M\subset \tilde{M}_1$ and $\tilde{\Gamma}:
  \dG\pfeil\tilde{\M}_1$ satisfy 
  the same universality property as in part 1.), then there exists a
  unique algebra isomorphism $f: \M_1\pfeil \tilde{\M}_1$ restricting
  to the identity on $\M$, such that $\tilde{\Gamma} = f\circ \Gamma$ 
\item The linear maps
  \begin{align}
    \label{1.03}
   \mu_L : \dG\tp\M \ni (\vi\tp m) &\mapsto \Gamma (\vi) m \, \in \M_1\\
    \label{1.04}
   \mu_R : \M\tp\dG \ni (m\tp \vi)&\mapsto m \Gamma (\vi) \, \in \M_1
  \end{align}
provide isomorphisms of vector spaces.
\end{enumerate}
  }
\end{sloppypar}\vspace{0.25cm}
Putting $\Ga : = e_\mu \tp \Gamma (e^\mu) \in \G\tp\M_1$ Theorem II
implies that $\Ga$ itself satisfies the defining relations
\eqref{1.01} and \eqref{1.02}, see also \no{0.2} and \no{0.4}. 
We call $\Ga$ the {\it universal
  $\l\rho$--intertwiner} in $\M_1$.
We emphasize that once being stated Theorem I almost appears
trivial. Its true power only arises when generalized to the
quasi-coassociative setting in Part II. 
Note that part 2.\ of Theorem I implies that the algebraic structures
induced on 
$\hG\o\M$ and $\M\o\hG$ via $\mu_{L/R}^{-1}$ from $\M_1$ are uniquely
fixed. We denote them as
\begin{align}
   \label{1.05}
    \dG\relilr \M &\equiv \mu^{-1}_L (\M_1) \\
   \label{1.06}
     \M\relilr \dG &\equiv \mu^{-1}_R (\M_1).
 \end{align}
and call them {\em left (right) diagonal crossed products},
respectively. (The right diagonal version \no{1.06} will be shown to
be the one described by \Eq{0.1}).

To actually prove Theorem I we will first construct the algebras
$\M\relilr \dG$ and $\dG\relilr \M$ explicitly as equivalent
extensions of $\M$ in \prop{1.7} and \cor{1.11}, respectively. The
description in terms of $(\l,\rho)$-intertwiners will then be given in
Sect.\ 5. In order to carefully prepare the much more complicated
quasi-coassociative scenario we deliberately present these arguments 
in rather elementary steps.

\sec{Coactions and crossed products}

To fix our conventions and notations we start with shortly reviewing some basic
notions on Hopf module actions, coactions and crossed products. For 
full textbook
treatments see e.g. [A,M3,Sw].

\bsn
Let $\G$ and $\hat\G$ be a dual pair of finite dimensional Hopf algebras. We
denote elements of $\G$ by Roman letters $a,b,c,\dots$ and elements of $\hat\G$
by Greek  letters $\varphi,\psi,\xi,\dots$. The units are denoted by
$\one\in\G$ and $\hat\one\in\hat\G$. Identifying
$\hat{\hat\G}=\G$,  the dual
pairing $\G\otimes\hat\G\to \CC$ is written as
$$
\bra a|\psi\ket \equiv \bra\psi|a\ket \in\CC,\quad
a\in\G,\psi\in\hat\G\ .
$$
We denote $\Delta :\G\to\G\otimes\G$ the coproduct, $\ep :\G\to\CC$
the counit and $S:\G\to\G$ the antipode. Similarly,
$\hat\Delta,\hat{\ep}$
and $\hat S$ are the structural maps on $\hat\G$. We will use the Sweedler
notation $\Delta(a)=a\1 \otimes a\2,\ (\Delta\otimes \id )(\Delta(a))\equiv
(\id \otimes \Delta)(\Delta(a))=a\1\otimes a\2\otimes a\3$, etc.  where the
summation symbol and the summation indices are suppressed. Together
with $\G$ we
have the Hopf algebras $\G_{op},\ \G^{cop}$ and $\G_{op}^{cop}$, where ``op''
refers to opposite multiplication and ``cop'' to opposite comultiplication. Note
that the antipode of $\G_{op}$ and $\G^{cop}$ is given by $S^{-1}$ and the
antipode of $\G^{cop}_{op}$ by $S$. Also, $\widehat{\G_{op}} = (\hat\G)^{cop},
\ \widehat{\G^{cop}}=(\hat\G)_{op}$ and $\widehat{\G^{cop}_{op}}
=(\hat\G)^{cop}_{op}$.
\\
A right coaction of $\G$ on an algebra $\M$ is an algebra map
$\rho:\M\to\M\otimes \G$ satisfying
\begin{eqnarray}
(\rho\otimes \id )\circ\rho &=& (\id \otimes \Delta)\circ \rho \label{1.1}\\
(\id \otimes \e)\circ\rho &=& \id  \label{1.2}
\end{eqnarray}
Similarly, a left coaction $\lambda$  is an algebra map
$\lambda:\M\to\G\otimes\M$
satisfying
\begin{eqnarray}
(\id \otimes\lambda)\circ \lambda &=& (\Delta\otimes \id )\circ\lambda \label{1.3}\\
(\e \otimes \id )\circ\lambda &=& \id  \label {1.4}
\end{eqnarray}
Obviously, after a permutation of tensor factors $\M\o\G
\leftrightarrow \G\o\M$
a left coaction of $\G$ may always be viewed as a right coaction by $\G^{cop}$
and vice versa.
Similarly as for coproducts we will also use the suggestive notations
\begin{eqnarray} \label{1.4a}
\rho(m)&=&m\0 \o m\1
\\\label{1.4b}
(\rho\o \id )(\rho(m)) \equiv (\id \o\Delta)(\rho(m))&=&m\0 \o m\1\o m\2
\\\label{1.4c}
\lambda(m) &=& m_{(-1)} \o m\0
\\
\label{1.4d}
(\Delta\o \id )\circ \lambda \equiv (\id \o\lambda)(\lambda(m))
&=&m_{(-2)} \o m\1 \o m\0
\end{eqnarray}
etc., where again summation indices and a summation symbol are suppressed.
In this way we will always have $m_{(i)}\in\G$ for $i\neq 0$ and
$m\0\in\M$.
Next, we recall that
there is a one-to-one correspondence between right (left) coactions of $\G$ on
$\M$ and left (right) Hopf module actions, respectively, of $\hat\G$
on $\M$ given for $\psi\in\hat\G$ and $m\in\M$ by
\begin{eqnarray}
\psi\re m &:=& (\id \o \psi)(\rho(m))\label{1.5}\\
m\li \psi &:=& (\psi\o \id )(\lambda(m)) \label{1.6}
\end{eqnarray}
One easily verifies the defining properties of Hopf module actions, i.e.
\begin{equation}
\begin{array}{rclcrcl}
\varphi\re(\psi\re m) &=& (\varphi\psi)\re m \qquad
&,& \qquad (m\li\varphi)\li\psi
&=& m\li(\varphi\psi)
\\
\hat\one \re m &=& m  &,& m\li\hat\one &=& m
\\
\varphi\re (mn) &=& (\varphi\1\re m)(\varphi\2\re n)
&,& (mn)\li\varphi &=&(m\li\varphi\1)(n\li\varphi\2)
\\
\varphi\re \one_\M &=& \hat\e (\varphi)\one_\M &,& \one_\M\li\varphi &=&
\one_\M \hat\e(\varphi)
\end{array}
\end{equation}
where $\varphi,\psi\in\hat\G$ and $m,n\in\M$.
As a particular example we recall the case $\M=\G$ with $\rho=\lambda=\Delta$.
In this case we denote the associated left and right actions of
$\psi\in\hat\G$
on $a\in\G$ by $\psi\arr a$ and $a\arl\psi$, respectively.

Given a right coaction $\rho:\M\to\M\o\G$ with dual left $\hat\G$-action $\re$
one defines the (untwisted) crossed product (also called smash product)
$\M\>cros \hat\G$ to be the vector space $\M\o\hat\G$ with associative algebra
structure given for $m,n\in\M$ and $\varphi,\psi\in\hat\G$ by
\begin{equation}\label{1.8}
(m\>cros\varphi)(n\>cros\psi)=(m(\varphi\1\re n)\>cros\varphi\2\psi)
\end{equation}
where we use the notation $m\cros\psi$ in place of $m\o\psi$ to
emphasize the new algebraic structure.
Then $\one_\M\>cros\hat\one$ is the unit in $\M\>cros \hat\G$ and
$m\mapsto (m\>cros\hat\one),\
\varphi\mapsto(\one_\M\cros\varphi)$ provide unital inclusions
$\M\to\M\>cros\hat\G$
and $\hat\G\to\M\>cros\hat\G$, respectively. Similarly if $\lambda:\M\to\G\o\M$ is
a left coaction with dual right action $\li$ then $\hat\G\<cros \M$ denotes
the associative algebra structure on $\hat\G\o\M$ given by
\begin{equation}\label{1.9}
(\varphi\<cros m)(\psi\<cros n) =(\varphi\psi\1\<cros(m\li\psi\2)n)
\end{equation}
containing again $\M$ and $\hat\G$ as unital subalgebras. If there are several
coactions under consideration we will also write $\M\>cros\!_\rho\, \hat\G$ and
$\hat\G\<cros\!_\l\, \M$, respectively.
We note that Eq. \no{1.8} implies that as an algebra
\begin{equation}\label{1.10}
\M\>cros\hat\G =\M\hat\G =\hat\G\M
\end{equation}
where we have identified $\M\equiv \M\cros\hat\one$ and
$\hat\G\equiv \one_\M\cros\hat\G$.
In fact using the antipode axioms one easily verifies from \no{1.8}
\begin{equation}\label{1.11}
m\cros\varphi=(m\cros\hat\one)(\one_\M\cros\varphi)
=(\one_\M\cros\varphi\2)((\hat S^{-1}(\varphi\1) \re m)\cros\hat\one)
\end{equation}
Similar statements hold in $\hat\G\<cros \M$. More generally we have

\Lemma{1.1}
{Let $\re :\hat\G\o\M\to \M$ be a left Hopf module action and let $\A$ be an
algebra containing $\M$ and $\hat\G$ as unital subalgebras. Then in $\A$ the
relation
\begin{equation}\label{1.12}
\varphi m =(\varphi\1\re m)\varphi\2,\quad
\,\forall \varphi\in\hat\G, \,\forall m\in\M
\end{equation}
is equivalent to
\begin{equation}\label{1.13}
m\varphi=\varphi\2 (\hat S^{-1}(\varphi\1)\re m),
\quad\,\forall \varphi\in\hat\G,\,\forall m
\in \M
\end{equation}
and if these hold then $\M\hat \G=\hat\G \M \subset\A$ is a subalgebra and
\begin{equation}\label{1.14}
\M\>cros\hat\G  \ni (m\cros\varphi) \mapsto m\varphi\in \M\hat\G,
\end{equation}
is an algebra epimorphism.}

The proof of \lem{1.1} is obvious from the antipode axioms and therefore omitted.
A similar statement of course holds for the crossed product
$\hat \G \<cros\M$. As an application we remark that any left crossed product
$\hat\G\<cros\!_\l\, \M$ can be identified with an associated right
crossed product
$\M \>cros\!_\rho\, \hat\G^{cop}$,
where $\hat \G^{cop} \equiv \widehat{(\G_{op})}$
and where $\rho:\M\to \M\o \G_{op} $ is the right coaction given by
\begin{equation}\label{1.15}
\rho= (\id \o S^{-1})\circ \tau_{\G,\M} \circ \lambda
\end{equation}
$\tau_{\G,\M} :\G\o \M\to \M\o\G$ being the permutation of tensor factors.
In fact we have

\Lemma{1.2}
{Let $\lambda:\M\to\G\o \M$ and $\rho:\M\to \M\o \G_{op}$ be a pair of left and
right coactions, respectively, related by \no{1.15}. Then
\begin{equation}\label{1.16}
f(m\cros \varphi):= (\hat\one \<cros m)(\varphi \<cros\one_\M)
\end{equation}
defines an algebra isomorphism
$f:\M \>cros\!_\rho\, \hat\G^{cop} \to \hat \G \<cros\!_\l\, \M$
with inverse given by}
\begin{equation}\label{1.17}
f^{-1} (\varphi\<cros m) =(\one_\M\cros \varphi)(m\cros \hat\one)
\end{equation}
\proof Let $\re :\hat\G^{cop} \o\M\to\M$ and $\li :\M\o\hat\G \to \M$
be the left and right actions dual to $\rho$ and $\lambda$, respectively. Then
\no{1.15} is equivalent to
$
\varphi \re m=m\li \hat S^{-1} (\varphi)
$
for all $\varphi\in\hat\G$ and $m\in\M$. Putting $\A=\hat\G \<cros\!_\l\, \M$
in \lem{1.1} we have from \no{1.9}
$
m\varphi =\varphi\1(\hat S(\varphi\2)\re m)
$
which is the relation \no{1.13} with $\hat\G$ replaced by
$\hat\G^{cop}$. Hence,
by \lem{1.1} $f$ is an algebra map whose inverse is easily verified to coincide
with \no{1.17}. \qed

\bsn
We conclude this introductory part by describing crossed products in
terms of the
``generating matrix'' formalism as advocated by the St. Petersburg school. Our
presentation will closely follow the review of [N1].
First we note that since $\G$ is finite dimensional we may identify
$\Hom_\CC(\hat\G,V)\cong\G\o V$ for any
$\CC$-vector space $V$. In particular, the relation
\begin{equation}\label{1.19a}
T(\varphi)=(\varphi\o \id )(\TT),\quad \,\forall \varphi\in\hat G,
\end{equation}
provides a one-to-one correspondence
between algebra maps $T:\hat\G\to\A$ into some target algebra $\A$ and
elements $\TT\in\G\o\A$ satisfying
\begin{eqnarray}
\TT^{13} \TT^{23} =(\Delta\o \id )(\TT)\label{1.19}
\end{eqnarray}
where Eq. \no{1.19} is to be understood as an identity in
$\G\o\G\o\A$, the upper
indices indicating the canonical embedding of tensor factors (e.g. $\TT^{23}=
\one_\G\o \TT$, etc.). 
Throughout, we will call elements $\TT\in\G\o\A$ {\em normal}, 
if
$$
(\e \o \id )(\TT) =\one_\A 
$$
which in \Eq{1.19a} is equivalent to $T:\hat\G\to\A$ being unit 
preserving. In what follows, the target algebra $\A$ may always be
arbitrary. In the particular case $\A=\End V$ we would be
talking of representations of $\hat\G$ on $V$, or
more generally, as discussed in \lem{1.4} below, of representations
of $\M\>cros \hat\G$ or $\hat\G\<cros \M$, respectively, on $V$.

\Definition{1.3}
{Let $\lambda :\M\to \G\o\M$ be a left coaction and let $\gam:\M\to\A$ be an
algebra map. An {\em implementer} of $\lambda$ in $\A$ (with respect
to $\gam$) is an element $\LL\in\G\o\A$
satisfying 
\begin{equation}\label{1.21}
\,[\one_\G \o \gam(m)]\,\LL=\LL\,[(\id _\G\o \gam)(\lambda(m))]\,
\end{equation}
for all $m\in\M$. Similarly, an implementer in $\A$
of a right coaction $\rho:\M\to\M\o\G$ is an element
$\RR\in\G\o\A$ satisfying 
\begin{equation}\label{1.22}
\RR\,[\one_\G\o \gam(m)] =[(\id \o \gam)(\rho^{op} (m))]\, \RR
\end{equation}
where $\rho^{op} =\tau_{\M\o\G} \circ \rho$.}

We remark that Eq. \no{1.22} suggests that for right coactions one
should rather
use the convention $\RR\in\A\o\G$ in place of introducing $\rho^{op}$ as
a left coaction of $\G^{cop}$. However, for later purposes it will be more
convenient to place the ``auxiliary copy'' of $\G$ in the definition
of implementers always to the left of $\M$. We now have

\Lemma {1.4}
{Let $\lambda:\M\to\G\o\M$ be a left coaction
and let $\gam:\M\to\A$ be an algebra map. Then the relation
\begin{equation}\label{1.23}
\gam_L(\varphi\<cros m):=(\varphi\o \id )(\LL)\gam(m)
\end{equation}
provides a one-to-one correspondence between algebra maps
$\gam_L:\hat\G\<cros\!_\l\, \M\to\A$ extending $\gam$ and 
normal implementers $\LL\in\G\o\A$
of $\lambda$ satisfying
\begin{equation}\label{1.24}
\LL^{13} \LL^{23} =(\Delta\o \id )(\LL)
\end{equation}
Similarly, if $\rho :\M\to\M\o\G$ is a right coaction, then the relation
\begin{equation}\label{1.25}
\gam_R(m\cros\varphi)=\gam(m)(\varphi\o \id )(\RR)
\end{equation}
provides a one-to-one correspondence between algebra maps 
$\gam_R:\M \>cros \hat\G\to \A$ extending $\gam$ and normal 
implementers $\RR\in\G\o\A$ of $\rho$ satisfying
\begin{equation}\label{1.26}
\RR^{13} \RR^{23} =(\Delta \o \id )(\RR)
\end{equation}
}
\proof Writing $R(\varphi):=(\varphi\o \id )(\RR)\equiv
\gam_R(\one_\M\cros\varphi)\in\A$ and using $\Hom_\CC(\hat\G,\A)\cong\G\o\A$
the relation $\RR\leftrightarrow \gam_R$ is one-to-one. The implementer
property \no{1.22} is then equivalent to
$R(\varphi) \gam(m)=\gam(\varphi\1\re m)R(\varphi\2)$ and $\RR$ 
is normal iff $\gam_\R$ is unit preserving.
Together with the remarks \no{1.19a} - \no{1.19} this is further
equivalent to $\gam_R$
defining an algebra map, similarly as in \lem{1.1}.
The argument for $\gam_L$ is analogous.
\qed

\bsn
Next, we note that the equivalence
\no{1.12} $\Leftrightarrow$ \no{1.13} can be
reformulated for  implementers as follows

\Lemma{1.5}
{Under the conditions of Lemma 1.4 denote $\lambda(m)=m_{(-1)} \o m_{(0)}$
and $\rho(m)=m\0\o m\1$ and let $\TT\in\G \o\A$ be normal.
Dropping the symbol $\gam$ we then have
\\
i) $\TT$ is an implementer of $\rho$ if and only if
\begin{equation}\label{1.22'}
\,[\one_\G \o m]\, \TT=[S^{-1}(m\1)\o \one_\A]\,\TT\,[\one_\G\o m\0],
\quad\,\forall m\in\M
\end{equation}
ii) $\TT$ is an implementer of $\lambda$ if and only if
\begin{equation}\label{1.21'}
\TT\,[\one_\G\o m] =[\one_\G\o m_{(0)}]\, \TT\,[S^{-1} (m_{(-1)}) \o\one_\A],
\quad\,\forall m\in\M
\end{equation}
}
\proof Suppose $\TT$ is an implementer of $\rho$. Then by \no{1.22}
$$
[S^{-1}(m\1)\o\one_\A]\,\TT\,[\one_\G\o m\0] =
 [S^{-1}(m\2)m\1\o m\0]\,\TT = [\one_\G\o m]\,\TT
$$
by \no{1.2} and the antipode axioms. Conversely, if $\TT$ satisfies \no{1.22'}
then
$$
\TT\,[\one_\G\o m] = [m\2 S^{-1} (m\1)\o\one_\A]\, \TT\,[\one_\G\o m\0]
= [m\1 \o m\0]\,\TT
$$
proving \no{1.22}. Part ii) is proven analogously. \qed

\sec{Two-sided coactions and diagonal crossed products}

In Part II we will give a straightforward
generalization of the notion of
coactions to quasi-Hopf algebras. However, in general an associated
notion of a crossed product extension $\M\>cros\hat\G$ will not
be well defined as an associative algebra,
basically because in the quasi-Hopf case the natural product in
$\hat \G$ is not associative.
We are now going to provide a new construction of what we call a
{\em diagonal crossed product}
which will allow to escape this obstruction when generalized to
the quasi-Hopf case. Our diagonal crossed products are always 
based on {\em two-sided coactions} or, equivalently, on pairs 
of commuting left and right coactions. These structures are 
largely motivated by the specific example
$\M=\G$ where our methods reproduce the quantum
double $\D(\G)$.

\Definition{1.5}
{A two-sided coaction of $\G$ on an algebra $\M$ is an algebra map
$\delta:\M\to \G\o\M\o \G$ satisfying
\bea
(\id _\G\o\delta\o \id _\G)\circ \delta
&=& (\Delta \o \id _\M\o\Delta)\circ\delta  \label{1.27}\\
(\e\o \id _\M\o\e) \circ \delta &=& \id _\M \label{1.28}
\eea}
An example of a two-sided coaction is given by
$\M=\G$ and $\delta:=D\equiv (\Delta\o \id )\circ \Delta$.
More generally let $\lambda:\M\to\G\o\M$ and $\rho:\M\to \M\o\G$
be a left and a right coaction, respectively.
We say that $\lambda$ and $\rho$ {\em commute}, if
\begin{equation}\label{1.29a}
(\lambda\o \id ) \circ\rho = (\id \o\rho) \circ\lambda
\end{equation}
It is straight forward to check that in this case
\begin{equation}\label{1.29b}
\delta:=(\lambda\o \id ) \circ\rho\equiv (\id \o\rho)\circ\lambda
\end{equation}
provides a two-sided coaction. Conversely, we have

\Lemma{1.6}
{Let $\delta:\M\to\G\o \M\o\G$ be a two-sided coaction and define $\lambda:=
(\id \o \id  \o\e)\circ\delta$ and $\rho:=(\e\o \id  \o \id )\circ \delta$.
Then $\lambda$ and $\rho$ provide a pair of commuting left and right
coactions, respectively, obeying Eq. \no{1.29b}
}
The proof \lem{1.6} is again straightforward and therefore omitted. Similarly
as in Eqs. \no{1.4a} - \no{1.4d} it will be useful to introduce the
following notations
\begin{equation}\label{1.29c}
\ba{rclll}
\delta(m) &=& m_{(-1)}\o m\0 \o m\1 &&
\\
(\id \o \id \o\Delta)(\delta(m)) &\equiv& (\id \o\rho\o \id )(\delta(m))
&=& m_{(-1)}\o m\0\o m\1\o m\2
\\
(\Delta\o \id \o \id )(\delta(m)) &\equiv& (\id \o\lambda\o \id )(\delta(m)) &=&
m_{(-2)} \o m_{(-1)} \o m\0 \o m\1
\\
(\Delta\o \id \o\Delta)(\delta(m)) &\equiv& (\id \o\delta\o \id )(\delta(m)) &=&
m_{(-2)} \o ....\o m\2
\ea
\end{equation}
etc., implying again the usual summation conventions.
We remark that in the quasi-coassociative setting of Part II the
relation between two-sided coactions and pairs $(\l,\rho)$ of left and
right coactions becomes more involved, justifying the treatment of
two-sided coactions as distinguished objects on their own right also
in the present setting.
Two-sided coactions may of course be considered
as ordinary (right or left)
coactions by rewriting them equivalently as 
\begin{align}\label{1.30}
\rho_\delta &:= \tau_{\G,(\M\o\G)} \circ \delta :\M\to \M\o(\G\o\G^{cop} )\\
\label{1.31}
\la_\delta &:= \tau_{(\G\o\M),\G} \circ\delta:\M\to (\G^{cop}\o \G)\o\M
\end{align}
where the Hopf algebra structures on $\G\o\G^{cop}$ and $\G^{cop} \o\G$ are
the usual tensor product operations induced from the structures on
$\G$ and $\G^{cop}$,
respectively. Conversely, any left coaction $\lambda$ (right coaction
$\rho$) may be considered as a two-sided coaction by putting $\delta(m)=
\lambda(m)\o\one_\G$ or $\delta(m)=\one_\G\o\rho(m)$, respectively.

Next, in view of \lem{1.6} we also have a one-to-one correspondence between
two-sided coactions $\delta$ of $\G$ on $\M$ and pairs of
mutually commuting left and right Hopf module actions, $\re$ and $\li$, of
$\hat\G$ on $\M$, the relation
being given by
\begin{equation}\label{1.32a}
(\varphi\o \id \o\psi)(\delta (m)) =\psi\re m\li \varphi,
\end{equation}
where $\varphi,\psi\in\hat\G$ and $m\in\M$. This allows to construct
as a new algebra the
{\em right diagonal crossed product } $\M\reli\hat\G$ as follows.

\Proposition{1.7}
{Let $\delta=(\l\o\id _\G)\circ\rho=(\id_\G\o\rho)\circ\l$ 
be a two-sided coaction of $\G$ on $\M$ and let $\re$
and $\li$ be
the associated commuting pair of left and right actions of $\hat\G$ on $\M$.
Define on $\M\o\hat\G$ the product
\begin{equation}\label{1.32}
(m\reli \varphi)(n\reli\psi):= (m(\varphi\1\re n\li
\hat S^{-1}(\varphi\3)) \reli\varphi\2\psi)
\end{equation}
where we write $(m\reli\varphi)$ in place of
$(m\o\varphi)$ to distinguish the new algebraic structure.
Then with this product $\M\o\hat\G$ becomes an associative algebra with unit
$(\one_\M\reli\hat\one)$  containing
$\M\equiv\M\reli\hat\one$ and
$\hat\G\equiv\one_\M\reli\hat\G$ as unital subalgebras.}

\proof For $m,m',n\in\M$ and $\varphi,\psi,\xi\in\hat\G$ we compute
\begin{eqnarray*}
&&\big[(m\reli\varphi)(m'\reli\psi)\big] (n\reli\xi)=
\big[m(\varphi\1\re m'\li \hat S^{-1}(\varphi\3))\reli\varphi\2\psi\big]
(n\reli\xi)
\\&&
=\big[m(\varphi\1\re m'\li \hat S^{-1}(\varphi\5))(\varphi\2\psi\1 \re
n\li \hat S^{-1}(\psi\3)\hat S^{-1}(\varphi\4))\big]\reli(\varphi\3\psi\2\xi)
\\&&
=m\big[\varphi\1\re[m'(\psi\1\re n\li \hat S^{-1} (\psi\3))]
\li \hat S^{-1} (\varphi\3)\big]\reli(\varphi\2\psi\2\xi)
\\&&
=(m\reli\varphi)\big[(m'\reli\psi)(n\reli\xi)\big].
\end{eqnarray*}
which proves the associativity.
The remaining statements follow trivially from
$\varphi\re\one_\M =\one_\M\li\varphi=\hat \e(\varphi)\one_\M$
and the counit axioms.
\qed

\bsn
We emphasize that while \prop{1.7} still is almost
trivial as it stands, its true power only appears when generalized to
quasi-Hopf algebras $\G$, which will be done in Part II.

\Definition{1.8}
{Under the setting of \prop{1.7} we define the {\em right diagonal
crossed product} $\M\relid\,\hat\G\equiv\M\relilr\hG$ to be the vector
space $\M\o\hat\G$ with associative multiplication structure \no{1.32}.}

\noindent
Left diagonal crossed products will be constructed in \cor{1.11} below.
In cases where the two-sided coaction $\delta$ is unambiguously
understood from the context we will also write $\M\reli\hat\G$.
We emphasize already at this place that in Part II not every
two-sided coaction will be given as $\delta=(\l\o\id_\G)\circ\rho$ (or
$\delta=(\id_\G\o\rho)\circ\l$ ), in which case the notations
$\M\relid\,\hat\G$ and $\M\relilr\hG$ will denote different (although
still equivalent) extensions of $\M$. Here we freely use either one of them.
If $\delta=\one_\G\o\rho$ or $\delta=\lambda \o \one_\G$ for
some left coaction $\lambda$ (right coaction $\rho$) then
$\M\relid\hat\G=\M\crosr \hat\G$ or $\M\relid
\hat\G=\hat\G\crosl \M$, respectively, where
in the later case we have to invoke the isomorphism \no{1.17}.
More generally, $\M\relilr \hat\G$ may be identified as a
subalgebra of
$\hat\G \crosl(\M\crosr\hat\G)\equiv(\hat\G\crosl \M)\crosr\hat\G$
using the injective algebra map
\begin{eqnarray}
\M\relilr\hat\G \ni(m\reli\varphi)
&\mapsto&
(\hat\one\<cros m\>cros\hat\one)
(\varphi\2\<cros \one_\M \>cros \varphi\1)\nonumber
\\
&\equiv&
\big[\varphi\2\<cros (m\li \varphi\3)\>cros \varphi\1\big]
\in\hat\G\crosl\M\crosr\hat\G\label{1.42a}
\end{eqnarray}
which we leave to the reader to check.
\Eq{1.42a} also motivates our choice of calling the crossed 
product $\M\reli\hat\G$ ``diagonal''.

In the case $\M=\G$ and
$\delta:=D\equiv(\Delta \o \id )\circ \Delta$
the formula \no{1.32} coincides with the multiplication rule in the quantum
double $\D(\G)$ [Dr1,M3], i.e.
\begin{equation}\label{1.32'}
\D(\G) =\G\reliD \hat\G
\end{equation}
It is well known, that $\D(\G)$ is itself again a Hopf algebra
with coproduct $\Delta_D$ given by
\begin{equation}\label{1.32b}
\Delta_D(a\reliD\varphi)=(a\1\reliD\varphi\2)\o(a\2\reliD \varphi\1)
\end{equation}
where $a\in\G$ and $\varphi\in\hat\G$.
It turns out that this result generalizes to diagonal crossed
products as follows

\Proposition{1.9}
{Let $\delta:\M\to\G\o\M\o\G$ be a two-sided coaction. Then
$\M\reli \dG$ admits a commuting pair of coactions
$\lambda_D:\M\reli\dG\to\D(\G)\o(\M\reli\dG)$ and \\
$\rho_D:\M\reli\dG\to(\M\reli\dG)\o\D(\G)$ given by
\begin{eqnarray}\label{1.32c}
\lambda_D(m\reli\varphi)=(m\m1\reliD\varphi\2)\o(m\0\reli\varphi\1)\\
\label{1.32c'}
\rho_D(m\reli\varphi)=(m\0\reliD\varphi\2)\o(m\1\reli\varphi\1)
\end{eqnarray}
where elements in
$\D(\G)$ are written as $(a\reliD\varphi),\ a\in\G,\
\varphi\in\hat\G$, and where we have used the notation
\no{1.4a}-\no{1.4d}.
}

\proof In view of \no{1.32b} the comodule axioms and the 
commutativity \no{1.29b} are obvious. We are left to prove that 
$\lambda_D$ and $\rho_D$
provide algebra maps. To this end we use the following
identities obviously holding for any
two-sided coaction and all $m\in\M,\varphi,\psi\in\hat\G$
\begin{eqnarray}\label{1.32d}
\rho(\varphi\re m\li\psi) &=& m\0\li\psi\o(\varphi\arr m\1)
\\
\label{1.32e}
m\0 \o(m\1\arl \psi) &=&(\psi\re m\0) \o m\1
\end{eqnarray}
With this we now compute
\begin{eqnarray*}
&&\rho_D(m\reli \varphi)\rho_D(n\reli\psi)=\\
&&=\big[m\0(\varphi\4\re n\0\li\hat S^{-1}(\varphi\6))
\reli\varphi\5\psi\2\big] \o\big[m\1(\varphi\1\arr n\1\arl
\hat S^{-1}(\varphi\3))\reliD\varphi\2\psi\1\big]
\\
&&=\big[m\0 (n\0\li\hat S^{-1}(\varphi\4))\reli\varphi\3\psi\2\big]
\o \big[m\1(\varphi\1\arr n\1)\reliD\varphi\2\psi\1\big]
\\
&&=\big[m\0(\varphi\1\re n\li \hat S^{-1} (\varphi\4))\0
\reli \varphi\3\psi\2\big]
\o \big[m\1 (\varphi\1\re n\li \hat S^{-1} (\varphi\4))\1
\reliD \varphi\2\psi\1\big]
\\
&&= \rho_D \big[m(\varphi\1\re n\li \hat S^{-1} (\varphi\3))
\reli\varphi\2\psi\big]
\end{eqnarray*}
where we have used \no{1.32} in the first equation, \no{1.32e}
and the antipode
axioms in the second equation and \no{1.32d} in the third equation. Hence
$\rho_D$ provides an algebra map. The argument for $\lambda_D$ 
is analogous.
\qed

\bsn
Before presenting further examples of diagonal crossed products
let us recall the well known Hopf algebra identity
$\D(\hat\G)= \D(\G)^{cop}$ where the algebra isomorphism is 
given by
\begin{equation}\label{1.32f}
\D(\hat\G)\ni (\varphi\reli\!\!_{\hat D}\,\, a)\mapsto
(\one \reliD\varphi)(a\reliD\hat\one)\in \D(\G)
\end{equation}
Again this generalizes to diagonal crossed products in the sense that they
may equivalently be modeled on the vector space $\hat\G\o\M$.
To see this, we first note that there is an analogue of \lem{1.1}

\Lemma{1.10}
{Let $\A$ be an algebra containing $\M$ and $\hat\G$ as unital
subalgebras and let $\re$ and $\li$ be a commuting pair of left
and right Hopf module actions of
$\hat\G$ on $\M$. Then in $\A$ the relations
\begin{equation}\label{1.33a}
\varphi m=[\varphi\1 \re m\li \hat S^{-1} (\varphi\3)] \varphi\2,\quad
\forall \varphi\in\hat\G,\, \forall m\in\M
\end{equation}
are equivalent to
\begin{equation}\label{1.33b}
m\varphi = \varphi\2 [\hat S^{-1} (\varphi\1)\re m \li \varphi\3],
\quad\forall\varphi\in\hat\G,\,\forall m\in\M
\end{equation}
and if these hold then $\M\hat\G=\hat\G\M \subset\A$ is a subalgebra and
$
\M\reli\hat\G\ni(m\reli \varphi)\mapsto m\varphi\in\M\hat\G
$
is an algebra epimorphism.}

\bsn
The proof of \lem{1.10} is again  straightforward from the antipode axioms
and is left to the reader. Similarly as in \lem{1.2} we now conclude

\Corollary{1.11}
{Under the setting of \prop{1.7} define 
$\hat\G\reli\M\equiv\hat\G\relid\M\equiv\hat\G\relilr\M$
to be the vector space $\hat\G\o\M$ with multiplication rule
\begin{equation}\label{1.34}
(\varphi\reli m)(\psi\reli n):=
\varphi\psi\2\reli (\hat S^{-1} (\psi\1)\re m\li \psi\3)n
\end{equation}
Then $\hat\G\reli \M$ is an associative algebra and
\begin{equation}\label{1.35}
\hat\G\reli \M\ni \varphi\reli m\mapsto
(\one_\M \reli\varphi)(m\reli\hat\one)\equiv
\varphi\1\re m\li\hat S^{-1}(\varphi\3)\reli\varphi\2 \in\M\reli\hat\G
\end{equation}
provides an isomorphism restricting to the identity on $\M$
with inverse given by
\begin{equation}\label{1.36}
\M\reli\hat\G\ni m\reli\varphi\mapsto
(\hat\one\reli m)(\varphi\reli\one_\M)\equiv
\varphi\2\reli\Big(\hat S^{-1}(\varphi\1)\re m\li\varphi\3\Big)
\in\hat\G\reli\M
\end{equation}}
\bsn
Here we have used the same symbol $\reli$ on either side in order not to
overload the notation.
The reader is invited to check that in the case $\M=\G$ and
$\delta=(\Delta\o \id ) \circ \Delta$ we recover
$\hat\G\relid \M=\D(\hat\G)$.
Also note that after a trivial permutation of tensor factors
the multiplication rule \no{1.34} implies
\begin{equation}\label{1.34'}
(\hat\G\relid\M)_{op} = \M_{op}\relidop\hat\G_{op}^{cop}\ ,
\end{equation}
where ``op'' refers to the algebra with opposite multiplication
and where $\delta_{op}:=\delta^{321}$.
We propose to call $\hat\G\reli\M$ the {\em left diagonal} and 
$\M\reli\hat\G$ the {\em right diagonal} crossed product.


\sec{Examples and applications}
\subsection{Double crossed products}

In this Subsect. we relate our diagonal crossed product with
the double crossed product construction of Majid [M3,4]. Here
we adopt the version [M3, Thm.7.2.3], according to which a
bialgebra $\B$ is a double crossed product, written as
\begin{equation}
\label{4.1.1}
\B=\M\doubcross\H,
\end{equation}
iff $\M$ and $\H$ are sub--bialgebras of $\B$ such that the
multiplication map
$$
\mu:\M\o\H\ni m\o h\mapsto mh\in\B
$$
provides an isomorphism of coalgebras.
In this case the bialgebras $\M$ and $\H$ become a
{\em matched pair} with mutual actions
$
\Re:\H\o\M\to\M
$
and
$
\Li:\H\o\M\to\H
$
given by
\begin{eqnarray}
\label{4.1.2}
h\Re m &:=& (\id_\M\o\e_\H)(\mu^{-1}(hm))\\
\label{4.1.3}
h\Li m &:=& (\e_\M\o\id_\H)(\mu^{-1}(hm)),
\end{eqnarray}
see [M3, Chap.7.2] for more details.
The guiding example is again given by the quantum double
satisfying
\begin{equation}
\label{4.1.4}
\D(\G)=\G\doubcross\hat\G^{cop}.
\end{equation}
More generally, any diagonal crossed product $\M\relilr\hG$
becomes a double crossed product
\begin{equation}
\M\relilr\hG = \M\doubcross\hG^{cop}
\end{equation}
provided $\M$ is equipped with a bialgebra structure
$\Delta_\M,\ \e_\M$ such that the tensor product coalgebra
structures
\begin{eqnarray}
\label{4.1.6}
\Delta_\B(m\reli\psi) &:=& (m^{(1)}\reli\psi\2)\o(m^{(2)}\reli\psi\1)\\
\label{4.1.7}
\e_\B(m\reli\psi) &:=& \e_\M(m)\hat\e(\psi)
\end{eqnarray}
give algebra maps $\Delta_\B:\B\to\B\o\B$ and $\e_\B:\B\to\CC$
(here we denote $\Delta(m)\equiv m^{(1)}\o m^{(2)}$).
In this way all {\em generalized quantum doubles}
$\M\doubcross\hG^{cop}$ in the sense of [M3, Ex.7.2.7] are
as algebras actually diagonal crossed products in our sense.
These may be described in terms of the bialgebra homomorphism
$\kappa:\M\to\G$ induced by the Hopf skew-pairing
$\M\o\hG^{cop}\to\CC$ required in [M3, Chap. 7.2], which gives
rise to the commuting pair of $\G$-coactions
\begin{equation}
\label{4.1.8}
\l:=(\kappa\o\id_\M)\circ\Delta_\M\qquad ,
\qquad\rho:=(\id_\M\o\kappa)\circ\Delta_\M.
\end{equation}
For these models the ``matched pair'' actions \no{4.1.2} and
\no{4.1.3} are given by the coadjoint actions
\begin{eqnarray}
\label{4.1.9}
\psi\Re m &=& \psi\1\re m\li\hat S^{-1}(\psi\2)\\
\psi\Li m &=& \e_\M(\psi\1\re m\li\hat S^{-1}(\psi\3))\,\psi\2
\nonumber\\
\label{4.1.10}
&=&\hat S^{-1}(\kappa(m^{(1)}))\arr\psi\arl\kappa(m^{(2)}).
\end{eqnarray}
Also note that the coactions \no{4.1.8} satisfy the
compatibility condition
\begin{equation}
\label{4.1.11}
(\rho\o\id_\M)\circ\Delta_\M=(\id_\M\o\l)\circ\Delta_\M
\end{equation}
 and the homomorphism $\kappa:\M\to\G$ may be recovered as
\begin{equation}
\label{4.1.12}
\kappa = (\id_\G\o\e_\M)\circ\l=(\e_\M\o\id_\G)\circ\rho.
\end{equation}
Conversely, we have
\Proposition{4.0}
{Let $(\M,\,\Delta_\M,\,\e_\M)$ be a bialgebra and let
$(\l,\rho)$ be a commuting pair of (left and right)
$\G$-coactions on $\M$. Then the compatibility condition
\no{4.1.11} is equivalent to
\begin{equation}
\label{4.1.13}
(\id_\G\o\e_\M)\circ\l = (\e_\M\o\id_\G)\circ\rho =:\kappa
\end{equation}
being a bialgebra homomorphism $\kappa:\M\to\G$ which satisfies
\Eq{4.1.8}.
If these conditions are satisfied, then the diagonal crossed
product $\M\relilr\hG$ is also a double crossed product with
respect to the bialgebra structure \no{4.1.6}, \no{4.1.7}, i.e.
in this case we have
$$
\M\relilr\hG= \M\doubcross\hG^{cop}.
$$
}
\proof
Applying $(\e_\M\o\id_\G\o\e_\M)$ to \no{4.1.11} proves that the
two expressions in \no{4.1.13} define the same algebra map
$\kappa:=\M\to\G$.
Clearly, we have $\e_\G\circ\kappa=\e_\M$. To prove \no{4.1.8}
we compute
\begin{eqnarray*}
(\kappa\o\id_\M)\circ\Delta_\M &=&
(\id_\G\o\e_\M\o\id_\M)\circ(\l\o\id_\M)\circ\Delta_\M
\\
&=&(\e_\M\o\id_\G\o\e_\M\o\id_\M)\circ(\id_\M\o\l\o\id_\M)\circ\DeltaM^{(2)}
\\
&=&(\e_\M\o\id_\G\o\e_\M\o\id_\M)\circ(\rho\o\idM\o\idM)\circ\DeltaM^{(2)}
\\
&=&(\e_\M\o\id_\G\o\idM)\circ(\rho\o\idM)\circ\DeltaM
\\
&=&(\e_\M\o\l)\circ\DeltaM =\l,
\end{eqnarray*}
where
$\DeltaM^{(2)}=(\DeltaM\o\id)\circ\DeltaM=(\id\o\DeltaM)\circ\DeltaM$
and where we have used \no{4.1.11} in the third and the fifth
line. The identity $\rho=(\id\o\kappa)\circ\DeltaM$ is proven
similarly. Finally, $\kappa:\M\to\G$ is a bialgebra map, since
\begin{eqnarray*}
(\kappa\o\kappa)\circ\DeltaM &=& (\id_\G\o\kappa)\circ\l
=(\id_\G\o\id_\G\o\e_\M)\circ(\id_\G\o\l)\circ\l
\\
&=&(\Delta_\G\o\e_\M)\circ\l=\Delta_\G\circ\kappa.
\end{eqnarray*}
by the $\G$-coaction property for $\l$.
According to the results of [M3, Sect. 7.2]  this implies that
the ``generalized quantum double'' $\M\doubcross\hG^{cop}$
as an algebra coincides with the diagonal crossed product
$\M\relilr\hG$ in our sense.
\qed

\bsn
Using Theorem I we will give a short direct proof in Sect. 5
that under the compatibility conditions \no{4.1.11} the
colgebra structures $\Delta_\B$ \no{4.1.6} and $\e_\B$
\no{4.1.7} are indeed algebra maps.

We also remark that it seems natural to conjecture that the
condition \no{4.1.11} is also necessary for $\M\relilr\hG$
being a double crossed product with respect to this coalgebra
structure. However we have not succeeded in proving this.


\subsection{Two--sided crossed products}
We now provide further examples of diagonal crossed products. 
A simple recipe
to produce two-sided $\G$-comodule algebras $(\M, \delta)$ is by taking a
right $\G$-comodule algebra $(\A,\rho)$ and a left $\G$-comodule algebra
$(\B,\lambda)$ and define $\M=\A\o\B$ and
\begin{equation}\label{1.36a}
\delta(A\o B):=B\m1 \o(A\0\o B\0)\o A\1
\end{equation}
where $A\in\A,\ B\subset\B,\ \rho(A)=A\0\o A\1$ and $\lambda(B)=B\m1 \o B\0$.
In terms of the $\hat\G$-actions $\re$ on $\A$ and $\li$ on $\B$ dual
to $\rho$ and $\lambda$, respectively, the $\hat\G$-actions $\re_\M$ and
$\li_\M$ dual to \no{1.36a} are given by
\begin{equation}\label{1.36b}
\varphi \re\!_\M\,(A\o B)\li\!_\M\,\psi =(\varphi \re A\o B\li \psi)
\end{equation}
where $\varphi,\psi\in\hat\G$.
Hence, we may construct the diagonal crossed product
$\M\reli\G$ as before. It turns out that this example may be
presented differently as a so-called {\em two-sided crossed product}.

\Proposition{1.12}
{Let $\re:\hat\G\o\A\to\A$ and $\li:\B\o\hat\G\to\B$ be a left and
a right Hopf module action, respectively. Define the 
``two-sided crossed product''
$\A\>cros\hat\G\<cros\B$ to be the vector space
$\A\o\hat\G\o\B$ with multiplication structure
\begin{equation}\label{1.37}
\big(A\>cros\phi\<cros B\big)\,\big(A'\>cros\psi\<cros B'\big)=
A(\phi\1\re A')\>cros\phi\2\psi\1\<cros(B\li\psi\2)B'\ .
\end{equation}
Then $\A\>cros\hat\G\<cros\B$ becomes an associative algebra
with unit $\one_\A\>cros\hat\one\<cros\one_\B$ and
\begin{equation}\label{1.38a}
f:\A\>cros\hat\G\<cros\B\ni A\>cros\phi\<cros B\mapsto
\big((A\o\one_\B)\reli\phi\big)\big((\one_\A\o B)\reli\hat\one\big)
\in(\A\o\B)\reli\hat\G
\end{equation}
provides an algebra isomorphism with inverse given by
\begin{equation}\label{1.38b}
f^{-1}\big((A\o B)\reli\phi\big)
=(\one_\A\>cros\hat\one\<crosB)(A\>cros\phi\<cros\one_B)
\end{equation}
}
\proof
To prove associativity of \no{1.37} we compute for
$A,A',A''\in\A,\ B,B',B''\in\B$ and $\varphi,\psi,\xi\in\hat\G$
\begin{eqnarray*}
\lefteqn{\big[ (A\>cros \varphi\<cros B)(A'\>cros \psi\<cros B')\big]
(A'' \>cros\xi \<cros B'')=}\\
&&=\big[ A(\varphi\1 \re A') \>cros \varphi\2 \psi\1 \<cros
(B\li \psi\2)B'\big] (A''\>cros\xi \<cros B'')\\
&&=A(\varphi\1\re A')(\varphi\2\psi\1 \re A'')\>cros \varphi\3\psi\2 \xi\1
\<cros (B\li \psi\3\xi\2)(B'\li \xi\3)B''\\
&&=(A\>cros \varphi\<cros B)\big[ A'(\psi\1\re A'')\>cros\psi\2\xi\1 \<cros
(B'\li\xi\2)B''\big]\\
&&=(A\>cros\varphi \<cros B)
\big[ (A'\>cros \psi\<cros B')(A'' \>cros \xi \<cros  B'')\big]
\end{eqnarray*}
Hence \no{1.37} is associative. Obviously $(\one_\A\>cros \hat\one \<cros
\one_B)$ is the unit in $\A\>cros \hat\G \<cros \B$.
To prove that $f$ is an algebra map we note  $f(A\>cros \varphi \<cros B)=
\big[A(\varphi\1\re \one_\A)\o B\li S^{-1} (\varphi\3)\big]\reli \varphi\2
=\big[A\o B\li S^{-1} (\varphi\2)\big]\reli \varphi\1$ implying
\begin{eqnarray*}
\lefteqn{f(A(\varphi\1 \re A')\>cros \varphi\2 \psi\1 \<cros
(B\li\psi\2)B')=}\\
&&= \big[A(\varphi\1\re A')\o \big( (B\li\varphi\3)B'\big)
\li S^{-1} (\varphi\3\psi\2)\big]\reli \varphi\2\psi\1\\
&&=\big[A(\varphi\1\re A')\o \big(B\li S^{-1} (\varphi\4)\big)
\big(B'\li S^{-1} (\psi\2)S^{-1}(\varphi\3)\big)\reli \varphi\2 \psi\1\\
&&=\big[(A\o B\li S^{-1} (\varphi\2))\reli \varphi\1\big]
\big[(A'\o B'\li S^{-1}(\psi\2))\reli \varphi\1\big]\\
&&=f(A\>cros\varphi \<cros B) f(A'\>cros \psi \<cros B')
\end{eqnarray*}
Hence $f$ is an algebra map. The proof that $f^{-1}$ is given
by \no{1.38b} is left to the reader. \qed

\bsn
As a particular example of the setting of \prop{1.12} we may
choose $\A=\B=\G$ with its canonical left and right $\hat
\G$-action. It turns out that in this case the two-sided
crossed product
$\G\>cros \hat\G \<cros \G \equiv (\G \o\G)\reli \hat \G$
is isomorphic to the iterated crossed product
$(\G \>cros \hat\G)\>cros \G$.
More generally we have

\Proposition{1.13}
{Let $\A$ be a right $\G$-comodule algebra and consider the
iterated crossed product $(\A\>cros \hat \G)\>cros \G$, where
$\G$ acts on $\A\>cros \hat\G$ in the usual way by $a\re
(A\>cros \varphi):= A\>cros (a\arr\varphi),\  A\in\A,\
a\in\G,\ \varphi \in\hat \G$. Then as an algebra
$(\A\>cros \hat\G)\>cros \G=\A\>cros \hat\G \<cros \G$ with
trivial identification. }

\proof The claim follows from
$$A(\varphi\1\re A') \>cros \varphi\2 \psi\1 \<cros (a\arl \psi\2)b=
A(\varphi\1\re A') \>cros \varphi\2 (a\1\arr \psi))\>cros a\2 b $$
as an identity in $\A \o \hat\G\o \G$, where we have used
\begin{equation}\label{1.38c} 
\psi\1 \o (a\arl \psi\2)=\psi\1\bra a\1|\psi\2\ket \o a\2=
(a\1\arr\psi)\o a\2
\end{equation}
as an identity in $\hat\G\o\G$.
\qed

It will be shown in Section 11.2 that being a particular 
example of a two-sided (and therefore of a diagonal) 
crossed product the analogue of $(\A\>cros \hat\G)\>cros \G\equiv \A\>cros
\hat\G\<cros\G$ may also be constructed for quasi-Hopf algebras
$\G$. However, in this case $\A\>cros \hat\G$ (if defined to be
the linear subspace $\A\o\hat\G\o \one_\G$ ) will no longer
be a subalgebra of $\A\>cros \hat\G \<cros \G$.
We will see in Section 11.4 that this fact is very much analogous to what
happens in the field
algebra constructions with quasi-Hopf symmetry as given by V.\ Schomerus [S].
%
%
\subsection{Hopf spin chains and lattice current algebras}
Next, we point out that \prop{1.12} and \prop{1.13} also apply to the
construction of Hopf algebraic quantum chains as considered in
[NSz1, AFFS]\footnote{
For earlier versions of lattice current algebras see also [AFSV, AFS, FG],
the relation with the model of [NSz1] being clarified in [N1].
}.
To see this let us shortly review the model of [NSz1], where one considers even
(odd) integers to represent the sites (links) of a one-dimensional lattice
and where one places a copy of $\G\cong\A_{2i}$ on each site and a copy
of $\hat\G\cong\A_{2i+1}$ on each link. Non-vanishing commutation relations are
then postulated only on neighboring site-link pairs, where one requires [NSz1]
\begin{eqnarray}
\label{1.39} A_{2i} (a) A_{2i-1} (\varphi)
&=& A_{2i-1} (a\1\arr\varphi)A_{2i}(a\2)\\
\label{1.40} A_{2i+1} (\varphi) A_{2i} (a)
&=& A_{2i} (\varphi\1\arr a)A_{2i+1}(\varphi\2)
\end{eqnarray}
Here $\G\ni a\mapsto A_{2i}(a) \in\A_{2i} \subset \A$
and $\hat\G \ni \varphi\mapsto A_{2i+1} (\varphi)\in\A_{2i+1}
\subset \A$ denote the embedding of the
single site (link) algebras into the global quantum chain $\A$. Denoting
$\A_{i,j} \subset \A$ as the subalgebra generated  by $\A_\nu, i\le\nu\le j$,
we clearly have from \no{1.39} and \no{1.40}
\begin{eqnarray}
\label{1.41} \A_{i,j+1} &=& \A_{i,j} \>cros \A_{j+1}\\
\label{1.42} \A_{i-1,j} &=& \A_{i-1} \<cros \A_{i,j}
\end{eqnarray}
Hence, by Proposition 4.2, we recognize the two-sided crossed products
\begin{equation}\label{1.43}
\A_{2i,2j+2} =\A_{2i,2j} \>cros \hat\G \<cros \G
\end{equation}
More generally for all $i\le \nu\le j-1$ we have
\begin{equation}\label{1.44}
\A_{2i,2j} = \A_{2i,2\nu} \>cros \hat\G\<cros \A_{2\nu +2,2j}
\end{equation}
where $\hat\G\equiv \A_{2\nu+1}$.
The advantage of looking at it in this way again comes from the
fact that the constructions \no{1.43} and \no{1.44} generalize
to quasi-Hopf algebras $\G$ whereas \no{1.41} and \no{1.42} do not.
Using the identifications given in Eqs.\ \no{1.3}-\no{1.8} of [N1] a
similar remark applies to the lattice current 
algebra of [AFFS]. This observation will be needed to formulate
a theory of Hopf spin models and lattice current algebras at 
roots of unity, see Section 11.3.

Also note that the
identification \no{1.44} together with \prop{1.12} and
\prop{1.9} immediately imply that quantum
chains of the type \no{1.39}, \no{1.40} admit localized commuting left and right
coactions of the quantum double $\D(\G$), which is precisely the result of
Theorem 4.1 of [NSz1]. In fact, applied to the example in \prop{1.12}
, \prop{1.9} 
gives (for \Eq{1.46} use the identification \no{1.32f})

\Corollary{1.14}
{Under the setting of \prop{1.12} we have two commuting right coactions
$\rho_D:\A\>cros\hat\G\<cros \B\to(\A\>cros \hat \G\<cros\B)\o \D(\G)$ and
$\hat\rho_D:(\A\>cros \hat\G\<cros \B)\o \D(\hat\G)$ given by
\begin{eqnarray}
\label{1.45} 
\rho_D(A\>cros \varphi \<cros B) &=& (A\1 \>cros\varphi\2)
\o (A\2\reliD\varphi\1)\\
\label{1.46} 
\hat\rho_D(A\>cros \varphi\<cros B) &=& (A\>cros\varphi\1\<cros B\0 )
\o (\varphi\2\relihD B\m1)
\end{eqnarray}
}
Next, we remark that the identification \no{1.44} may be
iterated in the obvious way. This observation also
generalizes to the situation where in \prop{1.12} $\A$ and $\B$
are both two-sided $\G$-comodules algebras with
dual $\hat\G$ actions
denoted$\re\!_\A\,,\li\!_A\,,\re\!_B\,,\li\!_B\,$, respectively. Then in the
multiplication rule \no{1.37} only $\re\!_\A\,$ and $\li\!_B\,$
appear and one easily
checks, that for $\varphi,\psi\in\hat\G$ and $A\in\A,\ B\in\B$
the definitions
\begin{eqnarray}
\label{1.44'}
\varphi\re (A\>cros\psi \<cros B) &:=& 
(A\>cros \psi\<cros (\varphi\re\!_\B\, B))\\
\label{1.45'}
(A\>cros\psi \<cros B)\li \varphi &:=& 
((A\li\!_\A\,\varphi) \>cros \psi\<cros B)
\end{eqnarray}
again define a two-sided $\G$-comodule structure on
$\A\<cros \hat\G \>cros \B$.
Hence, we have a multiplication law on two-sided $\G$-comodule
algebras which
is in fact associative, i.e. as a two-sided $\G$-comodule algebra
\begin{equation}\label{1.46'}
(\A\>cros \hat\G\<cros \B)\>cros \hat\G \<cros C = \A\>cros
\hat\G\<cros (B\>cros \hat\G \<cros C)
\end{equation}
which the reader will easily  check. Obviously, one may also consider mixed
cases, e.g. where in \no{1.46'} $\A$ is only a right
$\G$-comodule algebra, but
$\B$ and $C$ are two-sided, in which case \no{1.46'} would be an identity
between right $\G$-comodule algebras.

We conclude this subsection by mentioning that diagonal
crossed products also appear
when formulating periodic boundary conditions for the quantum
chain \no{1.39} - \no{1.40}. In this case, starting with the
open chain $\A_{0,2i}$ localized
on $[0,2i] \cap \ZZ$ one would like to add another copy of
$\hat\G$ sitting on the link $2i+1\equiv -1$ joining the sites
$2i$ and $0$ to form a periodic lattice.
Algebraically this means that $\A_{2i+1} \cong \hat\G$
should have non-vanishing
commutation relations with $\A_{2i} \cong \G$ and $\A_0\cong\G$
in analogy with
\no{1.39} and \no{1.40}, i.e.
\begin{eqnarray}\label{1.47}
A_{2i+1} (\varphi) A_{2i} (a)
&=& A_{2i} (\varphi\1\arr a) A_{2i+1} (\varphi\2)\\
\label{1.48}
A_0(a) A_{2i+1}(\varphi) &=& A_{2i+1} (\varphi\1) A_0 (a\arl\varphi\2)
\end{eqnarray}
Written in this way Eqs.\no{1.47}
and \no{1.48} are precisely the relations in
$$\A_{0,2i} \relid\hat\G$$
where $\delta:\A_{0,2i} \to\G\o\A_{0,2i} \o\G$ is the two--sided
coaction given by
$\delta|_{\A_0}=\Delta \o\one_\G,\ \delta|_{\A_{2i}}=\one_\G\o\Delta$
and $\delta|_{\A_{1,2i-1}}=\one_\G\o \id \o\one_\G$.
Hence, the periodic quantum  chain appears as a
diagonal crossed product of the open
lattice chain by a copy of $\hat\G$ sitting on the link joining the
end points.

\sec{Generating matrices}

Similarly as in \lem{1.4} we now describe the defining relations
of diagonal
crossed products in terms of a generating matrix.
However, whereas in \lem{1.4} the generating matrices $\LL$ and
$\RR$ had to fulfill the {\em implementer}
properties \no{1.21} or \no{1.22}, respectively, the natural
requirement here is that $\TT$ {\em intertwines} the left and
right coactions associated with $\delta$.

\Definition {1.15}
{Let $(\l,\rho)$ be a commuting
pair of left and right $\G$-coactions on $\M$ and let
$\gam:\M\to\A$ be an algebra map into some target algebra $\A$.
Then a {\em $\l\rho$-intertwiner} in $\A$
(with respect to $\gam$) is an element $\TT\in\G\o\A$ satisfying
\begin{equation}\label{1.49}
\TT\lambda_\A(m) =\rho_\A^{op}(m)\TT, \quad\forall m\in\M\ ,
\end{equation}
where $\lambda_\A\equiv(\gam\o \id )\circ \lambda$ and 
$\rho_\A\equiv(\id \o \gam)\circ\rho$.
A $\l\rho$-intertwiner is called {\em coherent} if in 
$\G\o\G\o\A$ it satisfies
\begin{equation}\label{1.51}
\TT^{13} \TT^{23} =(\Delta\o \id )(\TT)
\end{equation}
}
\bsn
Similarly as in \lem{1.5} we then have

\Lemma{1.16}
{Let $(\M,\delta)$ be a two-sided $\G$-comodule algebra with
associated commuting left and right $\G$-coactions $(\l,\, \rho)$,
and let $\gam:\M\to\A$
be an algebra map. Then for $\TT\in\G\o\A$ the following 
properties are equivalent:
\\[.1cm]
i) $\TT$ is a $\l\rho$-intertwiner
\\
ii) $\TT\,[\one_G\o \gam(m)]
=[m\1\o \gam(m\0)]\, \TT[S^{-1} (m\m1)\o\one_\A]$
\\
iii) $[\one_G\o \gam(m)]\,\TT=
[S^{-1} (m\1)\o\one_\A]\,\TT[m\m1 \o \gam(m\0)]$
.}
\bsn
\proof Suppose $\TT$ is a $\l\rho$-intertwiner. Then
$$
[m\1\o \gam(m\0)]\,\TT\,[S^{-1} (m\m1)\o \one_\A]=
\TT\,[m\m1 S^{-1} (m\n2 )\o \gam(m\0)]=\TT\,[\one_\G \o \gam(m)]
$$
by the antipode axiom. Conversely, if $\TT$ satisfies (ii) then
$$
\TT\,[m\m1 \o \gam(m\0)] = [m\1\o \gam(m\0)]\,\TT\,
[S^{-1}(m\m1)m\n2 \o \one_\A]= [m\1\o \gam(m\0)]\, \TT
$$
proving (i) $\Leftrightarrow$ (ii). The equivalence (i)
$\Leftrightarrow$ (iii) follows similarly.
\qed

\bsn
Note that by applying $(\varphi\o \id )$ to both sides the equivalence of (ii)
and (iii) in \lem{1.16} precisely reflects the equivalence of \no{1.33a} and
\no{1.33b} in \lem{1.10}. 

We now arrive at a proof of Theorem I by
concluding similarly as in \lem{1.4}.

\Proposition{1.17}
{Let $(\M,\delta)$ be a two-sided $\G$-comodule algebra with
associated commuting pair of coactions $(\l,\,\rho)$, and let $\gam:\M\to\A$
be an algebra map. Then the relation
\begin{equation}\label{1.50}
\gam_T (m\reli \varphi) =\gam(m) (\varphi\o \id ) (\TT)
\end{equation}
provides a one-to-one correspondence between normal coherent
$\l\rho$-intertwiners $\TT$ and unital algebra maps
$\gam_T:\M\relilr\hat\G\to \A$ extending $\gam$.
}
\proof
Let $T(\varphi):=(\varphi\o \id )(\TT)$. Then \Eq{1.51} 
together with normality is equivalent to
$\hat\G\ni \varphi\mapsto T(\varphi)
\equiv \gam_T(\one_\M\reli \varphi)\in\A$
being a unital algebra morphism and the correspondence
$\TT\leftrightarrow \gam_T|_{\one_\M\reli\hat\G}$ is one-to-one.
Clearly, $\gam_T$ extends $\gam$ 
and \lem{1.16} ii) implies $\forall\varphi\in\hat\G,\ m\in\M$
$$
T(\varphi)\gam(m) = \gam(\varphi\1\re m\li S^{-1} (\varphi\3)) T(\varphi\2)
$$
and therefore $\gam_T$ is an algebra map. Conversely, since
$(m\reli\varphi)=(m\reli\hat\one)(\one_\M\reli\varphi)$, any algebra map
$\gam:\M\relid\hat\G\to\A$ is of the form \no{1.50}.\qed

\bsn
We remark that one could equivalently have chosen to work with
\begin{equation}\label{1.51'}
\gam_T^{op} (\varphi \reli m) :=(\varphi\o \id )(\TT) \gam(m)
\end{equation}
to obtain algebra maps $\gam_T^{op}: \hat \G\relilr\M\to \A$.

Note that \prop{1.17} and \cor{1.11} prove part 1.) and 3.) of Theorem I by
putting $\M_1:=\M\relilr\hat\G,\ \mu_R:=\id_{\M\o\hG}$ and $\mu_L$ as
given in \Eq{1.35}. The uniqueness of $\M_1$ up to equivalence follows
by standard arguments.

Putting $\A =\M_1$ and $\gam_T = \id$ in Proposition \ref{Prop 1.17}
we obtain $\TT =\GG \equiv e_\mu\tp (\one_\M\reli e^\mu) \in
\G\tp\M_1$. Thus we call $\GG$ the {\it universal
  $\la\rho$--intertwiner in $\M_1$}.

Applying the above formalism to the case $\M=\G$ and
$\delta=D\equiv(\Delta\o \id )\circ\Delta$ we realize that
\no{1.49} becomes (suppressing the symbol $\gam$)
\begin{equation}\label{1.52}
\TT\Delta (a) = \Delta_{op} (a) \TT,~~\forall a\in\A
\end{equation}
As already remarked, in this case $\G\reliD\hat\G\equiv \D(\G)$
is the quantum double of $\G$, in which case \prop{1.17}
coincides with [N1, Lem.5.2]
describing $\D(\G)$ as the unique algebra generated by $\G$ and
the entries of
a generating Matrix ${\bf D} \equiv \GG_{\D(\G)} \in\G\o\D(\G)$
satisfying \no{1.51} and \no{1.52}. 
More generally, \prop{1.17} may be reformulated to say
that  in \lem{1.10}
the algebra embeddings $\hat\G\to\A$ satisfying \no{1.33a} and
\no{1.33b} are in one-to-one correspondence with
$\l\rho$-intertwiners $\TT\in\G\o\A$ satisfying \no{1.51}. We
will see in Section 9 that this is the
appropriate language to be generalized to the quasi-Hopf case.

\bigskip
We conclude this Section with demonstrating the good use of 
\prop{1.17} by giving a short alternative proof of how the 
compatibility condition \no{4.1.11} in \prop{4.0} makes 
diagonal crossed products into double crossed product 
bialgebras.

\Proposition{5.4}
{Under the conditions \eqref{4.1.8}, \eqref{4.1.11}, and 
\eqref{4.1.12} of \prop{4.0}
let $\B:=\M\relilr\hG$.
Then the coalgebra structures $\Delta_\B$ \no{4.1.6} and 
$\e_\B$ \no{4.1.7} become algebra maps.}

\proof
Since $\Delta_\B$ extends $\DeltaM$ we may use \prop{1.17} by 
putting $\A:=\B\o\B,\\ \gamma:=\DeltaM:\M\to\M\o\M\subset\A$ and
$$
\TT:=(\id_\G\o\Delta_\B)(\Ga)\equiv\Ga^{13}\Ga^{12}\in\G\o\A,
$$
where $\Ga\in\G\o\B$ is the universal \lri. Thus, 
$\Delta_\B:\B\to\A$ provides an algebra map if and only if 
$\TT$ is a normal coherent \lri\ with respect to 
$\DeltaM:\M\to\A$. Now normality and coherence of $\TT$ hold, 
since the restriction $\Delta_\B |_{\hG}\equiv\hat\Delta_{op}$ 
is a unital algebra map. To prove the \lri\ property for $\TT$ 
we use the identities \no{4.1.8} to compute for all $m\in\M$
\begin{eqnarray*}
\TT(\idG\o\DeltaM)(\l(m))
&=& \Ga^{13}\Ga^{12}[\kappa(m^{(1)})\o m^{(2)}\o m^{(3)}]
\\
&=& \Ga^{13}\Ga^{12}[\l(m^{(1)})\o m^{(2)}]
\\
&=&\Ga^{13}[\kappa(m^{(2)})\o m^{(1)}\o m^{(3)}]\Ga^{12}
\\
&=&[\kappa(m^{(3)}\o m^{(1)})\o m^{(2)}]\Ga^{13}\Ga^{12}
\\
&=&
(\idG\o\DeltaM)(\rho^{op}(m))\TT,
\end{eqnarray*}
where $ m^{(1)}\o  m^{(2)}\o  m^{(3)}\equiv \DeltaM^{(2)}(m)$.
Hence, $\TT$ is a \lri\ and therefore $\Delta_\B$ is an 
algebra map.
Similarly, we prove that $\e_\B:\B\to\CC$ is multiplicative by
putting $\A=\CC$ and $\TT=(\idG\o\e_\B)(\Ga)\equiv\one_\G$.
In this case the \lri\ property reduces to
$$
(\idG\o\e_\M)(\l(m))=(\e_\M\o\idG)(\rho(m)),
$$
which is precisely the condition \no{4.1.13}
\qed

%
%

\renewcommand{\e}{{\bf 1}}
\renewcommand{\eG}{{\e}_{\G}}
\renewcommand{\edG}{{\e}_{\dG}}
\renewcommand{\eM}{{\e}_{\M}}
\renewcommand{\eN}{{\e}_{\N}}
\renewcommand{\LL}[1][]{{\rm \bf L}^{#1}}
\renewcommand{\TT}[1][]{{\rm \bf T}^{#1}}
\renewcommand{\DD}[1][]{{\rm \bf D}^{#1}}
\renewcommand{\RR}[1][]{{\rm \bf R}^{#1}}

\part{The Quasi--Coassociative Setting}
\setcounter{section}{6}
\setcounter{equation}{0}
In Part I we have reviewed the notions of left and right
$\G$--coactions and crossed products and we have introduced as new
concepts the notions of two--sided $\G$--coactions and diagonal
crossed products, where throughout $\G$ had been supposed to be a
standard coassociative Hopf algebra. As examples and applications we
have mentioned the Drinfel'd double $\D(\G)$ and the constructions of
quantum chains based on a Hopf algebra $\G$.

We now proceed to generalize the above ideas to quasi--Hopf algebras
$\G$. In Section 6 we give a short review of the definitions and
properties of quasi--Hopf algebras as introduced by Drinfel'd
[Dr2]. In Section 7 we propose an obvious generalization of the
notion of right $\G$--coactions $\rho$ on an algebra $\M$ to the case
of quasi--Hopf algebras $\G$ (and similarly for left coactions
$\la$). As for the coproduct on $\G$, the basic idea here is that
$(\rho\tp\id)\circ \rho$ and $(\id\tp\cop)\circ \rho$ are still
related by an inner automorphism, implemented by a reassociator
$\Fi{\rho} \in \M\tp\G\tp\G$. Similarly as for Drinfel'd's
reassociator $\phi \in \G\tp\G\tp\G$, $\Fi{\rho}$ is required to obey
a pentagon equation to guarantee McLane's coherence condition under
iterated rebracketings. We also generalize Drinfel'd's
notion of a twist transformation from coproducts to coactions.

It is important to realize that $\Fi{\rho}$ has to be non--trivial, if
$\phi$ is 
non--trivial. On the other hand, $\Fi{\rho}$ might be non--trivial 
even if $\phi = \eG\tp\eG\tp\eG$, in which case the above mentioned
pentagon equation reduces to a cocycle condition for $\Fi{\rho}$ as
already considered by [DT,BCM,BM].

In Section 8 we pass to two--sided $\G$--coactions $(\delta,\Psi)$,
which as in Section 3 could alternatively be considered as right
$(\G\tp\G^{cop})$--coactions in the above sense. Correspondingly,
$\Psi \in \G\tp\G\tp\M\tp\G\tp\G$ is the reassociator for $\delta$,
which is again required to obey the appropriate pentagon equation. As
in Section 3, associated with any two--sided $\G$-coaction
$(\delta,\Psi)$ we have a pair $(\la,\Fi{\la})$ and $(\rho,\Fi{\rho})$
of left and right $\G$--coactions, respectively, which however in this
case only {\it quasi--commute}. This means that there exists another
reassociator $\Fi{\la\rho} \in \G\tp\M\tp\G$ such that
\begin{equation*}
  \Fi{\la\rho} \, (\la\tp\idG)(\rho(m)) = (\idG\tp\rho)(\la(m))\,
  \Fi{\la\rho},\quad \forall m \in \M.
\end{equation*}
Also, $\Fi{\la\rho}$ obeys in a natural way two pentagon identities
involving $(\la,\Fi{\la})$ and $(\rho,\Fi{\rho})$, respectively. We
show that twist equivalence classes of two--sided coactions are in
one--to--one correspondence with twist equivalence classes of
quasi--commuting pairs of coactions, i.e. any two--sided coaction
$\delta$ is twist--equivalent to $(\la\tp \id)\circ \rho$ (and also to
$(\id\tp\rho)\circ \la$) where $\la = (\idG\tp\idM\tp\ep)\circ \delta$
and $\rho = (\ep\tp\idM\tp\idG)\circ \delta$.

In Section 9 we give a representation theoretic interpretation of
the notions of left, right and two--sided coactions by showing that
they give rise to functors $\Rep_\G \times \Rep_\M \pfeil \Rep_\M$,    
$\Rep_\M \times \Rep_\G \pfeil \Rep_\M$ and 
$\Rep_\G \times \Rep_\M \times \Rep_\G \pfeil \Rep_\M$, respectively,
furnished with natural associativity isomorphisms, obeying the
analogue of McLane's coherence conditions for monoidal categories [ML]. 

In Section 10 we use our formalism to construct, for any two--sided
$\G$--coaction $(\delta,\Psi)$ on $\M$, the left and right diagonal
crossed products $\M\relid\dG \cong \dG\relid \M$ as equivalent
associative algebra extensions of $\M$. 
Up to equivalence, these extensions only depend on the
twist--equivalence class of $\delta$'s, and therefore on the
twist--equivalence class of quasi--commuting pairs $(\la,\rho)$.
The basic strategy for
defining the multiplication rules in these diagonal crossed products
is to generalize the generating matrix formalism of Section 5 to the
quasi--coassociative setting. In this way one is naturally lead to
define $\la\rho$--intertwiners $\TT$ as in Definition \ref{Def 1.15}, where
now the coherence condition \eqref{1.51} has to be replaced by
appropriately injecting the reassociators $\Fi{\la},\Fi{\la\rho}$ and
$\Fi{\rho}$ into the l.h.s., similarly as in Drinfel'd's definition of
a quasitriangular  $R$--matrix for quasi--Hopf algebras. With these
substitutions our main result is given by the following
\label{seite}
\vspace{0.25cm}
  \begin{sloppypar} \noindent{\bf Theorem II}
{\it 
Let $\G$ be a finite dimensional quasi--Hopf algebra and let
$(\la,\Fi{\la},\rho,\Fi{\rho},\Fi{\la\rho})$ be a quasi--commuting pair
of (left and right) $\G$--coactions on an associative algebra $\M$.
\begin{enumerate}
\item Then their exists a unital associative algebra extension $\M_1
  \supset \M$ together with a linear map $\Gamma : \dG\pfeil \M_1$
  satisfying the following universal property: \\
  $\M_1$ is algebraically generated by $\M$ and $\Gamma(\dG)$ and for
  any algebra map $\gamma :\M\pfeil\A$ into some target algebra $\A$
  the relation 
  \begin{equation}
    \label{2.00}
   \gam_T (\Gamma(\vi)) = (\vi\tp\id)(\TT)
  \end{equation}
provides a one--to--one correspondence between algebra maps
$\gam_T: \M_1 \pfeil \A$ extending $\gam$ and normal elements $\TT \in
\G\tp\A$ satisfying 
\begin{align}
  \label{2.01}
\TT \, \la_\A (m) &= \rho^{op}_\A (m)\,\TT, \quad \forall m \in \M \\
  \label{2.02}
 (\Fi[312]{\rho})_\A \, \TT[13] \, (\Fi[132]{\rho\la})^{-1}_\A \, \TT[23] \,
 (\Fi{\la})_\A &= (\cop \tp\idA)(\TT),
\end{align}
where $\la_\A (m) := (\id\tp\gam)(\la(m))$, $(\Fi{\la})_\A : =
(\idG\tp\idG\tp\gam)(\Fi{\la})$, etc.
\item If $\M\subset \tilde{M}_1$ and $\tilde{\Gamma}:
  \dG\pfeil\tilde{\M}_1$ satisfy 
  the same universality property as in part 1.), then there exists a
  unique algebra isomorphism $f: \M_1\pfeil \tilde{\M}_1$ restricting
  to the identity on $\M$, such that $\tilde{\Gamma} = f\circ \Gamma$ 
\item There exist elements $p_\la \in \G\tp\M$ and $q_\rho \in \M\tp
  \G$ such that the linear maps
  \begin{align}
    \label{2.03}
   \mu_L : \dG\tp\M \ni (\vi\tp m) &\mapsto (\id\tp\vi_{(1)})(q_\rho)
   \Gamma (\vi_{(2)}) m \, \in \M_1\\
    \label{2.04}
   \mu_R : \M\tp\dG \ni (m\tp \vi)&\mapsto m \Gamma (\vi_{(1)})
   (\vi_{(2)}\tp \id)(p_\la) \, \in \M_1
  \end{align}
provide isomorphisms of vector spaces.
\end{enumerate}
  }
\end{sloppypar}\vspace{0.25cm}
Putting $\GG : = e_\mu \tp \Gamma (e^\mu) \in \G\tp\M_1$ Theorem II
implies that $\GG$ itself satisfies the defining relations
\eqref{2.01} and \eqref{2.02}. As before, we call $\GG$ the {\it universal
  $\la\rho$--intertwiner} in $\M_1$.
We remark that it is more or less straightforward to check that the
relations \eqref{2.01} and \eqref{2.02} satisfy all associativity
constraints, such that the existence of $\M_1$ and its uniqueness up
to isomorphy may not be too much of a surprise to the experts. 
In this way part 1. and 2. of Theorem II could also be proven without
requiring an antipode on $\G$.
The main non--trivial content of Theorem II is stated in part 3.,
saying that $\M_1$ may still be modeled
on the underlying spaces 
$\dG\tp\M$ or $\M\tp\dG$, respectively 
\footnote{
  To define the elements $p_\la$ and $q_\rho$ one needs an invertible
  antipode, see Eqs. (\ref{2.5.14'}), (\ref{2.5.17'}) }.
However, as a warning against
likely misunderstandings we emphasize that in general (i.e. for
$\phi_{\la\rho}\neq \eG\tp\eM\tp\eG$) neither of the
maps 
\begin{align*}
  \M\tp\dG \ni (m\tp\vi) &\mapsto m \Gamma(\vi) \in \M_1 \\
  \dG\tp\M \ni (\vi\tp m) &\mapsto \Gamma(\vi) m \in \M_1
\end{align*}
 need to be injective (nor surjective)\footnote{
In fact, we don't even know whether the map $\Gamma : \dG \pfeil \M_1$
necessarily has to be injective.}. Also, in general neither of the
linear subspaces $\Gamma(\dG)$, $\mu_L(\dG \tp \eM)$ or $\mu_R
(\eM\tp\dG)$ will be a subalgebra of $\M_1$. 
Still, the invertibility of the
 maps $\mu_{L/R}$ guarantees that there exist well defined associative
 algebra structures induced on $\dG\tp\M$ and $M\tp\dG$ via
 $\mu^{-1}_{L/R}$ from $\M_1$. As in Part~I we denote these by 
 \begin{align}
   \label{2.05}
    \dG\relilr \M &\equiv \mu^{-1}_L (\M_1) \\
   \label{2.06}
     \M\relilr \dG &\equiv \mu^{-1}_R (\M_1).
 \end{align}
They are the analogues of the left and right diagonal crossed
products, respectively, constructed 
in Proposition \ref{Prop 1.7} and Corollary \ref{Cor 1.11}.

To actually prove Theorem II we go the opposite way,
i.e.\ for any two--sided coaction $(\delta,\Psi)$ we will first
explicitly construct left and right diagonal crossed 
products $\dG\relid \M$ and $\M\relid \dG$ as equivalent algebra
extensions of $\M$ in Section 10.1. 
As in Part~I these are defined on the underlying spaces $\dG\tp\M$ and
$\M\tp\dG$, respectively.
In Section 10.2 we describe these constructions in terms of so--called
left and right {\it diagonal 
 $\delta$--implementers} $\LL$ and $\RR$ obeying the relations of
Lemma \ref{Lem 1.16}(ii) and (iii), respectively, together with
certain coherence conditions reflecting the multiplication rules in
$\dG\reli \M$ and $\M\reli\dG$. In Section 10.3 we generalize Lemma
\ref{Lem 1.16} by showing that coherent (left or right) diagonal 
$\delta$--implementers are always in one--to--one correspondence with
(although not identical to) {\it coherent $\la\rho$--intertwiners} $\TT$,
i.e. generating matrices satisfying the relations \eqref{2.01} and
\eqref{2.02} of Theorem II. This will finally lead to a
proof of Theorem II in Section 10.4, where we show that for $\delta_l
: = (\la\tp\id)\circ \rho$ and $\delta_r : = (\id\tp \rho)\circ \la$ 
any of the four choices
 $\dG \relidl \M$,  $\dG \relidr \M$,
$\M\relidl \dG$, or $\M \relidr \dG$  explicitly solve
all properties claimed in Theorem~II. Moreover, in terms of the
notations \eqref{2.05}, \eqref{2.06} we will have 
\begin{align*}
  \dG\relilr \M &= \dG\relidl \M \\
  \M\relilr \dG &= \M\relidr \dG
\end{align*}
with trivial identification.

In Section 11 we extend the examples and applications of Section 4
to the quasi--coassociative setting. In particular, taking $\M = \G$
and $(\la,\Fi{\la}) = (\rho,\Fi{\rho}) = (\cop,\phi)$ we provide an
explicit construction of the Drinfel'd double $\D(\G) : = \G\reli\dG
\cong \dG\reli\G$ for quasi--Hopf algebras $\G$, which recently has
been defined in terms of a Tannaka--Krein like reconstruction
in [M2]. We also give generalizations of the Hopf Spin
models of [NSz] and the lattice current algebras of [AFFS] to the case
of quasi--Hopf algebras. 
\setcounter{section}{5}

\section{Definitions and properties of quasi-Hopf algebras}
\setcounter{equation}{7}
In this Section we review the basic definitions and properties of
quasi--Hopf algebras as introduced by Drinfel'd in [Dr2]. As 
before algebra morphisms are always supposed to be unital. 
A {\it quasi-bialgebra} $(\G,\cop,\ep,\phi)$ is an associative algebra
$\G$ with 
unit together with algebra morphisms $\cop:\, \G\pfeil \G\tp\G$ (the
coproduct) and $\ep: \,
\G\pfeil\mathbb{C}$ (the counit), and an invertible element $\phi \in
\G\tp\G\tp\G$, such that
\begin{gather}
  \label{eq11}
   (\id\tp\cop)(\cop(a)) \phi = \phi(\cop\tp\id)(\cop(a)), \quad
   a\in\G \\
   \label{eq12}
   (\id\tp\id\tp\cop)(\phi)(\cop\tp\id\tp\id)(\phi) =
   (\e\tp\phi)(\id\tp\cop\tp\id)(\phi)(\phi\tp\e), \\
   \label{eq13}
   (\ep\tp\id)\circ\cop = \id = (\id\tp\ep)\circ\cop , \\
   \label{eq14}
   (\id\tp\ep\tp\id)(\phi)=\e\tp\e
\end{gather}
It has been remarked by Drinfel'd that \eqref{eq11} - \eqref{eq14}
also imply the identities $(\ep\tp\id\tp\id)(\phi) =
(\id\tp\id\tp\ep)(\phi) = \e\tp\e$.
A coproduct with the above properties is called {\it quasi-coassociative} and
the element $\phi$ will be called the {\it reassociator}. 
As for Hopf algebras we will use the Sweedler notation $\cop(a) =
a_{(1)}\tp a_{(2)}$, but since $\cop$ is only quasi--coassociative we
adopt the further convention 
\begin{equation*}
  (\cop\tp\id)\circ \cop (a) = a_{(1,1)}\tp a_{(1,2)} \tp a_{(2)}
  \quad \text{and} \quad (\id\tp\cop)\circ \cop(a) = a_{(1)}\tp
  a_{(2,1)}\tp a_{(2,2)}, \,\,\text{etc.}.
\end{equation*}
Furthermore, here and throughout we use the notation
\begin{equation}
  \label{eq17a}
\phi=  X^j\tp Y^j \tp Z^j;\quad  \phi^{-1} = P^j\tp Q^j\tp R^j,
\end{equation}
where we have suppressed the summation symbol. To give an example,
Eq. \eqref{eq11} written with this notation looks like
\begin{equation*}
  a_{(1)} X^j \tp a_{(2,1)} Y^j \tp a_{(2,2)} Z^j =
  X^i a_{(1,1)} \tp Y^i a_{(1,2)} \tp Z^i a_{(2)}
\end{equation*}
As in the Hopf algebra case, one may define the dual $(\dG,
\hat{\mu},\hat{\cop}, \hat{\ep})$, but one should note that
$\dG$ is a different kind of object. Indeed, since $\cop$ is not
associative, the ``multiplication'' $\hat{\mu}$ on $\dG$ fails to be
associative. On the other 
hand $(\dG,\hat{\cop},\hat{\ep})$ is a strictly coassociative coalgebra,
in contrast to $\G$. 

A quasi-bialgebra $\G$ is called a {\it quasi-Hopf algebra}, if there
exists a linear 
antimorphism $S: \, \G \rightarrow \G$ and elements $\alpha, \beta \in
\G$ satisfying for all $a \in \G$
\begin{align}
  \label{eq16}
   S(a_{(1)})\alpha a_{(2)} &= \alpha \ep (a); \,\,\quad
   a_{(1)} \beta S(a_{(2)}) = \beta \ep(a)\quad \text{and} \\
  \label{eq17}
  \sum_j X^j \beta S(Y^j)\alpha Z^j &= 1 = \sum_j S(P^j)\alpha Q^j
  \beta S(R^j).
\end{align}
The map $S$ is called an antipode. 
We will also always suppose that $S$ is invertible. Note that as opposed to
ordinary Hopf algebras, an antipode is not uniquely determined,
provided it exists. 

Together with a quasi--Hopf algebra $\G\equiv
(\G,\cop,\ep,\phi,S,\alpha,\beta)$ we also have $\G_{op}, \G^{cop}$
and $\G_{op}^{cop}$ as quasi--Hopf algebras, where ``$op$'' means
opposite multiplication and ``$cop$'' means opposite
comultiplication. The quasi--Hopf structures are obtained by putting
$\phi_{op}:=\phi^{-1}$, $\phi^{cop}:=\Fii[321]{}$,
$\phi_{op}^{cop}:=\Fi[321]{}$, $S_{op}=S^{cop} =
(S_{op}^{cop})^{-1}:=\Si$, $\alpha_{op}:= \Si(\beta)$, $\beta_{op}
:=\Si(\alpha)$, $\alpha^{cop} := \Si(\alpha)$, $\beta^{cop} :=
\Si(\beta)$, $\alpha_{op}^{cop} :=\beta$ and $\beta_{op}^{cop} :=
\alpha$.

Next we recall that the definition of a  quasi-Hopf algebra is
`twist covariant' in the following sense:
An element $F\in \G\tp\G$ which is invertible and satisfies
$(\ep\tp\id)(F) = (\id\tp\ep)(F) = \e$, induces a so--called {\it
twist transformation}
\begin{align}
  \label{eq113}
  \cop_F(a) : &= F \cop(a) F^{-1}, \\
   \label{eq114}
  \phi_F : &= (\e\tp F)\, (\id\tp\cop)(F) \, \phi \,
  (\cop\tp\id)(F^{-1})\, (F^{-1}\tp\e)
\end{align}
It has been noticed by Drinfel'd [Dr2] that
$(\G,\cop_F,\ep,\phi_F)$ is again a quasi bialgebra. Moreover,
setting 
\begin{align*}
  \alpha_F &:= S(h^i)\alpha k^i, \quad \text{where}\, \sum_i h^i\tp k^i
  = F^{-1}\\
  \beta_F &:= f^i \beta S(g^i), \quad \text{where}\, \sum_i f^i\tp g^i
  = F,
\end{align*}
$(\G,\cop_F,\ep,\phi_F, S, \alpha_F, \beta_F)$ is also a quasi-Hopf
algebra. 
This means that a twist preserves the class of quasi-Hopf algebras.

For Hopf algebras one knows that the antipode is an anti coalgebra
morphism, i.e. $\cop(a) = (S\tp S)\big(\cop^{op}(\Si(a))\big)$. For
quasi-Hopf algebras this is true only up to a twist: Following
Drinfel'd we define the elemnts $\gamma, \delta \in \G\tp\G$ by setting
\begin{align}
  \label{eq116}
   \gamma & := \sum_i (S(U^i)\tp S(T^i))\cdot (\alpha\tp\alpha)\cdot (V^i\tp W^i)
   \\
   \label{eq117}
   \delta & := \sum_j (K^j\tp L^j)\cdot (\beta\tp\beta)\cdot (S(N^j)\tp S(M^j))
\end{align}
where
\begin{align}
  \sum_i T^i\tp U^i \tp V^i \tp W^i & =
               (\e\tp\phi^{-1})\cdot (\id\tp\id\tp\cop)(\phi), \\
  \sum_j K^j \tp L^j \tp M^j \tp N^j & =
              (\cop\tp\id\tp\id)(\phi)\cdot (\phi^{-1}\tp\e).
\end{align}
With these definitions Drinfel'd has shown, that $f\in \G\tp\G$ given by
\begin{equation}
  \label{eq118}
  f:=  (S\tp S)(\cop^{op}(P^i)) \cdot \gamma \cdot \cop(Q^i\beta R^i).
\end{equation}
defines a twist with inverse given by 
\begin{equation}
  \label{eq119}
  f^{-1} = \cop(S(P^j)\alpha Q^j) \cdot \delta \cdot (S\tp
  S)(\cop^{op}(R^j)), 
\end{equation}
such that for all $g\in\G$
\begin{equation}
  \label{eq120} f\cop(g)f^{-1} = (S\tp S)\big(\cop^{op}(\Si(g))\big).
\end{equation}
The elements $\gamma, \delta$ and the twist $f$ fulfill the relation
\begin{equation}
  \label{eq119b}
  f\, \cop(\alpha) = \gamma, \quad \cop(\beta)\, f^{-1} = \delta.
\end{equation}
Furthermore the corresponding twisted reassociator (s.(\ref{eq114})) is
given by
\begin{equation}
  \label{eq121}
  \phi_f = (S\tp S\tp S)(\phi^{321}).
\end{equation}
Taking 
$\G^{cop} \equiv
(\G,\cop^{op},\ep,\Fii[321]{},\Si,\Si(\alpha),\Si(\beta))$ instead of
$\G$, the above definition with the corresponding substitutions
furnishes a twist 
$\hat{f}$ in place of (\ref{eq118}) which is given by $\hat{f} =
(\Si\tp\Si)(f)$. Hence,
\begin{equation}
 \label{defh}
  h: =\hat{f}^{21} \equiv (\Si\tp\Si)(f^{21})
\end{equation}
 has the following properties 
\begin{align}
  \label{eq122}
  h\cop(a) h^{-1} &= (\Si\tp\Si)\big(\cop^{op}(S(a))\big)  \\
  \label{eq123}
   \phi_h &= (\Si\tp\Si\tp\Si)(\phi^{321})\\
  \label{eq122a}
  h\,\cop(\Si(\alpha)) &= (\Si\tp\Si)(\gamma^{21}).
\end{align}
These identities will be used frequently below as well as the following
\begin{corollary}
\label{corh}
  For $a \in \G$ let  $\cop_L(a) : = h \cop(a)$ and $\cop_R(a) : =
  \cop(a) h^{-1}$ where $h \in 
  \G\tp\G$ is the twist \eqref{defh}.
  Then 
  \begin{align}
    \label{eq122b}
    (\id\tp\cop_L )(\cop_L(a)) \,\phi &=
    (\Si\tp\Si\tp\Si)(\phi^{321}) \,
    (\cop_L \tp\id)(\cop_L (a))\\
    \label{eq122c}
     \phi \, (\cop_R \tp\id)(\cop_R(a)) &= (\id\tp\cop_R )(\cop_R (a))
     \, (\Si\tp\Si\tp\Si)(\phi^{321})
    ,\quad \forall a \in \G
  \end{align}
\end{corollary}
\begin{proof}
  Writing Eq. \eqref{eq123} as 
  \begin{equation*}
    (\e\tp h)\, (\id\tp\cop)(h) \, \phi =
    (\Si\tp\Si\tp\Si)(\phi^{321})\, (h\tp\e) \, (\cop\tp\id)(h),
  \end{equation*}
multiplicating both sides from the right with
$(\cop\tp\id)(\cop(a))$ and using \eqref{eq11} yields
Eq.~\eqref{eq122b}. Eq.~\eqref{eq122c} is proven analogously.
\end{proof}


\section{Coactions of quasi--Hopf algebras}
The generalization of the definition of coactions as given in
(\ref{1.1} - \ref{1.4})  
to the quasi--Hopf case is straightforward:
\begin{definition}
\label{def2.3}{\rm 
  A {\it left coaction} of a quasi-bialgebra $(\G,\eG,\cop,\ep,\phi)$ on a
  unital algebra $\M$  
  is an 
  algebra morphism $\la : \M \pfeil \G\tp\M $ 
   together with an  invertible element  
  $\Fi{\la} \in \G\tp\G\tp\M$  satisfying  
  \begin{gather}
     (\id\tp\la)(\la(m))\,\phi_{\la}  = \Fi{\la}\,(\cop\tp\id)(\la(m)),
     \qquad \forall m \in \M   \label{eq21} 
        \\ 
     (\eG\tp\Fi{\la}) (\id\tp\cop\tp\id)(\Fi{\la})(\Fi{}\tp\eM)
       =
     (\id\tp\id\tp\la)(\Fi{\la})(\cop\tp\id\tp\id)(\Fi{\la}), \label{eq22}
     \\
       (\ep\tp\id)\circ\la = \id  \label{eq22b} \\
        (\id\tp\ep\tp\id)(\Fi{\la}) = 
        (\ep\tp\id\tp\id)(\Fi{\la})  =\eG\tp\eM  \label{eq22c}
   \end{gather}
  Similarly a {\it right coaction} of $\G$ on $\M$  is
  an algebra morphism $\rho : \M \pfeil \M\tp\G $ together with
  $\Fi{\rho} \in \M\tp\G\tp\G$ such that 
  \begin{gather}
   \phi_{\rho} \, (\rho\tp\id)(\rho(m)) = 
   (\id\tp\cop)(\rho(m))\,\Fi{\rho}, 
           \qquad  \forall m \in \M   \label{eq23}
     \\       
    (\eM\tp\phi) (\id\tp\cop\tp\id)(\Fi{\rho})(\Fi{\rho}\tp\eG)
       =
      (\id\tp\id\tp\cop)(\Fi{\rho})(\rho\tp\id\tp\id)(\Fi{\rho}),\label{eq24}
      \\
       (\id\tp\ep)\circ\rho = \id \label{eq23b} \\
        (\id\tp\epsilon\tp\id)(\Fi{\rho}) =
          (\id\tp\id\tp\ep)(\Fi{\rho}) = \eM\tp\eG \label{eq23c}
  \end{gather}
  The triple $(\M,\la,\phi_\la)$ $[(\M,\rho,\phi_\rho)]$ is called a
  left [right] {\it comodule algebra} over $\G$. 
  }
\end{definition}
We remark, that of the two counit conditions in Eqs. \eqref{eq22c} and
\eqref{eq23c}, respectively, actually either one of them already
implies the other. Clearly, if $\G$ is a Hopf algebra, $\phi =
\eG\tp\eG\tp\eG $ and $\Fi{\la} = 
\eG\tp\eG\tp\eM$, one recovers the definitions given in (\ref{1.1} -
\ref{1.4}). Also, particular examples are given by $\M=\G$ and
$\la = \rho = \cop$, $\Fi{\la} = \Fi{\rho} = \phi$. In the general case
 equations (\ref{eq22}),(\ref{eq24}) 
may be understood as a  
generalized pentagon equation, whereas  (\ref{eq21}),(\ref{eq23}) mean,
that $\la$, $\rho$ respect the 
quasi-coalgebra structure of $\G$. One should notice, that because of
the pentagon equations (\ref{eq22}) and (\ref{eq24}), 
$\Fi{\la}$ and $\Fi{\rho}$ have to be nontrivial if $\phi$ is nontrivial 
(i.e. if $\G$ is not a Hopf algebra). On the other hand $\Fi{\la}$ or
$\Fi{\rho}$  
may be nontrivial even if $\phi = \eG\tp\eG\tp\eG$, i.e. if $\G$ is a
Hopf algebra. In fact, such a restricted setting has been investigated
before, see [DT, Lemma 10], [BCM, Lemma 4.5] or [BM,
Eqs. (1.2-3)]. In [BCM,BM] Eq. (\ref{eq23}) is called a ``twisted
module condition'' and Eq. \eqref{eq24} (for $\phi = \e\tp\e\tp\e$) a
``cocycle condition''. We will see in Section 10 that the twisted
crossed products considered in [DT,BCM,BM] are in fact special types
of our diagonal crossed products to be given in Definition
\ref{def2.9} below.

As has been remarked for ordinary Hopf
algebras in Section 2, after a permutation of tensor factors
$\M\tp\G \leftrightarrow \G\tp\M$ a left coaction of $\G$ may always
be considered  as a right coaction of $\G^{cop}$ (and vice versa)
where one would have to identify $\Fi{\rho}:= (\phi^{-1}_{\la})^{321}$.

As with Hopf algebras, a left coaction
$\la : \M\pfeil\G\tp\M$ induces a map $ \li : \M\tp\dG \pfeil \M$ by
\begin{equation}
  \label{eq24x}
  m \li \vi : = (\varphi \tp\id)(\la(m)), \quad \vi \in \dG, \, m \in \M
\end{equation}
which by convenient abuse of notation and terminology we still call a
``right action'' of $\dG$ on 
$\M$, despite of the fact that $\dG$ may not be an  associative
algebra. Similarly, we put $\vi \re m : = (\id\tp\vi)(\rho(m))$ and
call this a left $\dG$--action on $\M$.

Similarly as for the coproduct $\cop$ there is a natural notion of
{\it twist equivalence} for coactions of quasi--Hopf algebras. 
\begin{lemma}
  \label{lemtwist}
Let $\rho : \M\pfeil \M\tp\G$ be a right coaction of a
quasi--bialgebra $(\G,\cop,\ep,\phi)$ and let $U\in \M\tp\G$ be
invertible such that $(\id\tp\ep)(U) = \eM$. Then
\begin{equation*}
  \rho'(m) : = U\,\rho(m) \, U^{-1}
\end{equation*}
again defines a coaction $\rho' : \M\pfeil \M\tp\G$ with respect to
the same quasi--bialgebra structure on $\G$ with twisted reassociator 
\begin{equation*}
  \phi'_{\rho} = (\idM \tp\cop)(U) \, \Fi{\rho}\, (\rho\tp\idG)(U^{-1})
  \, (U^{-1}\tp \eG).
\end{equation*}
\end{lemma}
The proof of Lemma \ref{lemtwist} is straightforward and therefore
omitted. A similar statement holds for left coactions $\la$, where one
would have to take $U\in \G\tp\M$ and 
\begin{equation*}
  \phi'_{\la} = (\eG\tp U) \, (\idG\tp\la)(U) \, \Fi{\la} \,
  (\cop\tp\idM)(U^{-1}). 
\end{equation*}
Note that twisting indeed defines an equivalence relation for
coactions. Similarly, if $\cop_F$
and $\phi_F$ are given by Eqs.~\eqref{eq113} and \eqref{eq114}, then
any right (left) $\G$--coaction on $\M$ may also be considered as a
coaction with respect to the $F$--twisted structures on $\G$ by putting
$\rho_F = \rho$ $(\la_F = \la)$ and 
\begin{equation}
  \label{eq24y}
  (\Fi{\rho})_F := (\eM \tp F) \, \Fi{\rho}, \quad (\Fi{\la})_F :=
  \Fi{\la} \, (F^{-1}\tp \eM).
\end{equation}
The reader is invited to check that with these definitions Eqs
\eqref{eq22} and \eqref{eq24} are indeed also twist covariant.


\section{Two--sided coactions}
As already mentioned before, the fact that the dual $\dG$ fails to be
an associative algebra is the reason why there is no generalization of the
definitions of ordinary crossed products to the quasi--Hopf algebra case.
Nevertheless this will be possible for our diagonal crossed product
constructed from two--sided coactions. First we need 
\begin{definition}{\rm
 \label{def2.4}
  A {\it two--sided coaction} of a quasi--bialgebra
  $(\G,\cop,\ep,\phi)$ on an algebra $\M$ is an algebra map 
  $\delta :\,\M \pfeil \G\tp\M\tp\G$ together with an invertible
  element $\Psi \in \G\tp\G\tp\M\tp\G\tp\G$ satisfying
  \begin{align}
    \label{eq24c} 
    &(\id_\G\tp\delta\tp\id_\G)(\delta(m)) \,\Psi  = \Psi
    \,(\cop\tp\id_\M \tp \cop)( \delta(m)), \quad \forall m \in \M \\ 
      \notag  (\eG\tp\Psi\tp\eG) \,
   & (\id_\G\tp\cop\tp\id_\M\tp\cop\tp\id_\G)(\Psi) \,
    (\phi\tp\eM\tp\phi^{-1})
    \\ &=
    (\id_\G\tp\id_\G\tp\delta\tp\id_\G\tp\id_\G)(\Psi)\, 
    (\cop\tp\id_\G\tp\id_\M\tp\id_\G\tp\cop)(\Psi) \label{eq24e} \\ 
   & (\ep\tp\id_\M\tp\ep)\circ \delta = \id_\M \label{eq24d} \\ 
    (\id_\G\tp\ep  \tp\id_\M & \tp\ep\tp\id_\G)(\Psi) =
    (\ep\tp\id_\G\tp\id_\M\tp\id_\G\tp\ep)(\Psi) =
    \eG\tp\eM\tp\eG. \label{eq24f} 
    \end{align} }
\end{definition}
Again we remark, that either one of the two  counit axioms in
Eq. \eqref{eq24f} already implies the other. Moreover, applying $\ep_i
\tp\idM \tp \ep_j, \, 1\leq i,j \leq 3,$ to both sides of
\eqref{eq24e}, where $\ep_i : \G^{\tp^3} \rightarrow \G^{\tp^2}$ is
given by acting with $\ep$ on the $i$--th tensor factor, one gets additional
identities, which are collected in the Appendix.
 We also note,
that as in Eqs. \eqref{1.30} and \eqref{1.31} two--sided coactions
could of course be considered as right coactions of $\G\tp\G^{cop}$,
or left coactions of $\G^{cop}\tp\G$, respectively. Moreover, if
$(\delta,\Psi)$ is a two--sided coaction of $\G$ on $\M$, then
$(\delta,\Psi^{-1})$ is a two--sided coaction of $\G_{op}$ on
$\M_{op}$ and $(\delta_{op},\Psi_{op})$ is a two--sided coaction of
$\G^{cop}$ on $\M$, where 
\begin{equation}
  \label{deltaop}
  \delta_{op} : = \delta^{321}, \quad \Psi_{op} := \Psi^{54321}.
\end{equation}
An example of a two--sided coaction is given by $\M = \G$, $\delta =
(\cop\tp\id)\circ \cop$ and 
\begin{equation*}
  \Psi : = [(\id\tp\cop\tp\id)(\phi) \tp
  \e][\phi\tp\e\tp\e][(\delta\tp\id\tp\id)(\phi^{-1})] 
\end{equation*}
Similarly we could choose $\delta' = (\id\tp\cop)\circ \cop$ and 
\begin{equation*}
  \Psi' : = [\e\tp (\id\tp\cop\tp\id)(\phi^{-1})]
    [\e\tp\e\tp\phi^{-1}][(\id\tp\id\tp\delta')(\phi)].
\end{equation*}
From this example one
already realizes that in the present context the relation between
two--sided coactions and pairs of commuting left and right coactions
gets somewhat more involved as compared to Section 3, where we had
$\delta =\delta'$. First, one
easily checks that for any two--sided coaction $(\delta, \Psi)$ the
definitions
  \begin{align}
  \label{la1}
  \la &:= (\idG\tp\idM\tp\ep)\circ\delta &
      \Fi{\la}:& = (\idG\tp\idG\tp\idM\tp\ep\tp\ep)(\Psi), \\
   \label{rho1}
  \rho &:= (\ep\tp\idM\tp\idG)\circ \delta&
      \phi_{\rho}^{-1}:& = (\ep\tp\ep\tp\idM\tp\idG\tp\idG)(\Psi)
  \end{align} 
  provide us again with a left coaction $(\la,\Fi{\la})$ and a right
  coaction $(\rho,\Fi{\rho})$. Moreover, putting $\delta^{(2)} : =
  (\id\tp\delta\tp\id) \circ \delta$ we have 
  \begin{align}
    \label{la2}
   (\la\tp\idG)\circ \rho &= (\ep\tp\idG\tp\idM\tp\ep\tp\idG)\circ
   \delta^{(2)} \\
    \label{rho2}
   (\idG\tp\rho)\circ \la &= (\idG\tp\ep\tp\idM\tp\idG\tp\ep)\circ
   \delta^{(2)}
  \end{align}
However, due to the appearance of the reassociator $\Psi$ in the axiom
\eqref{eq24c}, the two expressions \eqref{la2} and \eqref{rho2} are in
general unequal, and neither one needs to coincide with
$\delta$. Instead, as we will show now, $\la$ and $\rho$ {\it
  quasi--commute} and both expressions, $(\la\tp\idG)\circ \rho$ and
$(\idG\tp\rho)\circ \la$, define two--sided coactions which are {\it
  twist equivalent} to $(\delta,\Psi)$. First we give 
\begin{definition}{\rm 
  \label{df22}
 Let $(\G,\cop,\ep,\phi)$ be a quasi--bialgebra. By a {\it quasi--commuting
 pair} of $\G$--coactions on an algebra $\M$ we mean a quintuple
 $(\la,\rho,\Fi{\la},\Fi{\rho},\Fi{\la\rho})$, where $(\la,\Fi{\la})$
 and $(\rho,\Fi{\rho})$ are left and right $\G$--coactions on $\M$,
 respectively, and where $\Fi{\la\rho} \in \G\tp\M\tp\G$ is invertible
 and satisfies 
  \begin{gather} 
    \Fi{\la\rho}\, (\la\tp\id)(\rho(m))  =  
    (\id\tp\rho)(\la(m))\, \phi_{\la\rho},  
     \qquad \forall m \in \M  \label{eq25} \\
    (\eG\tp\Fi{\la\rho}) (\id\tp\la\tp\id)(\Fi{\la\rho}) 
     (\Fi{\la}\tp\eG) 
       = 
    (\id\tp\id\tp\rho)(\Fi{\la})(\cop\tp\id\tp\id)(\Fi{\la\rho})
        \label{eq26}
     \\ 
     (\eG\tp\Fi{\rho}) (\id\tp\rho\tp\id)(\Fi{\la\rho}) 
     (\Fi{\la\rho}\tp\eG)
       = 
   (\id\tp\id\tp\cop)(\Fi{\la\rho})(\la\tp\id\tp\id)(\Fi{\rho})
      \label{eq27}
  \end{gather} }
\end{definition}
Obviously, the conditions \eqref{eq25} - \eqref{eq27} apply to the
case $\M = \G$, $\la=\rho=\cop$ and $\Fi{\la}=\Fi{\rho}=\phi$. Also
note, that acting with $(\ep\tp\ep\tp\idM\tp\idG)$ on Eq. \eqref{eq26}
and with $(\idG\tp\idM\tp\ep\tp\ep)$ on Eq. \eqref{eq27} and using the
invertibility of $\Fi{\la\rho}$ one concludes the further identities
\begin{equation}
  \label{27a}
  (\idG\tp\idM\tp\ep)(\Fi{\la\rho})=\eG\tp\eM, \quad
  (\ep\tp\idM\tp\idG)(\Fi{\la\rho})=\eM\tp\eG.  
\end{equation}
Let us also observe that quasi--commutativity is stable under
twisting. Indeed, if $U_\la \in \G\tp\M$ is a twist from
$(\la,\Fi{\la})$ to $(\la',\phi'_\la)$ and $U_\rho \in \M\tp\G$ is a
twist from $(\rho,\Fi{\rho})$ to $(\rho',\phi'_\rho)$, then 
$(\la,\rho,\Fi{\la},\Fi{\rho},\Fi{\la\rho})$ is a quasi--commuting
pair if and only if $(\la',\rho',\phi'_\la,\phi'_{\rho},\phi'_{\la\rho})$
is a quasi--commuting pair, where $\phi'_{\la\rho}$ is given by 
\begin{equation}
  \label{27b}
  \phi'_{\la\rho} = (\eG\tp U_\rho )\,(\idG\tp\rho)(U_\la )\,
  \Fi{\la\rho}\, (\la\tp\idG)(U^{-1}_\rho ) \,(U_\la \tp \eG)^{-1}.
\end{equation}
In this case we say that $(\la,\rho,\Fi{\la},\Fi{\rho},\Fi{\la\rho})$
and $(\la',\rho',\phi'_\la,\phi'_{\rho},\phi'_{\la\rho})$ are twist
equivalent as quasi--commuting pairs of coactions.

Next we point out that two--sided coactions $(\delta,\Psi)$ may be
twisted in the same fashion as one--sided ones.
\begin{definition}{\rm 
  Let $(\delta,\Psi)$ and $(\delta', \Psi')$ be two--sided coactions
  of $(\G,\cop,\ep,\phi)$ on $\M$. Then $(\delta',\Psi')$ is called
  {\it twist equivalent} to $(\delta,\Psi)$, if there exists $U\in
  \G\tp\M\tp\G$ invertible such that
  \begin{gather}
    \label{eq2.31}
    \delta' (m) = U\,\delta(m)\,U^{-1} \\
    \label{eq2.32}
    \Psi' = (\eG\tp U\tp \eG)\,(\idG\tp \delta\tp\idG)(U) \, \Psi \,
    (\cop\tp\idM\tp\cop)(U^{-1}) \\
    \label{eq2.33}
  (\ep\tp\idM\tp\ep)(U) = \eM
  \end{gather} }
\end{definition}
The reader is invited to check that for any two--sided coaction
$(\delta,\Psi)$ and any invertible $U$ satisfying \eqref{eq2.33} the
definitions \eqref{eq2.31} and \eqref{eq2.32} indeed produce another
two--sided coaction $(\delta',\Psi')$. It is also easy to see that
twisting does provide an equivalence relation between two--sided
coactions. Moreover, similarly as for one--sided coactions one readily
verifies that if $(\delta,\Psi)$ is a two--sided coaction of
$(\G,\cop,\ep,\phi)$ on $\M$, then for any twist $F\in \G\tp\G$ the
pair $(\delta,\Psi_F)$ is a two--sided coaction of
$(\G,\cop_F,\ep,\phi_F)$ on $\M$, where $\cop_F$ and $\phi_F$ are the
twisted structures on $\G$ given by Eqs. \eqref{eq113} and
\eqref{eq114}, and where 
\begin{equation}
  \label{PsiF}
  \Psi_F : = \Psi \, (F^{-1} \tp \eM\tp F^{-1})
\end{equation}
We now show that quasi--commuting pairs of $\G$--coactions
$(\la,\rho,\Fi{\la},\Fi{\rho},\Fi{\la\rho})$ always give rise to
two--sided $\G$--coactions as follows
\begin{proposition}
  \label{prop20}
  Let $(\la,\Fi{\la})$ and $(\rho,\Fi{\rho})$ be a pair of
  left and right $\G$--coactions, respectively, on $\M$, let
  $\Fi{\la\rho} \in \G\tp\M\tp\G$ be invertible and put  
  \begin{align}
    \label{eq27a}
     \delta_l&:= (\la\tp\id)\circ \rho \\
  \label{eq215}
  \Psi_l &:=   (\idG\tp\la\tp\idG^{\tp^2})\Big((\Fi{\la\rho}\tp\eG) 
      (\la\tp\idG^{\tp^2})(\phi_{\rho}^{-1})\Big)\, 
        [\Fi{\la}\tp\eG\tp\eG] \\
    \label{eq27b}
     \delta_r&:= (\id\tp\rho)\circ \la \\
  \label{eq215a}
  \Psi_r &:= 
       (\idG^{\tp^2}\tp\rho\tp\idG)\Big((\eG\tp\phi^{-1}_{\la\rho})\,
   (\idG^{\tp^2}\tp\rho)(\phi_{\la})\Big)      
        \,[\eG\tp\eG\tp\phi^{-1}_\rho] 
\end{align}
1. Consider the following conditions
\begin{itemize}
\item[(i)] $(\la,\rho,\Fi{\la},\Fi{\rho},\Fi{\la\rho})$ is a
  quasi--commuting pair of $\G$--coactions.
\item[(ii)] $(\delta_l,\Psi_l)$ is a two--sided $\G$--coaction
\item[(iii)] $(\delta_r,\Psi_r)$ is a two--sided $\G$--coaction
\end{itemize}
Then $(i)\Leftrightarrow (ii) \Leftrightarrow (iii)$, and under these
conditions $\Fi{\la\rho}$ is a twist equivalence from
$(\delta_l,\Psi_l)$ to $(\delta_r,\Psi_r)$. \\ \\
2. Under the conditions of 1.) let 
$(\la',\rho',\phi'_\la,\phi'_{\rho},\phi'_{\la\rho})$ be a
quasi--commuting pair of coactions obtained by twisting with $U_\la
\in \G\tp\M$ and $U_\rho \in \M\tp\G$, and let
$(\delta'_{l/r},\Psi'_{l,r})$ be the associated two--sided coactions. Then
$(U_\la\tp\eG)(\la\tp\idG)(U_\rho)$ is a twist from $(\delta_l,\Psi_l)$
to $(\delta'_l ,\Psi'_l)$ and $(\eG\tp U_\rho)(\idG\tp\rho)(U_\la)$ is a
twist from $(\delta_r,\Psi_r)$ 
to $(\delta'_r ,\Psi'_r)$.
\end{proposition}
\begin{proof}
  The proofs of the implications $(i)\Rightarrow (ii)$ and
  $(i)\Rightarrow (iii)$ as well as the proof of part 2.) are lengthy
  but straightforward and therefore omitted. The proof of the implication
  $(ii)\Rightarrow (i)$ is given in Appendix A, the implication $(iii)
  \Rightarrow (i)$ being analogous.
\end{proof}
Using this result we are now in the position to show that
twist--equivalence classes of quasi--commuting 
pairs of coactions $(\la,\rho,\Fi{\la},\Fi{\rho},\Fi{\la\rho})$ are in
one-to-one correspondence with twist equivalence classes of two--sided
coactions $(\delta,\Psi)$, since
up to twist equivalence any two--sided coaction is
of the type $(\delta_{l/r},\Psi_{l/r})$ given in \eqref{eq27a} -
\eqref{eq215a}. 
\begin{proposition}
  \label{prop21}
Let $(\delta,\Psi)$ be a two--sided $\G$--coaction on $\M$ and let
$(\la,\Fi{\la})$ and $(\rho,\Fi{\rho})$ be the pair of associated left
and right $\G$-coactions, respectively, given in Eqs. \eqref{la1} 
and \eqref{rho1}. Define $U_{l/r}, \Fi{\la\rho} \in \G\tp\M\tp\G$ by
\begin{align}
  \label{2.36.1}
 U_l &:= (\ep\tp\idG\tp\idM\tp\ep\tp\idG)(\Psi) \\
  \label{2.36.2}
 U_r &:= (\idG\tp\ep\tp\idM\tp\idG\tp\ep)(\Psi) \\
  \label{2.36.3}
 \Fi{\la\rho} &:= U_r U_l^{-1}
\end{align}
and let $(\delta_{l/r},\Psi_{l/r})$ be given in terms of
$(\la,\rho,\Fi{\la},\Fi{\rho},\Fi{\la\rho})$ by Eqs. \eqref{eq27a} - 
\eqref{eq215a}. Then
\begin{itemize}
\item[(i)] $(\la,\rho,\Fi{\la},\Fi{\rho},\Fi{\la\rho})$ is a
  quasi--commuting pair of $\G$--coactions.
\item[(ii)] $U_{l/r}$ provides a twist equivalence from $(\delta,\Psi)$
  to $(\delta_{l/r},\Psi_{l/r})$.
\item[(iii)] If $(\delta',\Psi')$ is twist equivalent to
  $(\delta,\Psi)$ then the associated quasi--commuting pair\\
  $(\la',\rho',\phi'_{\la},\phi'_{\rho},\phi'_{\la\rho})$ is twist
 equivalent to $(\la,\rho,\Fi{\la},\Fi{\rho},\Fi{\la\rho})$.
\end{itemize}
\end{proposition}
\begin{proof}
  Part (iii) is straightforward and omitted. Part (i) follows from
  part (ii) by part 1. of Proposition \ref{prop20}, since the twist
  equivalence in part (ii) already guarantees that $(\delta_{l/r},
  \Psi_{l/r})$ provide two--sided coactions.
The  proof of part (ii) is given in Appendix A.
\end{proof}


\section{The representation theoretic interpretation}
In this Section we provide a representation theoretic
interpretation of our notions of left, right and two--sided coactions.
From now on let $\G$ be a quasi--Hopf algebra with invertible antipode
$S$.
Let $\Rep\M$ and $\Rep\G$ be the category of unital representations of
$\M$ and $\G$, respectively, where in $\Rep\G$ we only mean to speak
of finite dimensional representations. We denote the objects in
$\Rep\G$ by $(U,\pi_U),(V,\pi_V),(W,\pi_W),\dots$, where $U,V,W \dots$
denote the underlying representation spaces and $\pi_V: \G \pfeil
\End_\Co (V)$ the representation maps. Similarly, we denote the objects in
$\Rep\M$ by
$(\gotH,\gamma_\gotH),(\gotK,\gamma_\gotK),(\gotL,\gamma_\gotL),\dots$,
where the Gothic symbols denote the representation spaces and where $
\gamma_\gotH : \M \pfeil \End_\Co(\gotH)$, etc. We will also freely use the
 $\G$--module notation by writing $ a \cdot v
: = \pi_V(a) v$ and  $V \equiv
(V,\pi_V)$ (and analogously for $\M$--modules $\gotH$). 
The set of morphisms
$\Hom_\G(U,V)$ (also called intertwiners) is given by the linear maps $f
: U \pfeil V$ satisfying  
$f\, \pi_U(a) = \pi_V(a) \, f,\,\, \forall a \in \G$.

It is well known (see e.g. [Dr2]) that for quasi--Hopf algebras $\G$
the category
$\Rep\G$ becomes a rigid monoidal category, where the tensor
product $(V\bo W, \pi_V \bo \pi_W)$ of two representations $(V,\pi_V)$ and
$(W,\pi_W)$ is given by  
\begin{equation}
  \label{2.4.1}
  V \bo W : = V\tp W \, \text{and}\, \pi_V \bo \pi_W : = (\pi_V \tp
  \pi_W)\circ \cop 
\end{equation}
whereas for morphisms ( $\equiv$ $\G$--module intertwiners) $f,g$ one has $f\bo
g : = f \tp g$. (The symbol $\tp$ always denotes the usual tensor
product in the category of vector spaces.) The associativity
isomorphisms are given in terms of the reassociator $\phi$ by the natural
family of $\G$--module isomorphisms 
\begin{equation}
  \label{2.5.2}
 \phi_{UVW} : (U\bo V)\bo W \pfeil U\bo (V\bo W), \quad \,
  \Fi{UVW} : = (\pi_U \tp \pi_V \tp \pi_W)(\phi)
\end{equation}
The unit object in $\Rep\G$ is given by $({\mathbb
  C},\ep)$. Throughout, if $\Co$ is viewed as a $\G$--module it is
always meant to be equipped the module structure given by the one
dimensional representation $\ep$. The
left and right dual of any representation $(V,\pi_V)$ are defined by
${}^*V = V^* = \hat{V} : = \Hom_{\Co}(V,\Co)$ and 
\begin{equation}
  \label{2.5.3}
 \pi_{{}^*V} : = \pi^t_V \circ S, \quad  \pi_{V^*}: = \pi_V^t \circ \Si,
\end{equation}
where ${}^t$ denotes the transposed map. The (left) rigidity structure
is given by the family of morphisms ($\G$--module intertwiners)
\begin{align}
  \label{2.5.5}
   a_V &: {}^*V\bo V \pfeil \Co, \quad\quad \hat{v} \tp v \longmapsto
   \lpa \hat{v} \mid \alpha_V \, v\rpa  \\
   \label{2.5.6} 
  b_V &: \Co \pfeil V\bo {}^* V, \quad\quad 1  \longmapsto \beta_V \,
  v_i \tp v^i,
\end{align}
where $\alpha_V \equiv \pi_V(\alpha)$, $\beta_V \equiv \pi_V(\beta)$
and where $v_i \in V$ and $v^i \in \hat{V}$ are a choice of dual
basses. Drinfel'd's antipode axioms for $\G$ precisely reflect the
fact that $a_V$ and $b_V$ are morphisms in $\Rep\G$ fulfilling 
the {\it rigidity identities}
\begin{align}
  \label{2.5.7}
  (\id_V \bo a_V) \circ \Fi{V({}^*V)V}\circ (b_V \bo \id_V) &= \id_V \\
  \label{2.5.8.}
  (a_V \bo \id_{{}^*V})\circ \phi^{-1}_{(\lV) V (\lV)} \circ (\id_{\lV}
  \bo b_V) &= \id_\lV.
\end{align}
Also note that one has ${}^*(V^*) =
({}^*V)^* = V$ with trivial identification.

Next, we recall that in any left--rigid monoidal category one has
natural isomorphisms ${}^*(U\bo V) \cong {}^*V \bo {}^*U$. It has been
mentioned by Drinfel'd [Dr2] that in our case these isomorphisms are
given by 
\begin{equation}
  \label{2.5.8a}
 f_{UV} : U\bo V \pfeil ({}^*V \bo {}^*U)^*, \quad
  u\tp v \mapsto (\pi_V\tp \pi_U)(f)\, (v\tp u)
\end{equation}
where we trivially identify the vector spaces $V\tp W \equiv
(\hat{V}\tp \hat{W})^\wedge$ and 
where the twist $f \in \G\tp\G$ is given by \eqref{eq118}. The fact
that $f_{UV}$ is indeed a morphism in $\Rep \G$ follows from
\eqref{eq120}. Similarly, we have a natural family of isomorphisms 
\begin{equation}
  \label{2.5.8b}
 h_{UV} : U\bo V \pfeil {}^*(V^* \bo U^*), \quad
  u\tp v \mapsto (\pi_V\tp \pi_U)(h)\, (v\tp u)
\end{equation}
see Eqs. \eqref{defh}, \eqref{eq122}. \\

Now a left $\G$--coaction $(\la,\Fi{\la})$ on $\M$ naturally induces a
{\it left action} of $\Rep \G$ on $\Rep \M$. By this we mean a
functor 
\begin{equation}
  \so : \, \Rep \G \times \Rep \M \pfeil \Rep \M,
\end{equation}
where for $(V,\pi) \in \Rep \G$ and $(\gotH,\gamma) \in \Rep \M$ we define $(V\so\gotH,\pi
\so \gamma) \in \Rep \M$ by 
\begin{equation}
  \label{2.5.9}
  V\so \gotH = V\tp \gotH, \quad \pi \so \gamma : = (\pi \tp \gamma)\circ \la,
\end{equation}
whereas for morphisms we put $f \so g : = f \tp g$.
The counit axiom for $\la$ implies $\ep \so \gamma_\gotH =
\gamma_\gotH$ for all 
$(\gotH,\gamma_\gotH) \in \Rep \M$ and the axioms for $\Fi{\la}$ imply the
quasi--associativity relations
\begin{equation*}
  (\pi_V \bo \pi_W) \so \gamma_\gotH \cong \pi_V\bo (\pi_W\so\gamma_\gotH)
\end{equation*}
where the isomorphism is given by 
\begin{equation}
  \label{2.5.10}
\phi_{VW\gotH}:=(\pi_V\tp\pi_W\tp \gamma_\gotH)(\phi_\la).
\end{equation}
Finally, the pentagon axiom \eqref{eq22} provides us with the analogue
of McLane's coherence conditions, i.e. the following commuting
diagram
 \unitlength0.5cm
\begin{equation}  
\label{2.5.11}
  \begin{picture}(30,8)
  \put(3.5,4){\makebox(0,0){$((U\bo V)\bo W)\so \gotH$}}
  \put(14,7){\makebox(0,0){$(U\bo V)\so(W\so \gotH)$}}
  \put(23.5,4){\makebox(0,0){$U\so (V\so(W\so \gotH))$}}
  \put(7.5,0){\makebox(0,0){$(U \bo (V \bo W)) \so \gotH$}}
  \put(20,0){\makebox(0,0){$U\so ((V \bo W) \so \gotH)$}}
  \put(11.5,0){\vector(1,0){4}}
  \put(13.5,1){\makebox(0,0){$\phi_{U(V\bo W)\gotH}$}}
    \put(3.5,3){\vector(2,-1){4}}
  \put(4.5,1.5){\makebox(0,0)[r]{$\phi_{UVW}\so\id_\gotH$}}
    \put(20.5,1){\vector(2,1){4}}
  \put(23.5,1.5){\makebox(0,0)[l]{$\id_U\so \phi_{VW\gotH}$}}
    \put(18.5,7){\vector(3,-1){6}}
    \put(5.5,6.5){\makebox(0,0)[r]{$\phi_{(U\bo V)W\gotH}$}}
    \put(22.5,6.5){\makebox(0,0)[l]{$\phi_{U V(W\so \gotH)}$}}
     \put(3.5,5){\vector(3,1){6}}
  \end{picture}
\end{equation}
With the obvious substitutions analogue statements hold for right
coactions $(\rho,\phi_\rho)$, where now these induce a right action 
$\so: \, \Rep \M \times \Rep \G \pfeil \Rep \M$ given for
$(\gamma,\gotH) \in \Rep \M$ and $(\pi,V)\in \Rep \G$ by 
\begin{equation}
  \label{2.5.12}
   \gotH\so V : = \gotH \tp V, \quad \gamma \so \pi : = (\gamma \tp
   \pi)\circ \rho.
\end{equation}
Finally, a quasi--commuting pair $(\la,\rho,\phi_\la, \phi_\rho,
\phi_{\la\rho})$ provides us with both, a left and a right action of
$\Rep \G$ on $\Rep \M$, together with a further family of
associativity equivalences \\
$(\pi_U\so \gamma_\gotH)\so \pi_V \cong \pi_U \so
(\gamma_\gotH\so\pi_V)$, where now the isomorphisms are given by 
\begin{equation}
  \label{2.5.13}
  \phi_{U\gotH V} : = (\pi_U \tp \gamma_\gotH \tp \pi_V)(\phi_{\la\rho}).
\end{equation}
Again, the conditions \eqref{eq26} and \eqref{eq27} imply further
pentagon diagrams of the type \eqref{2.5.11} with objects of the type
$UV\gotH W$ or $U\gotH VW$, respectively, in appropriate bracket
positions.

In the obvious way the above may also be generalized to arbitrary
two--sided coactions $\delta$, in which case we would obtain a functor
\begin{equation*}
  \Rep \G \times \Rep \M \times \Rep \G \pfeil \Rep \M
\end{equation*}
denoted as 
\begin{equation}
  \label{2.5.14}
 U\dol \gotH \dor V : = U\tp \gotH \tp V, \quad
 \pi_U \dol \gamma_\gotH \dor \pi_V : = (\pi_U \tp \gamma_\gotH \tp
 \pi_V)\circ \delta
\end{equation}
together with associativity isomorphisms \\
$(\pi_U\bo \pi_V) \dol \gamma_\gotH \dor (\pi_W \bo \pi_Z) \cong \pi_U
\dol (\pi_V  \dol \gamma_\gotH \dor \pi_W) \dor \pi_Z$ given by 
\begin{equation}
  \label{2.5.15}
  \Psi_{UV\gotH WZ} : = (\pi_U\tp\pi_V \tp \gamma_\gotH \tp \pi_W \tp
  \pi_Z)(\Psi) 
\end{equation}
and obeying analogue ``two--sided'' pentagon diagrams. 
Note that the operations $\dol$ and $\dor$ are {\it not} defined
individually, i.e. only the two--sided operation $\dol \cdot \dor$
makes sense. According to
Proposition \ref{prop21} the relation between two--sided $\Rep
\G$--actions and one--sided $\Rep \G$--actions is given by 
\begin{equation}
  \label{2.5.16}
 \pi_V\so \gamma_\gotH : = \pi_V \dol\gamma_\gotH\dor \ep, \quad
 \gamma_\gotH \so \pi_V : = \ep \dol\gamma_\gotH\dor \pi_V
\end{equation}
implying 
\begin{equation}
  \label{2.5.18'}
 (\pi_V\so \gamma_\gotH)\so \pi_U \cong \pi_V\dol \gamma_\gotH\dor
 \pi_U \cong \pi_V \so (\gamma_\gotH\so \pi_U) 
\end{equation}
where here the intertwiners are given by $(\pi_V \tp \gamma_\gotH \tp
\pi_U)(U_{l/r})$, respectively, see Proposition \ref{prop21}. \\

Motivated by this representation theoretic interpretation we now
introduce various natural algebraic objects associated with left,
right and two--sided coactions, respectively, which will play a
central role later, when constructing and identifying different left and
right versions of
our diagonal crossed products. Associated with any left
$\G$--coaction $(\la,\phi_\la)$ on $\M$ we define elements
$p_\la,q_\la \in \G\tp\M$ by 
\begin{align}
  \label{2.5.14'}
   p_\la :&= \Fi[2]{\la}\Si(\Fi[1]{\la}\beta)\tp\Fi[3]{\la}, \quad
    \text{where} \quad
      \phi_{\la} =
     \Fi[1]{\la}\tp\Fi[2]{\la}\tp\Fi[3]{\la}, \\
    \label{2.5.15'}
   q_\la :&= S(\Fib[1]{\la})\alpha\Fib[2]{\la}\tp\Fib[3]{\la},
   \quad \text{where} \quad
   \phi_{\la}^{-1}=\Fib[1]{\la}\tp\Fib[2]{\la}\tp\Fib[3]{\la},
   \end{align}
and where as before we have dropped summation indices and summation symbols.
Similarly, associated with any right $\G$--coaction $(\rho,\phi_\rho)$
on $\M$ we define elements $p_\rho,q_\rho \in \M\tp\G$ by 
 \begin{align}
  \label{2.5.16'}
  p_\rho :&= \Fib[1]{\rho}\tp\Fib[2]{\rho}\beta S(\Fib[3]{\rho}), \quad
   \text{where} \quad  \phi_{\rho}^{-1} =
     \Fib[1]{\rho}\tp\Fib[2]{\rho}\tp\Fib[3]{\rho}\\
  \label{2.5.17'}
  q_\rho :&= \Fi[1]{\rho}\tp\Si(\alpha \Fi[3]{\rho})\Fi[2]{\rho},
  \quad \text{where} \quad
  \Fi{\rho}=\Fi[1]{\rho}\tp\Fi[2]{\rho}\tp\Fi[3]{\rho} 
\end{align}
Here $\alpha,\beta$ are the elements introduced in
Eq. \eqref{eq16}. In the case $\M = \G$ and $(\la,\phi_\la) =
(\rho,\phi_\rho) = (\cop,\phi)$ analogues of these elements have also
been considered by [Dr2,S]. They allow for
generalizations of formulas like $\mo\tp\Si(\mz)\me = m\tp \eG$,
etc. which will be needed to generalize Lemma \ref{Lem 1.16}. 
\newpage
\begin{lemma}$\,$ 
  \label{lem22}
  \begin{enumerate}
  \item Let $(\la,\phi_\la)$ be a left $\G$--coaction on $\M$ and let
    $p_\la ,q_\la$ be given by \eqref{2.5.14'},\eqref{2.5.15'}. Then
    the following identities hold for all $m \in \M$, where $\la(m)
    \equiv \mme \tp \mo$
    \begin{align}
      \label{2.5.18}
\la(\mo) \, p_\la\, [\Si(\mme) \tp \eM] &=  p_\la \, [\eG \tp m]  \\
     \label{2.5.19}
 [S(\mme) \tp \eM] \, q_\la \, \la(\mo) &=  [\eG \tp m] \, q_\la  \\
      \label{2.5.20}
 \la(q_\la^2) \, p_\la \, [\Si(q_\la^1) \tp \eM] &=\eG\tp \eM \\
    \label{2.5.21} 
    [S(p_\la^1) \tp \eM]\, q_\la \, \la(p^2_\la) & = \eG\tp \eM 
    \end{align}
Moreover, with $f,h \in \G\tp\G$ being the twists given by
\eqref{eq118},\eqref{defh}, the following identities are valid
    \begin{align}
      \label{2.5.23}
&  \phi^{-1}_\la \, (\idG \tp \la)(p_\la) \, (\eG \tp p_\la) \notag\\
  &\quad \quad\quad =
     (\cop \tp \idM)\big(\la(\phi^3_\la) p_\la\big) \, [h^{-1} \tp \eM] \,
   [\Si(\phi^2_\la ) \tp \Si(\phi^1_\la)\tp \eM]\\
       \label{2.5.24}
& (\eG \tp q_\la) \, (\idG\tp \la)(q_\la) \, \phi_\la \notag \\
&\quad\quad\quad  = 
    [S(\Fib[2]{\la}) \tp S(\Fib[1]{\la}) \tp \eM] \, [f\tp \eM]\,
   (\cop\tp\idM)\big(q_\la \la(\Fib[3]{\la})\big)
    \end{align}
\item  Similarly, let $(\rho, \phi_\rho)$ be a right $\G$--coaction on
  $\M$ and let $p_\rho,q_\rho$ be given by \eqref{2.5.16'} and
  \eqref{2.5.17'}. Then the following identities hold for all $m \in
  \M$, where $\rho(m) \equiv \mo\tp \me$.
  \begin{align}
    \label{2.5.25}
   \rho(\mo) \, p_\rho \, [\eM \tp S(\me)] & =  p_\rho \, [m\tp\eG] \\
    \label{2.5.26}
   [\eM \tp \Si(\me)] \, q_\rho \, \rho(\mo) & = [m\tp \eG ] \, q_\rho 
    \\
  \label{2.5.27}
 \rho(q^1_\rho) \, p_\rho \, [\eM \tp S(q_\rho^2)] &= \eM\tp \eG \\
  \label{2.5.28}
  [\eM \tp \Si (p_\rho^2)]\, q_\rho \, \rho(p^1_\rho) & = \eM\tp \eG
  \end{align}
  \begin{align}
    \label{2.5.29}
& \phi_\rho \, (\rho \tp \idG)(p_\rho) \, (p_\rho \tp \eG)\notag \\
 & \quad\quad\quad = 
   (\idM\tp \cop)\big( \rho(\Fib[1]{\rho}) p_\rho\big) \,
   [\eM\tp f^{-1}]\, [\eM\tp S(\Fib[3]{\rho}) \tp S(\Fib[2]{\rho})]\\
  \label{2.5.30}
&  (q_\rho \tp \eG) \, (\rho \tp \idG)(q_\rho)\,
\phi_\rho^{-1}\notag\\
 &\quad\quad\quad = 
  [\eM\tp \Si(\Fi[3]{\rho})\tp \Si(\Fi[2]{\rho})] \, [\eM\tp h]\, (\idM\tp
  \cop)\big(q_\rho \rho(\Fi[1]{\rho})\big).
  \end{align}
  \end{enumerate}
\end{lemma}
Lemma \ref{lem22} is proven in the Appendix.  Note that part 2. of Lemma
\ref{lem22} is functorially equivalent to part 1., since
$(\rho,\phi_\rho)$ is a right $\G$--coaction if and only if
$(\rho^{op}, \Fii[321]{\rho})$ is a left $\G^{cop}$--coaction. Also,
considering $(\rho, \phi^{-1}_\rho)$ as a right $\G_{op}$--coaction on
$\M_{op}$, the roles of $q_\rho$ and $p_\rho$ interchange, which makes
it enough to just prove Eqs. \eqref{2.5.26}, \eqref{2.5.28} and
\eqref{2.5.30} or the corresponding sets of equations in part 1. \\

We now give a representation theoretic interpretation of these
elements by defining for $(V,\pi_V) \in \Rep \G$ and $(\gotH,
\gamma_\gotH) \in \Rep \M$ the natural family of morphisms\footnote{
Again a summation is understood, where $\{v_i\}$ is a basis of $V$
with dual basis $\{v^i\}$}  
\begin{align}
\label{2.5.31}
  P_{V\gotH} &: \gotH \pfeil V^* \so (V\so \gotH), \quad 
   \goth \mapsto v^i \tp p_\la \cdot (v_i \tp\goth) \\
\label{2.5.32}
  Q_{V\gotH} &: {}^*V \so (V\so \gotH) \pfeil \gotH, \quad
 \hat{v} \tp v \tp\goth \mapsto (\hat{v}\tp\id)\Big(
    q_\la \cdot (v\tp\goth )\Big)
\end{align}
Eqs.~\eqref{2.5.18} and \eqref{2.5.19} then say that these are indeed
morphisms in $\Rep \M$, and Eqs. \eqref{2.5.20} and \eqref{2.5.21}
imply the ``generalized left rigidity'' identities
\begin{align}
  \label{2.5.33}
Q_{V^* (V\so\gotH)} \circ (\id_V \so P_{V\gotH}) & = \id_{V\so \gotH} \\
(\id_V \so Q_{V\gotH}) \circ P_{{}^*V (V \so \gotH)} & = \id_{V\so \gotH}
\end{align}
Finally, Eqs. \eqref{2.5.23} and \eqref{2.5.24} imply the coherence
conditions given by the following commuting diagrams:
\begin{equation}
  \label{2.5.34}
  \begin{picture}(40,12)(-3,0)
    \put(3,1){\makebox(0,0){$(V^* \bo U^*) \so [(U\bo V) \so \gotH]$}}
    \put(8,1){\vector(1,0){11}}
    \put(13.5,2){\makebox(0,0){$\phi_{V^*U^*[(U\bo V) \so \gotH]}$}}
    \put(24,1){\makebox(0,0){$V^* \so \big[ U^* \so ((U\bo V) \so
        \gotH)\big]$}}
    \put(3,5){\vector(0,-1){3}}
   \put(2.5,3.5){\makebox(0,0)[r]{$f^{-1}_{V^*U^*}\so \id_{(U\bo V)\so
         \gotH}$}} 
    \put(3,6){\makebox(0,0){$(U\bo V)^* \so [(U\bo V) \so \gotH]$}}
    \put(3,10){\vector(0,-1){3}}
   \put(2.5,8.5){\makebox(0,0)[r]{$P_{(U\bo V )\gotH}$}}
    \put(3,11){\makebox(0,0){$\gotH$}}
    \put(4,11){\vector(1,0){3}}
   \put(5.5,12){\makebox(0,0){$P_{ V \gotH}$}}
\put(10,11){\makebox(0,0){$V^* \so (V \so \gotH)$}}
 \put(13,11){\vector(1,0){6}}
  \put(16,12){\makebox(0,0){$\id_{V^*} \so P_{U(V \so \gotH)}$}}
 \put(24,11){\makebox(0,0){$V^* \so\big[ U^* \so (U \so(V \so
     \gotH))\big]$}}
    \put(24,10){\vector(0,-3){8}}
 \put(23.5,6){\makebox(0,0)[r]{$\id_{V^*} \so (\id_{U^*} \so
     \phi^{-1}_{UV\gotH})$}}  
  \end{picture}
\end{equation}
\begin{equation}
  \label{2.5.35}
  \begin{picture}(40,12)(-3,0)
    \put(3,1){\makebox(0,0){$(\lV \bo \lU) \so [(U\bo V) \so \gotH]$}}
    \put(19,1){\vector(-1,0){11}}
    \put(13.5,2){\makebox(0,0){$\phi^{-1}_{\lV \lU[(U\bo V) \so \gotH]}$}}
    \put(24,1){\makebox(0,0){$\lV \so \big[ \lU \so ((U\bo V) \so
        \gotH)\big]$}}
    \put(3,2){\vector(0,1){3}}
   \put(2.5,3.5){\makebox(0,0)[r]{$h_{\lV \lU}\so \id_{(U\bo V)\so
         \gotH}$}} 
    \put(3,6){\makebox(0,0){${}^*(U\bo V) \so [(U\bo V) \so \gotH]$}}
    \put(3,7){\vector(0,1){3}}
   \put(2.5,8.5){\makebox(0,0)[r]{$Q_{(U\bo V )\gotH}$}}
    \put(3,11){\makebox(0,0){$\gotH$}}
    \put(7,11){\vector(-1,0){3}}
   \put(5.5,12){\makebox(0,0){$Q_{ V \gotH}$}}
\put(10,11){\makebox(0,0){$\lV \so (V \so \gotH)$}}
 \put(19,11){\vector(-1,0){6}}
  \put(16,12){\makebox(0,0){$\id_{\lV} \so Q_{U(V \so \gotH)}$}}
 \put(24,11){\makebox(0,0){$\lV \so\big[ \lU \so (U \so(V \so
     \gotH))\big]$}}
    \put(24,2){\vector(0,1){8}}
 \put(23.5,6){\makebox(0,0)[r]{$\id_{\lV} \so (\id_{\lU} \so
     \phi_{UV\gotH})$}}  
  \end{picture}
\end{equation}
Similar statements hold of course for the natural family of morphisms 
\begin{align}
  \label{2.5.36}
 P_{\gotH V} &: \gotH \pfeil (\gotH \so V)\so \lV, \quad 
 \goth \mapsto p_\rho \cdot (\goth \tp v_i)  \tp v^i \\
  \label{2.5.37}
 Q_{\gotH V} &:  (\gotH \so V) \so V^* \pfeil \gotH, \quad
  \goth \tp  v\tp \hat{v} \mapsto (\id\tp \hat{v})\Big(
  q_\rho\cdot(\goth\tp v)\Big) 
\end{align}
Later we will also need some additional identities in the case where
$(\la,\phi_\la, \rho, \phi_\rho,\phi_{\la\rho})$ is a quasi--commuting
pair of coactions.
\begin{lemma}
  \label{lem22a}
Let $(\la,\phi_\la, \rho, \phi_\rho,\phi_{\la\rho})$ be a
quasi--commuting pair of $\G$--coactions on $\M$ and let $p_{\la
  /\rho}, \, q_{\la / \rho}$ be given by Eqs. \eqref{2.5.14'} - 
\eqref{2.5.17'}. Then putting $\Fib{\la\rho} \equiv
\phi^{-1}_{\la\rho}$
\begin{align}
  \label{2.5.38}
 \phi^{-1}_{\la\rho} \, (\idG \tp \rho)(p_\la) &= [\la(\Fi[2]{\la\rho})
 p_\la \tp \Fi[3]{\la\rho}]\,[\Si(\Fi[1]{\la\rho}) \tp \eM\tp \eG]\\
  \label{2.5.39}
  (\idG\tp \rho)(q_\la) \, \Fi{\la\rho} &= 
   [S(\Fib[1]{\la\rho}) \tp \eM\tp\eG] \, [q_\la \la(\Fib[2]{\la\rho})
   \tp \Fib[3]{\la\rho}] \\
 \label{2.5.40}
  \Fi{\la\rho} \, (\la\tp \idG)(p_\rho) & = 
  [\Fib[1]{\la\rho} \tp \rho(\Fi[2]{\la\rho}) p_\rho ]\,
 [\eG\tp\eM\tp S(\Fib[3]{\la\rho})] \\
 \label{2.5.41}
  (\la\tp \idG)(q_\rho) \, \phi^{-1}_{\la\rho} &= 
   [\eG\tp \eM\tp\Si(\Fi[3]{\la\rho})]\, [\Fi[1]{\la\rho} \tp q_\rho
   \rho(\Fi[2]{\la\rho})] 
\end{align}
\end{lemma}
Lemma \ref{lem22a} is also proven in the appendix. Again we remark that
for functorial reasons Eqs.\ \eqref{2.5.38} - \eqref{2.5.41} are all
equivalent, see the arguments after Lemma \ref{lem22}. Note that for
example Eq.\ \eqref{2.5.38} gives rise to the commuting diagram
\begin{equation}
  \label{2.5.42}
  \begin{picture}(22,8)
   \put(3,1){\makebox(0,0){$\big[U^* \so (U \so \gotH)\big] \so V$}}
   \put(7,1){\vector(1,0){6}}
   \put(10,2){\makebox(0,0){$\phi_{U^*(U\so \gotH )V}$}}
   \put(17,1){\makebox(0,0){$U^*\so \big[(U \so \gotH) \so V\big]$}}
      \put(3,5){\vector(0,-1){3}}
    \put(2.5,3.5){\makebox(0,0)[r]{$P_{U\gotH}\so \id_V$}}
       \put(17,5){\vector(0,-1){3}}
    \put(17.5,3.5){\makebox(0,0)[l]{$\id_{U^*} \so \phi^{-1}_{U \gotH V}$}}
   \put(3,6){\makebox(0,0){$ \gotH \so V$}}
   \put(7,6){\vector(1,0){6}}
  \put(10,7){\makebox(0,0){$P_{U(\gotH \so V)}$}}
   \put(17,6){\makebox(0,0){$U^*\so \big[U \so (\gotH \so V)\big]$}}
  \end{picture}
\end{equation}
Similar diagrams follow from \eqref{2.5.39} - \eqref{2.5.41}.

We conclude this section with giving the analogue of the elements
$p_{\la /\rho },\, q_{\la /\rho}$ for two--sided coactions. 
\begin{lemma}
  \label{lem22b}
  Let $(\delta,\Psi)$ be a two--sided $\G$--coactions on $\M$ and
  define $p_\delta ,\,q_\delta \in \G\tp \M\tp \G$ by 
  \begin{align}
  \label{2.5.43}
    p_\delta &: = \Psi^2 \Si(\Psi^1 \beta) \tp \Psi ^3 \tp \Psi^4
    \beta S(\Psi^5) \\
  \label{2.5.44}
   q_\delta &: = S(\Pb^1) \alpha \Pb^2 \tp \Pb^3 \tp \Si (\alpha
   \Pb^5) \Pb^4 
  \end{align}
where $\Pb \equiv \Psi^{-1}$. Then the following identities hold for
all $m \in \M$, where $\delta(m) \equiv \mme\tp\mo\tp\me$ and where $f,h \in
\G\tp\G$ are the twists defined in \eqref{eq118} and \eqref{defh}.
\begin{align}
  \label{2.5.45}
& p_\delta \, [\eG\tp m \tp \eG ]  =
 \delta(\mo) \, p_\delta \, [\Si(\mme)\tp \eM \tp S(\me)] \\
  \label{2.5.46}
 & [\eG\tp m \tp \eG]\, q_\delta = [S(\mme) \tp \eM \tp \Si(\me)] \,
 q_\delta\, \delta(\mo) \\
  \label{2.5.47}
 & \delta(q^2_\delta) \, p_\delta \, 
  [\Si(q^1_\delta) \tp \eM \tp S(q^3_\delta)] = \eG\tp \eM \tp \eG \\
  \label{2.5.48}
&  [S(p_\delta^1) \tp \eM \tp \Si(p_\delta^3)] \, q_\delta \,
  \delta(p^2_\delta) = \eG\tp \eM\tp \eG
\end{align}
\begin{align}
  \label{2.5.49}
  &\Psi^{-1} \, (\idG\tp \delta \tp \idG)(p_\delta) \, [\eG\tp
  p_\delta \tp \eG] \\
   & \quad = (\cop\tp \idM\tp\cop)\big(\delta(\Psi^3) p_\delta\big)\,
  [h^{-1} \tp\eM\tp f^{-1}]\,
  [\Si(\Psi^2) \tp \Si(\Psi^1) \tp \eM \tp S(\Psi^5) \tp S(\Psi^4)]
  \notag\\
    \label{2.5.50}
  & [\eG\tp q_\delta \tp\eG]\, (\idG\tp \delta \tp \idG)(q_\delta)\,
  \Psi \\
  &  \quad = [S(\Pb^2)\tp S(\Pb^1) \tp \eM \tp \Si(\Pb^5)\tp
  \Si(\Pb^4)]\,
  [f\tp \eM\tp h]\,
 (\cop \tp \idM\tp \cop)\big( q_\delta \delta(\Pb^3)\big) \notag
\end{align}
\end{lemma}
\begin{proof}
  This follows immediately from Lemma \ref{lem22} by noting that
  after a permutation of tensor factors $\idG \tp \tau_{\M,\G} :
  \G\tp\M\tp\G \pfeil \G\tp\G\tp\M$ any two--sided $\G$--coaction
  becomes a left $(\G\tp \G^{cop})$--coaction.
\end{proof}
We remark that similarly as before the elements $p_\delta$ and
$q_\delta$ give rise to natural families of morphisms 
\begin{align*}
  P_{U\gotH V} &: \gotH \pfeil U^* \dol (U\dol\gotH \dor V)\dor \lV \\
  Q_{U\gotH V} &: \lU \dol (U\dol\gotH\dor V) \dor V^* \pfeil \gotH
\end{align*}
obeying analogue ``rigidity identities'' and ``coherence diagrams'' as
before. We leave the details to the reader.

We also remark without proof, that in the cases $\delta = \delta_l
\equiv (\la\tp\id)\circ \rho$ or $\delta = \delta_r \equiv
(\id\tp\rho)\circ \la$ for a quasi--commuting pair of $\G$--coactions
$(\la,\rho,\phi_\la,\phi_\rho,\phi_{\la\rho})$ the elements
$q_\delta$ and $p_\delta$ may be expressed in terms of the
$q_{\la/\rho}$'s and $p_{\la/\rho}$'s by
\begin{align}
  \label{5.57}
q_{\delta_l} &= (\eG\tp q_\rho)\, (\idG\tp \rho)(q_\la)\,
\phi_{\la\rho}\\
   \label{5.58}
p_{\delta_l} &= \phi^{-1}_{\la\rho}\, (\idG\tp\rho)(p_\la)\,(\eG\tp
p_\rho)\\
\label{5.59}
q_{\delta_r} &= (q_\la \tp \eG)\, (\la\tp\idG)(q_\rho)\,
\phi^{-1}_{\la\rho}\\
\label{5.60}
p_{\delta_r} &= \phi_{\la\rho} \, (\la\tp\idG)(p_\rho)\, (p_\la\tp\eG)
\end{align}


\section{Diagonal crossed products}
Having developed our theory of two--sided $\G$--coactions $\delta$ for
quasi--bialgebras $\G$ we are now in the position to generalize the
construction of the left and right diagonal crossed products $\dG\relid\M$ and
$\M\relid\dG$ to the quasi--coassociative setting. Before writing down
the concrete multiplication rules we would like to draw the reader's
attention to some important conceptual differences in comparison with
the results of Section 3, specifically those of Proposition \ref{Prop 1.7},
Lemma \ref{Lem 1.10} and Corollary \ref{Cor 1.11}. 

First, as already remarked, the natural
``multiplication'' $\hat{\mu} : \dG\tp\dG \pfeil \dG$ given as the
transpose of the coproduct $\cop : \G\pfeil \G\tp\G$ is {\it not}
associative. Nevertheless, we will still write $\vi\psi :=
\hat{\mu}(\vi\tp\psi), \,\vi,\psi\in\dG$, i.e. 
\begin{equation}
  \label{24.1}
 \lpa \vi\psi\mid a\rpa := \lpa \vi\tp\psi\mid\cop(a)\rpa,\quad a \in\G.
\end{equation}
Also note that the counit $\ep \equiv \edG \in \dG$ still is a unit
for $\hat{\mu}$. On the other hand, being the dual of an associative
unital algebra, $\dG$ is a coassociative counital coalgebra with
coproduct $\hat{\cop}(\vi) \equiv \vi_{(1)}\tp\vi_{(2)}$ given by
\begin{equation}
  \label{24.2}
\lpa \hat{\cop}(\vi)\mid a\tp b \rpa := \lpa \vi\mid ab\rpa , \quad a,b\in\G
\end{equation}
We also note the identities
\begin{align}
  \label{24.3}
  \hat{\cop}(\edG) &= \edG\tp\edG \\
  \label{24.4}
  \hat{\cop}(\vi\psi) &=\hat{\cop}(\vi)\hat{\cop}(\psi)
\end{align}
Denoting the natural left and right actions of $\G$ on $\dG$ by 
$a\arr \vi \equiv \vi_{(1)}\lpa\vi_{(2)}\mid a\rpa$ and $\vi\arl a
\equiv \lpa \vi_{(1)}\mid a \rpa \vi_{(2)}, \, a\in\G,\,\vi\in\dG$, this
also implies
\begin{align}
  \label{24.5}
a\arr (\vi\psi) &= (a_{(1)}\arr \vi)(a_{(2)} \arr \psi) \\
  \label{24.6}
(\vi\psi )\arl a &= (\vi \arl a_{(1)})(\psi \arl a_{(2)}).
\end{align}
The second important warning concerns the fact that although we will
have $\dG \relid \M = \dG\tp\M$ and $\M\relid\dG = \M\tp\dG$ as linear
spaces, the subspaces $\dG \tp \eM$ and $\eM\tp\dG$ will {\it not} be
subalgebras in the diagonal crossed product. On the other hand, $\M$
will naturally be embedded as the unital subalgebra $\M\cong \edG\tp\M
\cong \M\tp\edG$.

The third warning concerns the fact that $\dG\relid\M \cong
\M\relid\dG$ will still be equivalent algebra extensions of $\M$,
similarly as in Corollary \ref{Cor 1.11}. However the
subspaces $\dG\reli\eM$ and $\eM\reli\dG$ will {\it not} be mapped
onto each other under this isomorphism. (Recall that this was the case
in Eqs. \eqref{1.35} and \eqref{1.36}). 
%
%
%
\subsection{The algebras 
  $\hat{{\cal G}} \relid {\cal M}$ and ${\cal M} \relid \hat{{\cal G}}$}
We now  proceed to the details. Given a two--sided
$\G$--coaction $(\delta,\Psi)$ on 
$\M$ we still write as before 
\begin{equation}
  \label{24.7}
 \vi\re m\li\psi : =
(\psi\tp\idM\tp\vi)(\delta(m)), \quad m\in\M, \, \vi,\psi \in\dG,
\end{equation}
disregarding the fact that $\delta$ might be neither of the form
\eqref{eq27a} nor \eqref{eq27b}. We also introduce the elements $\Om_L,
\Om_R \in \G\tp\G\tp\M\tp\G\tp\G$ given by
\begin{align}
  \label{24.8}
  \Omega_L &\equiv \Omega^1_L \tp \Omega^2_L \tp
  \Omega^3_L \tp \Omega^4_L \tp \Omega^5_L 
          := (\idG\tp\idG\tp\idM\tp\Si\tp\Si)(\Psi^{-1}) \cdot h^{54} \\
  \label{24.9}
  \Om_R &\equiv  \Om^1_R \tp \Om^2_R \tp
  \Om^3_R \tp \Om^4_R \tp \Om^5_R 
          := (h^{-1})^{21}\cdot
         (\Si\tp\Si\tp\idM\tp\idG\tp\idG)(\Psi),
\end{align}
where $h \equiv (\Si\tp\Si)(f^{21}) \in \G\tp\G$ has been introduced
in \eqref{defh}. As before, we drop all summation symbols and
summation indices.
\begin{definition}{\rm
 \label{def2.9}
  Let $(\delta,\Psi)$ be a two--sided coaction of a quasi--Hopf
  algebra $\G$ on an algebra $\M$. We define the {\it left diagonal crossed
  product} $\dG\relid\M$ to be the vector space $\dG\tp\M$ equipped with
  the multiplication rule
\begin{equation}
   \label{def261}
  (\vi\reli m)(\psi\reli n) 
   := \Big[(\Om_L^1 \arr\vi\arl\Om_L^5)(\Om_L^2\arr\psi_{(2)}\arl\Om_L^4)\Big] \reli
  \Big[ \Om_L^3 (\hat{S}^{-1}(\psi_{(1)})\re m \li \psi_{(3)})\, n \Big]
\end{equation}
and we define the {\it right diagonal crossed product} $M\relid\dG$ to be
the vector space  $\M\tp\dG$ with the multiplication rule 
\begin{equation}
   \label{def262}
  (m\reli \vi)( n \reli\psi) := 
   \Big[ m\, (\vi_{(1)}\re n \li \hat{S}^{-1}(\vi_{(3)}))\,\Om_R^3 
       \Big] \reli  \Big[(\Om_R^2  \arr\vi_{(2)}\arl
       \Om_R^4)(\Om_R^1\arr\psi\arl\Om_R^5)  \Big]. 
\end{equation}}
\end{definition}
A representation theoretic interpretation of these definitions will be
given in Section 10.2, starting with Definition \ref{def2.11}.
In cases where the two--sided coaction is unambiguously understood
from the context we also write $\dG\reli\M$ and $\M\reli\dG$.
As in Section 3 the choice of placing $\dG$ either to the left or to
the right of $\M$ stems from the fact that in $\dG\reli \M$ we
have 
\begin{equation}
  \label{24.10}
(\vi\reli m)= (\vi\reli \eM)(\edG\reli m)
\end{equation}
whereas in $\M\reli\dG$ we have
\begin{equation}
  \label{24.11}
(m\reli \vi)= (m\reli \edG)(\eM\reli \vi),
\end{equation}
as one easily checks from the definitions. Also note that for $\Om_{L/R} =
\eG\tp\eG\tp\eM\tp\eG\tp\eG$ we recover the definitions of Section~3.
Furthermore, under the trivial permutation of tensor factors we
have as in \eqref{1.34'} 
\begin{equation}
  \label{24.11'}
 (\M\relid\dG)_{op} = \dG^{cop}_{op} \relidop \M_{op}
\end{equation}
where $(\M\relid \dG)_{op}$ denotes the diagonal crossed
product with opposite multiplication, and where we recall our remark
that with the definition \eqref{deltaop}
$(\delta_{op},\Psi^{-1}_{op})$ defines a two--sided coaction of 
$\G^{cop}_{op}$ on $\M_{op}$. 
We now formulate our first main result.
\begin{theorem}$\,$ 
  \label{thm2.5}
  \begin{itemize}
  \item[(i)] The diagonal crossed products $\dG\reli\M$ and
    $\M\reli\dG$ are associative algebras with unit $\edG\reli
    \eM$ and $\eM\reli\edG$, respectively.
  \item[(ii)] $\M \equiv \edG\reli \M \, \subset \dG\reli\M$ and 
              $\M\equiv  \M\reli\edG \, \subset \M\reli\dG$  are
              unital algebra inclusions.
  \item[(iii)] The algebras $\dG\reli\M$ and $\M\reli\dG$ 
    provide equivalent extensions of $\M$, the isomorphism being
    given by 
\begin{align}
\label{2.5.1.1} f: 
  \dG\reli \M \ni (\vi\reli m) &\mapsto \big(q^2 \reli (\Si(q^1) \arr
  \vi\arl q^3)\big) \cdot (m\reli \hat{\e}) \, \in \M\reli\dG\\
\label{2.5.1.2}
f^{-1}: 
  \M\reli \dG \ni (m\reli \vi) &\mapsto (\hat{\e} \reli m)\cdot
   \big((p^1 \arr
  \vi\arl \Si(p^3)) \reli p^2\big)\, \in \dG\reli\M,
\end{align}
where $p \equiv p_\delta$ and $q \equiv q_\delta$ are given in
Eqs. \eqref{2.5.43} and \eqref{2.5.44}.
  \end{itemize}
\end{theorem}
\begin{proof}
  One trivially checks the unit properties in part (i) and also the
  identities
  \begin{align}
    \label{24.14}
  (\vi\reli m)(\edG\reli n) &= (\vi \reli mn) \\
   \label{24.15}
  (n\reli \edG) (m \reli \vi) &=(nm\reli \vi)
  \end{align}
 for all $m,n \in \M$ and all $\vi \in \dG$, thereby proving part
 (ii). The proof of part (iii) is postponed to Section 10.2.

We now prove the associativity of the product in $(\dG\reli\M)$,
the case $(\M\reli\dG)$ being analogous by the remark
\eqref{24.11'}. First note that \eqref{def261} and \eqref{24.14}
immediately imply 
\begin{equation}
  \label{24.16}
  [XY](\edG\reli m) = X[Y(\edG\reli m)]
\end{equation}
for all $X,Y \in \dG\reli\M$ and all $m\in\M$. Next we show
that 
\begin{equation}
  \label{24.17}
  [X(\edG\reli m)]Y = X[(\edG\reli m)Y],\quad \forall X,Y \in
  \dG\reli\M,\, m\in\M
\end{equation}
To this end we use $(\id\tp\ep)(h) = \eG$ and therefore
$(\ep\tp\idG\tp\idM\tp\idG\tp\ep)(\Om_L) = \eG\tp\eM\tp\eG$ to conclude 
for $m,m',n \in \M$ and $\psi \in \dG$
\begin{equation}
  \label{24.18}
 (\edG\reli m)(\psi\reli n) = \psi_{(2)}\reli (\Si(\psi_{(1)}) \re m
 \li \psi_{(3)} )n
\end{equation}
and hence also
\begin{align}
 (\edG\reli m') \big[(\edG\reli m)(\psi\reli n)\big]
 & =  
  \psi_{(3)}\reli \big[(\Si(\psi_{(2)}) \re m' \li \psi_{(4)} )
   (\Si(\psi_{(1)}) \re m \li \psi_{(5)} )n \big] \notag \\
    & = 
     \psi_{(2)}\reli  (\Si(\psi_{(1)}) \re m'm \li \psi_{(3)}) \notag
     \\ \label{24.19}
    &= 
     (\edG\reli m'm)(\psi\reli n)
\end{align}
where we have used that $\delta$ is an algebra map. Moreover,
\eqref{24.18} also implies for all $\vi \in \dG$
\begin{align}
   &(\vi\reli \eM)\big[ (\edG\reli m)(\psi \reli n)\big]\notag \\
 & \quad\quad\quad =
  \big[ (\Om_L^1 \arr \vi\arl\Om_L^5)(\Om_L^2 \arr\psi_{(2)} \arl
  \Om_L^4)\big] \reli\big[ \Om_L^3 (\Si(\psi_{(1)})\re m\li \psi_{(3)}) n
  \big] \notag\\ 
 \label{24.20} & \quad\quad\quad =  (\vi\reli m)(\psi \reli n)
\end{align}
Putting \eqref{24.14}, \eqref{24.16}, \eqref{24.19} and \eqref{24.20}
together, we have proven \eqref{24.17}. 

In view of \eqref{24.14}, \eqref{24.16} and \eqref{24.17}, to finish
the proof of associativity we are now left with proving the following
two identities 
\begin{align}
  \label{24.21}
 (\vi\reli\eM)\big[(\psi\reli\eM)(\chi\reli \eM)\big] &= 
 \big[(\vi\reli\eM)(\psi\reli\eM)\big](\chi\reli \eM) \\
 \label{24.22}
 (\edG\reli m)\big[ (\vi\reli\eM)(\psi\reli \eM)\big] &=
 \big[ (\edG\reli m) (\vi\reli\eM)\big](\psi\reli \eM)
\end{align}
for all $\vi,\psi,\chi \in \dG$ and all $m\in\M$.
To prove these remaining identities we rewrite them using the
generating matrix formalism. 
Let $\LL \in \G\tp(\dG\reli\M)$ be given by 
$\LL = e_\mu \tp (e^\mu \reli \eM)$, where $\{e_\mu\}$ is a basis in
$\G$ with dual basis $\{e^\mu\}$ in $\dG$. We also abbreviate our
notation by identifying $m \equiv (\edG\reli m),\, m\in\M$. Then Eqs.
\eqref{24.21} and \eqref{24.22} are equivalent, respectively, to
\begin{align}
    \label{eq213}
    \LL[14](\LL[24]\LL[34]) &= (\LL[14]\LL[24])\LL[34] \\
    \label{eq214}
    [\eG\tp\eG\tp m] (\LL[13]\LL[23]) &= ([\eG\tp\eG\tp m] 
      \LL[13])\LL[23],
  \end{align}
where \eqref{eq213} is understood as an identity in $\G^{\tp^3} \tp
(\dG\reli\M)$ and \eqref{eq214} as an identity in 
$\G^{\tp^2} \tp (\dG\reli \M)$.
We now use that \eqref{def261} and \eqref{24.10} imply
  \begin{align}
   \label{216a}
    [\eG\tp \eM] \, \LL &= \LL \\ 
    \label{eq216b}
     [\eG\tp m]\, \LL &=[\Si(\me)\tp \eM ]\,\LL[]\,
   [\mme\tp \mo],   \quad   \forall m \in \M      \\
    \label{eq216c}
   \LL[13]\LL[23] &= [\Si(\bar{\Psi}^5) \tp \Si(\bar{\Psi}^4) \tp \eM]\,
      \, [(\cop_L\tp\id )(\LL)]\, 
                   [\bar{\Psi}^1 \tp \bar{\Psi}^2 \tp \bar{\Psi}^3] 
  \end{align}
where we have introduced the notation $\delta(m) = \mme\tp\mo\tp\me$
and $\Psi^{-1} \equiv \bar{\Psi} = \bar{\Psi}^1 \tp \bar{\Psi}^2 \tp
\bar{\Psi}^3 \tp \bar{\Psi}^4 \tp \bar{\Psi}^5$, and where $\cop_L(a)
: = h \cop(a), \, a \in \G$, has been introduced 
in Corollary \ref{corh}.

To prove (\ref{eq214}) we use
  (\ref{eq216b}) twice together with \eqref{24.17} to get for the
  r.h.s. of \eqref{eq214} 
  \begin{equation*} 
  ([\eG\tp\eG\tp m]\,\LL[13])\LL[23] =
     [\Si(\mz)\tp\Si(\me)\tp\eM]\,\LL[13]\LL[23]\,
     [\mmz\tp\mme\tp\mo],
   \end{equation*}
 where we have used the notation 
  $(\id\tp\delta\tp\id)\circ \delta(m) = \mmz\tp\mme\tp\mo\tp\me\tp\mz$.
 On the other hand, by the intertwiner property \eqref{eq122} together
 with (\ref{eq216b}),(\ref{eq216c}) the l.h.s. of \eqref{eq214} gives
   \begin{align*}
    & [\eG\tp\eG\tp m] \,(\LL[13]\LL[23]) \\
    & \quad = 
     [(\Si\tp\Si)\big(\cop^{op}(\me)(\bar{\Psi}^5\tp\bar{\Psi}^4)\big)\tp\eM]\,
     [(\cop_L\tp\id)(\LL)]\,      
       [\cop(\mme)\tp\mo]\,[\bar{\Psi}^1\tp\bar{\Psi}^2\tp\bar{\Psi}^3].
   \end{align*}
   Using again \eqref{eq216c} to rewrite the r.h.s. of this formula
   Eq. (\ref{eq214}) follows from the defining property 
   (\ref{eq24c}) of $\bar{\Psi}\equiv \Psi^{-1}$. \\ 

 To prove (\ref{eq213}), we use \eqref{eq216c} to compute for the l.h.s 
 (writing $\hat{\Psi}$ for another copy of $\Psi$)
  \begin{align}
  \LL[14](\LL[24]\LL[34])&= [\eG\tp\Si(\bar{\Psi}^5)\tp\Si(\bar{\Psi}^4)\tp\eM]
   [(\id\tp\cop_L\tp\id)(\LL[13]\LL[23])]
   [\eG\tp\bar{\Psi}^1\tp\bar{\Psi}^2\tp\bar{\Psi}^3]     
    \notag \\ &=
    [(\Si\tp\Si\tp\Si)\Big((\Phh^5\tp\cop^{op}(\Phh^4))
    (\eG\tp\bar{\Psi}^5\tp\bar{\Psi}^4) \Big)\tp\eM] \label{LLL1} \\ 
     & \quad\times \notag
    [(\id\tp\cop_L \tp\id)\circ (\cop_L\tp\id)(\LL)]\,
     [\Phh^1\tp\cop(\Phh^2)\tp\Phh^3] \,
     [\eG\tp\bar{\Psi}^1\tp\bar{\Psi}^2\tp\bar{\Psi}^3], 
  \end{align}
where for the second equality we have used the identity 
\begin{equation*}
\cop_L
(\Si(a)bc) = (\Si\tp\Si)(\cop^{op}(a)) \cop_L(b)\cop(c)
\end{equation*}
 following from (\ref{eq122}). For the
r.h.s. of (\ref{eq213}) we get:
\begin{align}
  &(\LL[14]\LL[24])\LL[34] \notag\\ &\quad = 
    [\Si(\bar{\Psi}^5)\tp\Si(\bar{\Psi}^4)\tp\Si(\bar{\Psi}^3_{(1)})\tp\eM]
    [\cop_L\tp\id\tp\id)(\LL[13]\LL[23])] 
    [\bar{\Psi}^1\tp\bar{\Psi}^2\tp\bar{\Psi}^3_{(-1)}\tp\bar{\Psi}^3_{(0)}]
   \notag \\&\quad=
   [(\Si\tp\Si\tp\Si)\Big( (\cop^{op}(\Phh^5)\tp\Phh^4) 
   (\bar{\Psi}^5\tp\bar{\Psi}^4\tp\bar{\Psi}^3_{(1)})\Big)\tp\eM] \label{LLL2} \\ 
    &\quad\quad \times  \notag
   [(\cop_L\tp\id \tp\id)\circ(\cop_L\tp\id)(\LL)]\,
   [ \cop(\Phh^1)\tp\Phh^2\tp\Phh^3 ]
   [\bar{\Psi}^1\tp\bar{\Psi}^2\tp\bar{\Psi}^3_{(-1)}\tp\bar{\Psi}^3_{(0)}],
\end{align}
  where for the first equality we have used (\ref{eq216b}) to move
  $\bar{\Psi}^3$ to the right of $\LL[34]$ and in the second equality again
  (\ref{eq122}). Now we use that by Corollary 2.1
  \begin{equation*}
    (\id\tp\cop_{L})(\cop_{L}(a)) = 
    (\Si\tp\Si\tp\Si)
    (\Fi[321]{})\,((\cop_{L}\tp\id)(\cop_{L}(a))\,\phi^{-1},
   \quad \forall a \in \G.
  \end{equation*}
Hence  \eqref{LLL1} and \eqref{LLL2}  are equal due to the
pentagon identity (\ref{eq24e}) for $\Psi$, which proves (\ref{eq213}). This
concludes the proof of parts (i) and (ii) of Theorem
\ref{thm2.5}. Part (iii) will be proven in Subsection 10.2 after Proposition
\ref{prop2.9}.
\end{proof}

Before proceeding let us shortly discuss how in the present context
one can see that ordinary crossed products $\M\>cros\!_\rho \dG$ (or
$\dG\<cros\!_\la \M$) in general cannot be defined as  associative
algebras any more. In 
the strictly  coassociative setting of Section 3 these could be
considered as special types of diagonal crossed products, where
$\delta = \eG\tp \rho$ (or $\delta = \la \tp \eG$). In the present
setting it is not clear whether such $\delta$'s give well defined
two--sided coactions, since in fact the maps 
\begin{equation*}
  \la_0 (m) := \eG\tp m, \quad \rho_0 (m) : = m\tp\eG
\end{equation*}
need not even be one--sided coactions. For this one would also need
the existence of reassociators $\phi_{\la_0}, \, \phi_{\rho_0}$
satisfying the axioms of Definition 2.1. On the other hand, suppose
there exists $\phi_{\la_0}$ such 
that $(\la_0,\phi_{\la_0})$ is a well defined left coaction and let
$(\rho,\phi_\rho)$  be a right coaction such that 
\begin{equation*}
  (\idG\tp\idG\tp\rho)(\phi_{\la_0})= \phi_{\la_0}\tp\eG
\end{equation*}
Then one immediately checks that
$(\la_0,\rho,\phi_{\la_0},\phi_\rho,\phi_{\la\rho}: =
\eG\tp\eM\tp\eG)$ provides a quasi--commuting pair of coactions and
therefore in this case  $\M\>cros\!_\rho\, \dG : =
\M\relid\dG$ would indeed be a well defined associative
algebra, where $\delta = \eG\tp\rho$ and $\Psi = (\eG\tp\eG\tp
\phi^{-1}_\rho)(\phi\!_{\la_0}\tp \eG\tp\eG)$. An analogous statement holds for
$\dG\<cros\!_{\la_0}\M$. Such a scenario may of course be produced
trivially by starting with $(\delta = \eG\tp \rho, \Psi =
\eG\tp\eG\tp\eM\tp\eG\tp\eG)$ as a two--sided coaction on a strictly
coassociative  Hopf algebra $(\G,\cop)$ and subsequently passing to a
twist equivalent quasi--Hopf algebra structure $\cop_F$ on $\G$ with
$\Psi_F$ given by \eqref{PsiF}. As another example one may take a
strictly coassociative Hopf algebra $(\G,\cop)$ with
$\phi,\phi_{\la_0},\alpha,\beta,h$ being trivial, but
$(\rho,\phi_\rho)$ being a $\G$--coaction on $\M$ in the sense of Definition
\ref{def2.3}, with a nontrivial cocycle $\phi_\rho$ as considered in
[DT,BCM,BM]. Hence $\delta := \eG\tp \rho, \Psi : = \eG\tp\eG\tp
\phi^{-1}_\rho$ would give a well defined two--sided $\G$--coaction in
the sense of our Definition \ref{def2.4}. In this case $\Om_R$ in
\eqref{24.9} would be given by
\begin{equation*}
  \Om_R = \eG\tp\eG\tp\phi^{-1}_\rho
\end{equation*}
and defining $\sigma : \dG\tp\dG \pfeil \M$ by 
\begin{equation*}
  \sigma (\vi\tp \psi) : =(\id\tp\vi\tp\psi)(\phi^{-1}_\rho)
\end{equation*}
one immediately verifies from Definition \ref{def2.9} that in this
special case  our right diagonal crossed product satisfies 
\begin{equation*}
  \M\relid \dG = \M \#_\sigma \dG
\end{equation*}
where the right hand side is the twisted crossed product considered in [DT,BCM,BM].


\subsection{Left and right diagonal $\delta$--implementers}
From the associativity proof of Theorem \ref{thm2.5} in terms of the
``generating matrix'' $\LL$ we immediately read off an analogue of
Proposition \ref{Prop 1.17} describing the conditions under which an
algebra map $\gam : \M\pfeil\A$ into some target algebra $\A$ extends
to an algebra map from the diagonal crossed products into
$\A$. First, in view of Eq. \eqref{eq216b}, given $\gam: \M\pfeil \A$,
we consider a left action $\succ$ and a right action $\prec$ of
$\G\tp\M\tp\G$ on $\G\tp\A$ given for $\XX \in \G\tp\A$, and $a,b \in
\G, m\in\M$ by 
\begin{align}
  \label{252.1}
  (a\tp m\tp b) \succ \XX &:= [b\tp\gam(m)]\, \XX\, [\Si(a)\tp \eA]\\
  \label{252.2}
  \XX \prec (a\tp m\tp b) & := [\Si(b)\tp \eA]\, \XX\, [a\tp \gam(m)]
\end{align}
Note that for $\A = \dG\reli \M$ and $\gam: \M\pfeil \A$ the
canonical inclusion, Eq.~\eqref{eq216b} now reads
\begin{equation*}
  [\eG\tp m]\,\LL = \LL \prec \delta(m), \quad \forall m \in \M
\end{equation*}
Also note that the action $\succ$ commutes with right multiplication
by $\eG\tp \gam(m)$ and the action $\prec$ commutes with left
multiplication by $\eG \tp \gam(m)$.
Thus we are lead to the following
\begin{definition}{\rm 
  \label{def2.7}
Let $\gam: \M\pfeil \A$ be an algebra map into some target
algebra $\A$ and let $(\delta,\Psi)$ be a two--sided $\G$--coaction on
$\M$. A {\it left diagonal $\delta$-- implementer} in $\A$ (with respect to
$\gam$) is an element $\LL \in \G\tp \A$ satisfying for all $m\in\M$
\begin{equation}
  \label{252.3}
  [\eG\tp \gam(m)]\,\LL = \LL \prec \delta(m)
\end{equation}
Similarly, a {\it right diagonal $\delta$--implementer} in $\A$ (with respect
to $\gam$) is an element $\RR \in \G\tp \A$ satisfying for all $m\in\M$
\begin{equation}
  \label{252.4}
\RR\,[\eG\tp \gam(m)] = \delta(m) \succ \RR 
\end{equation}
A left $\delta$--implementer $\LL$ (right $\delta$--implementer $\RR$)
is called {\it coherent} if, respectively,
\begin{align}
  \label{252.5}
\LL[13]\,\LL[23] & = [\Om_L^5 \tp \Om^4_L\tp\eA]\, (\cop\tp\id)(\LL)
\, [\Om^1_L\tp\Om_L^2\tp \gam(\Om^3_L)] \\
  \label{252.6}
\RR[13]\,\RR[23] & = [\Om_R^4 \tp \Om^5_R\tp\gam(\Om^3_R)]\, (\cop\tp\id)(\RR)
\, [\Om^2_R\tp\Om_R^1\tp \eA)], 
\end{align}
where $\Om_{L/R}$ have been defined in Eqs. \eqref{24.8} and
\eqref{24.9}.
}
\end{definition}
To unburden our terminology from now by a left (right)
$\delta$--implementer we will always mean  a left (right) diagonal
$\delta$--implementer in the sense of the above definition. We trust
that the reader will not be confused by this slight inconsistence of
terminology (which arises in comparison with Definition \ref{Def 1.3},
since two--sided coactions might also be looked upon as one--sided
ones).

As before, we also call $\LL / \RR$ {\it normal}, if
$(\ep\tp\id)(\LL/\RR) = \eA$.  Note that in the coassociative setting
of Lemma \ref{Lem 1.16} left and 
right $\delta$--implementers always coincide.
In the present context we will still have
one--to--one correspondences between left and right
$\delta$--implementers, however the
identifications will not be the trivial ones. This observation will
eventually lead to a proof of Theorem \ref{thm2.5}(iii). 
Before approaching this goal let us first note the immediate
\begin{corollary}
  \label{cor2.8}
Let $(\M,\delta,\Psi)$ be a two--sided $\G$--comodule algebra and let
$\gam : \M \pfeil \A$ be an algebra map into some target algebra $\A$. 
Then the relations 
\begin{align}
  \label{24.25}
  \gam_L (\vi\reli m) & = (\vi\tp\id)(\LL)\, \gam(m) \\
    \label{24.26}   \gam_R
 (m\reli \vi) &= \gam(m) \, (\vi\tp\id)(\RR)
\end{align}
provide one--to--one correspondences between algebra maps $\gam_L
: \dG\reli\M \pfeil \A\,$  $(\gam_R : M\reli\dG \pfeil
\A)$ extending $\gam$ and coherent left $\delta$--implementers $\LL$
(coherent right
$\delta$--implementers $\RR$), respectively, where
$\gam_L$/$\gam_R$  
is unital if and only if $\LL / \RR$ is normal.
\end{corollary}
\begin{proof}
  This follows immediately from Eqs.~\eqref{eq216b} and \eqref{eq216c}
  and the analogue relations in $\M\relid \dG$.
\end{proof}
To approach the proof of part (iii) of Theorem \ref{thm2.5}
we now show that there is a coherence preserving isomorphism between
the spaces of left and right $\delta$--implementers. 
At this point the reader may find it helpful to also consult the
representation theoretic interpretation of left and right
$\delta$--implementers given at the end of this subsection (starting
with Definition \ref{def2.11}), where Proposition \ref{prop2.9} is
expressed by the commuting diagram \eqref{252.14}.
\begin{proposition} 
  \label{prop2.9}
Under the setting of Corollary \ref{cor2.8} let $p \equiv p_\delta$
and $q \equiv q_\delta$ be given by Eqs.~\eqref{2.5.43} and
\eqref{2.5.44}. Then the map 
\begin{equation}
  \label{252.7}
 \LL \mapsto \RR : = \LL \prec p 
\end{equation}
provides a coherence and normality preserving isomorphism from the
space of left $\delta$--implementers onto the space of right
$\delta$--implementers with inverse given by
\begin{equation}
  \label{252.8}
 \RR \mapsto \LL : = q \succ \RR
\end{equation}
\end{proposition}
\begin{proof}
  Throughout we are going to drop the symbol $\gam$, in particular
  $\eM \equiv \gam(\eM) \equiv \eA$. By \eqref{eq24f}
  we have $(\ep\tp\idM\tp\ep)(p) = (\ep\tp\idM\tp\ep)(q) = \eM$ and
  therefore normality of $\LL / \RR$ implies normality of $\LL \prec
  p$ and $q \succ \RR$, respectively. If $\LL$ is a left
  $\delta$--implementer then, using \eqref{252.1} - \eqref{252.3}  
  \begin{align*}
    \delta(m) \succ (\LL \prec p) &\equiv [(\eG\tp \mo)\,\LL] \prec
    [p\, (\Si(\mme) \tp \eM\tp S(\me))] \\
    & =\LL \prec [\delta(\mo) \, p\, (\Si(\mme) \tp \eM\tp S(\me))] \\
   & = \LL \prec [p\, (\eG\tp m \tp \eG)]\\
  &\equiv (\LL \prec p) \, [\eG\tp m] 
  \end{align*}
for all $m \in \M$,
where in the third line we have used \eqref{2.5.45}. Hence $\LL \prec
p$ is a right $\delta$--implementer. Analogously, if $\RR$ is a right
$\delta$--implementer, then \eqref{2.5.46} implies that $q \succ \RR$
is a left $\delta$--implementer. Next, if $\LL$ is a left
$\delta$--implementer, then  
\begin{align*}
  q\succ(\LL\prec p) & =
  [(\eG\tp q^2)\,\LL] \prec [p \,(\Si(q^1)\tp \eM\tp S(q^3))] \\
  &= \LL \prec [\delta(q^2)\,p \, (\Si(q^1) \tp \eM \tp S(q^3))]\\
  &= \LL
\end{align*}
by \eqref{2.5.47}. Similarly, using \eqref{2.5.48}, one shows for a
right $\delta$--implementer $\RR$
\begin{equation*}
  (q\succ\RR) \prec p = [(S(p^1)\tp \eM\tp \Si(p^3))\, q \,
  \delta(p^2)] \succ \RR = \RR.
\end{equation*}
We are left to show that $\LL$ is coherent if and only if $\RR : = \LL
\prec p$ is coherent. So if $\LL$ is a coherent left
$\delta$--implementer we compute for $\RR = \LL \prec p$ from
\eqref{252.5} and \eqref{24.8}
\begin{align}
\RR[13]\,\RR[23] & = \big[\Si(p^3)\Om_L^5 \tp
\Si(p^2_{(1)}\hat{p}^3)\Om^4_L\tp\eM\big]\, (\cop\tp\id)(\LL) 
\, \big[\Om^1_L\, p^1\tp\Om_L^2\, p^2_{(-1)} \hat{p}^1\tp \Om^3_L\,
p^2_{(0)}\hat{p}^2\big] \notag\\   
     & = (\Si\tp \Si\tp \idM)\Big(f^{21}(\Pb^5 \tp \Pb^4)(p^3\tp p^2_{(1)}
     \hat{p}^3) \tp \eM\Big)\,
      \notag
\\ \label{252.9} &\quad\quad\quad (\cop\tp\id)(\LL) \,
     \big[\Pb^1 p^1 \tp \Pb^2 p^2_{(-1)} \hat{p}^1 \tp \Pb^3 p^2_{(0)}
     \hat{p}^2\big]
\end{align}
where $\hat{p}$ is another copy of $p$. On the other hand, using 
\begin{equation}
\label{252.10}
  (\cop\tp\id)(\RR) = [\cop(\Si(p^3)) \tp \eM]\, (\cop\tp\id)(\LL)\,
  [\cop(p^1) \tp p^2]
\end{equation}
one computes  
\begin{align}
&[\Om_R^4 \tp \Om^5_R\tp \Om^3_R]\, (\cop\tp\id)(\RR)
\, [\Om^2_R\tp\Om_R^1\tp \eM)] \notag  \\
  &\quad = \big[(\Om_R^4 \tp \Om^5_R)\cop(\Si(\Om^3_{R(1)}p^3)) \tp \eM\big]\,
  (\cop\tp\id)(\LL) 
\, \big[\cop(\Om^3_{R(-1)} p^1) \tp
\Om^3_{R(0)} p^2\big]\,\big[\Om^2_R\tp\Om_R^1\tp 
\eA \big]\notag \\
&\quad = (\Si\tp\Si\tp \id_\M)\Big( f^{21} \cop^{op}(\Psi^3_{(1)} p^3 )
(f^{-1})^{21} \big(S(\Psi^4)\tp S(\Psi^5)\big) \tp \eM\Big) \notag\\
  &\quad\quad\quad\quad\quad\quad \times  \label{252.11}
   (\cop\tp\id)(\LL)\, \Big[ \cop(\Psi^3_{(-1)}p^1) h^{-1} \big( \Si
   (\Psi^2)\tp \Si(\Psi^1)\big) \tp \Psi^3_{(0)} p^2\Big]
\end{align}
where in the first equation we have used \eqref{252.10} and
\eqref{252.3}, in the second equation \eqref{24.9} and \eqref{eq120}.
Comparing the r.h.s. of \eqref{252.9} and \eqref{252.11} we conclude
that they are equal due to Eq. \eqref{2.5.49}.
Hence the coherence condition \eqref{252.6} is satisfied.  This proves
that $\RR : 
= \LL \prec p$ is coherent provided $\LL$ is coherent. By
functoriality, the inverse conclusion follows analogously  from $\LL :
= q \succ \RR$ and the remark \eqref{24.11'}. This concludes the proof
of Proposition \ref{prop2.9}.
\end{proof}

\begin{proof}
[Proof of Theorem \ref{thm2.5} (iii):] 
We are now in the position to give the proof of Theorem
\ref{thm2.5} (iii).  Using the notation \eqref{252.1} the map $f:
\dG\reli \M \pfeil 
\M\reli \dG$ in \eqref{2.5.1.1} is of the form 
\begin{equation*}
  f(\vi\reli m) = (\vi\tp\id)\big( q\succ \RR\big)\, (m\reli \edG)
\end{equation*}
where $\RR = e_i\tp (\eM\reli e^i)$ is the generating matrix in
$\M\reli \dG$. By Proposition \ref{prop2.9} $\LL : = q\succ \RR$ is a
coherent left $\delta$--implementer in $\M\reli\dG$ and therefore, by
Corollary \ref{cor2.8}, $f$ is an algebra map. Also by Proposition
\ref{prop2.9}, $f$ is bijective with inverse given by 
\begin{equation*}
  f^{-1} (m\reli \vi) = (\edG\reli m)(\vi\tp \id)\big( \LL\prec p\big),
\end{equation*}
where $\LL = e_i \tp (e^i \reli \eM)$ is the generating matrix in
$\dG\reli\M$. This proves \eqref{2.5.1.2} and therefore concludes the
proof of Theorem \ref{thm2.5} (iii). 
\end{proof}
Next, we show that the diagonal crossed products associated with twist
 equivalent two--sided coactions are equivalent algebra extensions.
\begin{proposition}$\,$ 
  \label{prop2.10}
  \begin{enumerate}
  \item Let $(\delta,\Psi)$ and $(\delta',\Psi')$ be twist equivalent
    two--sided $\G$--coactions on $\M$. Then the diagonal crossed
    products $\M\relid\dG$ and $\M\relidp\dG$ are equivalent
    extensions of $\M$.
  \item Let $(\delta,\Psi)$ be a two--sided $\G$--coaction on $\M$ with
    respect to the coproduct $\cop: \G\pfeil \G\tp\G$, and let
    $(\delta, \Psi_F)$ be the two--sided coaction with respect to a
    twist equivalent coproduct $\cop_F$ on $\G$, see
    Eq. \eqref{PsiF}. Denote the associated diagonal crossed
    products by $\M\reli\dG$ and $\M\reli\dG_F$, respectively. Then
    $\M\reli\dG = \M\reli\dG_F$ with trivial identification.
  \end{enumerate}
\end{proposition}
\begin{proof}
  1. Let $U \in \G\tp\M\tp\G$ be a normal twist transformation from
  $(\delta,\Psi)$ to $(\delta',\Psi')$ and let $\RR\in \G\tp (\M\relid
  \dG)$ and $\RRp \in \G \tp (\M\relidp \dG)$ be the generating
  matrices. By Corollary \ref{cor2.8}, to provide a homomorphism 
  \begin{equation*}
    f : \M\relid \dG \pfeil \M\relidp \dG
  \end{equation*}
restricting to the identity on $\M$ we have to find a coherent normal
right $\delta$--implementer $\RRt \in \G\tp (\M\relidp \dG)$. We claim
that the canonical choice
\begin{equation}
 \label{6.45a}
  \RRt := U^{-1} \succ \RRp
\end{equation}
will do the job. Indeed, $\RRt$ obviously is a normal right 
$\delta$--implementer and one is left with checking the coherence
condition with respect to $(\delta,\Psi)$. Using \eqref{eq2.32} this
is straight forward and is left to the reader.

To prove part 2. we note that $\Psi_F = \Psi\,(F^{-1}\tp \eM\tp
F^{-1})$ implies by \eqref{24.9} $(\Om_R)_F = F^{21}\,
\Om_R\,(F^{-1})^{45}$ since the element $h\in\G\tp\G$ transforms under
a twist according to $h_F = (\Si \tp\Si )(F_{op}^{-1})\, h \,
F^{-1}$. Hence, by Definition \ref{def2.7} $\RR$ is coherent with respect to
$(\delta,\Psi,\cop)$ if and only if it is coherent with respect to
$(\delta,\Psi_F,\cop_F)$.
\end{proof}
Of course, analogous statements hold for the left diagonal crossed
products. 
We emphasize that Proposition \ref{prop2.10} implies that - in the
semisimple case - the
diagonal crossed product $\M\relid\dG$ is actually completely
determined (up to equivalence) by the ``two--sided fusion rules''
$\Rep \G\times\Rep \M \times 
\Rep \G \pfeil \Rep \M$ induced by $(\delta,\Psi)$. A more detailed
exploration of this observation will be given elsewhere. 

Moreover, in view of Proposition \ref{prop21} we may from now on
restrict ourselves to two--sided coactions of the form $(\delta,\Psi)
= (\delta_{l/r},\Psi_{l/r})$ for a quasi--commuting pair
$(\la,\phi_\la,\rho,\phi_\rho,\phi_{\la\rho})$, where $\delta_l = (\la\tp\id)\circ \rho$ and $\delta_r =
(\id\tp\rho)\circ \la$, see Eqs. \eqref{eq27a}
- \eqref{eq215a}.

In this light
it will also be appropriate to introduce as an alternative notation
consistent with \eqref{2.05},\eqref{2.06} 
\begin{align}
  \label{lr1}
  \dG \relilr \M &: = \dG \relidl \M \\
   \label{lr2}
   \M\relilr \dG &:= \M\relidr \dG
\end{align}
By Theorem \ref{thm2.5}(iii) these are equivalent to $\M\relidl\dG$
and $\dG\relidr\M$, respectively, and by Proposition \ref{prop2.10}(1.)
we have $\M\relidr\dG \cong \M\relidl\dG$ and $\dG\relidl\M \cong
\dG\relidr\M$, where according to \eqref{6.45a} and Prop. \ref{prop20}(1.)
\begin{align}
  \label{6.45b}
\LL_{\delta_l} &= \LL_{\delta_r} \prec \phi_{\la\rho} \\
  \label{6.45c}
\RR_{\delta_r} &= \phi_{\la\rho} \succ \RR_{\delta_l}
\end{align}
Thus we get four equivalent versions of diagonal crossed products
associated with any quasi--commuting pair
$(\la,\rho,\phi_\la,\phi_\rho,\phi_{\la\rho})$ of $\G$--coactions on
$\M$, all of which will be shown in Section 10.4 to be a
realization of the abstract algebra $\M_1$ in Theorem II. \\

We conclude this subsection by giving a representation theoretic
interpretation of our notion of left and right
$\delta$--implementers. To this end let $(\gotH,\gamma_\gotH)$ be a
fixed 
representation of $\M$. Putting $\A : = \End_\Co (\gotH)$, we consider
$\gamma \equiv \gamma_\gotH : \M\pfeil \A$ as an algebra map. This
leads to 
\begin{definition}{\rm
  \label{def2.11}
 Let $(\delta,\Psi)$ be a two--sided $\G$--coaction on $\M$. A
 representation $(\gotH,\gamma)$ of $\M$ is called {\it
   $\delta$--coherent} if there exists a normal coherent left
 $\delta$--implementer $\LL \in\G\tp \End_\Co (\gotH)\,$ (equivalently
right $\delta$--implementer $\RR \in \G\tp \End_\Co(\gotH)$). 
  }
\end{definition}
Corollary \ref{cor2.8} then says that a representation of $\M$ is
$\delta$--coherent if and only if it extends to a representation of
the diagonal crossed products $\dG\relid \M \cong \M\relid \dG$.

We now provide a category theoretic description of $\delta$--coherent
representations $(\gotH,\gamma)$. Associated with a left
$\delta$--implementer $\LL \in \G\tp\A$, $\A : = \End_\Co(\gotH)$, we
define a natural family of $\M$--linear morphisms 
\begin{equation}
   \label{252.12}
  l_V : V \dol \gotH \dor V^* \pfeil \gotH, \quad 
 v\tp \goth \tp \hat{v} \mapsto \Big((\hat{v} \tp \id)\circ \LL_V
 \Big)(v\tp \goth)
\end{equation}
where $(V,\pi_V) \in \Rep \G$ and $\LL_V : = (\pi_V \tp \id)(\LL)$ and
where we have used the notation \eqref{2.5.14}. Eq. \eqref{252.3}
guarantees that $l_V$ is in fact a morphism in $\Rep \M$, and
naturality means\footnote{
Actually, in the terminology of [ML,Chap.IX.4] $l_V$ and $r_V$ are
{\it dinatural} ($\equiv$ diagonal--natural) transformations. We thank
T.~Kerler for pointing this out to us.}
\begin{equation*}
  l_V \circ (g \tp \id_\gotH\tp \id_{W^*}) = l_V \circ (\id_V \tp
  \id_\gotH \tp g^t)
\end{equation*}
for any morphism $g : V \pfeil W$ in $\Rep \G$. The normality
condition for $\LL$ implies $l_\Co = \id_\gotH$ and the coherence
condition \eqref{252.5} for $\LL$ translates into the following
coherence condition 
for the $l_V$'s
\begin{equation*}
  l_{V\bo W} = l_V \circ (\id_V \tp l_W \tp \id_{V^*}) \circ
  \Om^L_{VW\gotH W^*V^*} 
\end{equation*}
where $\Om^L_{VW\gotH W^*V^*} : (V\bo W) \dol \gotH \dor (V\bo W)^*
  \pfeil V\dol (W\dol\gotH \dor W^*) \dor V^*$ is the natural
  $\M$--linear isomorphism given by 
  \begin{equation*}
    \Om^L_{VW\gotH W^* V^*} = \Psi_{VW\gotH W^* V^*} \circ
    (\id_V\tp\id_W \tp \id_\gotH \tp f^{-1}_{W^*V^*})
  \end{equation*}
see \eqref{2.5.8a}, \eqref{2.5.15} and \eqref{24.8}. Similarly, a right
$\delta$--implementer $\RR \in \G\tp\A$ gives rise to a natural
coherent family of $\M$--linear morphisms 
\begin{equation}
  \label{252.13}
 r_V: \gotH \pfeil V^* \dol \gotH \dor V, \quad \goth \mapsto v^i \tp
 \RR[21]_V(\goth \tp v_i)
\end{equation}
where $\RR_V : = (\pi_V \tp \id)(\RR)$. 
\addtocounter{footnote}{-1}
As above, naturality means\footnotemark 
\stepcounter{footnote}
\begin{equation*}
(\id_V^* \tp \id_\gotH \tp g)\circ r_V = (g^t \tp \id_\gotH \tp \id_W)
\circ r_W 
\end{equation*}
for all ($\G$--linear) morphisms $g : V \pfeil W$, and coherence means
\begin{equation*}
  r_{U\bo V} = \Om^R_{V^*U^*\gotH UV} \circ (\id_{V^*} \tp r_U \tp
  \id_V) \circ r_V
\end{equation*}
where $\Om^R_{V^*U^*\gotH UV} : V^*\dol (U^* \dol \gotH \dor U)\dor V
\pfeil (U\bo V)^* \dol \gotH \dor (U\bo V)$
is given by 
\begin{equation*}
  \Om^R_{V^*U^*\gotH UV} : = (f_{V^*U^*}\tp \id_\gotH \tp \id_U \tp \id_V)
  \circ \Psi^{-1}_{V^*U^*\gotH UV},
\end{equation*}
see Eqs.~\eqref{24.9} and \eqref{252.6}.
The transformation from left $\delta$--implementers to right
$\delta$--implementers given in Proposition \ref{prop2.9} may now be
described as follows. In analogy to Eqs. \eqref{2.5.31} and
\eqref{2.5.32} we define $\M$--linear morphisms
\begin{align*}
  P_{V\gotH W} : \quad&\gotH \pfeil V^* \dol (V\dol \gotH\dor W) \dor \lW \\ 
   &\goth  \longmapsto v^i \tp p_\delta\cdot 
   (v_i \tp \goth\tp w_j) \tp w^j \\
  Q_{V\gotH W}: \quad& \lV \dol (V\dol \gotH \dor W)\dor W^* \pfeil \gotH \\
 & \hat{v} \tp v \tp \goth \tp w \tp \hat{w} \longmapsto
 (\hat{v}\tp\id_\gotH \tp
 \hat{w})\Big( q_\delta\cdot (v\tp \goth \tp w)\Big)
\end{align*}
obeying
\begin{align*}
Q_{V^* (V\dol\gotH\dor W) \lW} \circ (\id_V \tp P_{V\gotH W}\tp \id_W) &
= \id_{V\dol \gotH \dor W} \\ 
(\id_V \tp Q_{V\gotH W} \tp \id_W) \circ P_{{}^*V (V \dol \gotH\dor W)W^*} & =
\id_{V\dol \gotH\dor W} 
\end{align*}
where we have used Eqs. \eqref{2.5.45} - \eqref{2.5.48}. Proposition 
\ref{prop2.9} then implies that the transition from $l_V$'s to $r_V$'s
and back is given by the following commuting diagram
\begin{equation}
  \label{252.14}
  \begin{picture}(25,8)
   \put(3,1){\makebox(0,0){$V\dol(V^* \dol\gotH \dor V)\dor V^*$}}
   \put(20,1){\makebox(0,0){$V^* \dol(V \dol\gotH \dor V^*)\dor V$}}
  \put(3,6){\vector(0,-1){4}}
   \put(2.5,4){\makebox(0,0)[r]{$\id_V\dol r_V\dor \id_{V^*}$}}    
\put(3,7){\makebox(0,0){$V \dol\gotH \dor V^*$}}
  \put(6,7){\vector(1,0){4}}
 \put(8,7.5){\makebox(0,0){$l_V$}}  
\put(11.5,7){\makebox(0,0){$\gotH$}}
    \put(13,7){\vector(1,0){4}}  
 \put(15,7.5){\makebox(0,0){$r_V$}}  
\put(20,7){\makebox(0,0){$V^* \dol\gotH \dor V$}}
    \put(20,2){\vector(0,1){4}}
 \put(20.5,4){\makebox(0,0)[l]{$\id_{V^*}\dol l_V\dor \id_V$}}      
\put(5,2){\vector(3,2){6}}
 \put(7.5,3){\makebox(0,0)[l]{$Q_{V^* \gotH V}$}}    
\put(16,3){\makebox(0,0)[r]{$P_{V \gotH V^*}$}}    
\put(12,6){\vector(3,-2){6}} 
  \end{picture}
\end{equation}


\subsection{Coherent $\la\rho$--intertwiners}
In this subsection we are going to generalize Lemma \ref{Lem 1.16} by
providing a normality and coherence preserving one--to--one
correspondence between right $\delta_r$--implementers $\RR$ or left
$\delta_l$--implementers $\LL$, respectively, and $\la\rho$--intertwiners
$\TT$, where $\delta_r : = (\id\tp\rho)\circ \la$ and $\delta_l : =
(\la\tp \id)\circ \rho$. This will finally lead to a proof of Theorem II.
 We start with a generalization of Definition \ref{Def 1.15}
\begin{definition}{\rm
  \label{def2.12} 
  Let $(\la,\phi_\la,\rho,\phi_\rho,\phi_{\la\rho})$ be a
  quasi--commuting pair of $\G$--coactions on $\M$ and let $\gam :
  \M\tp \A$ be a unital algebra map into some target algebra $\A$. A 
{\it $\la\rho$--intertwiner} in $\A$ (with respect to $\gam$) is an element
$\TT \in \G\tp \A$ satisfying 
\begin{equation}
  \label{253.1}
\TT\, \la_\A (m) = \rho^{op}_\A (m)\, \TT, \quad \forall m \in \A
\end{equation}
A $\la\rho$--intertwiner is called $normal$ if $(\ep\tp\id)(\TT) = \eA$ and
it is called {\it coherent}, if 
\begin{equation}
  \label{253.2}
  (\Fi[312]{\rho})_\A \, \TT[13] \, \Fii[132]{\la\rho}_\A \, \TT[23]\,
  (\Fi{\la})_\A = (\cop \tp\id)(\TT)
\end{equation}
where the index $\A$ refers to the image $\gam(\M) \subset \A$, see
also Eqs. \eqref{2.01} and \eqref{2.02}.
 }
\end{definition}
We first point out that Eq.~\eqref{253.2} is consistent with
\eqref{253.1} in the following sense
\begin{lemma}
  \label{lem2.13} 
Under the conditions of Definition \ref{def2.12} let $\TT$ be a
$\la\rho$--intertwiner in $\A$ and define $B\in \G\tp\G\tp\A$ by 
\begin{equation}
  \label{253.3}
 B: = (\Fi[312]{\rho})_\A \, \TT[13] \, \Fii[132]{\la\rho}_\A \, \TT[23]\,
  (\Fi{\la})_\A
\end{equation}
Then we have for all $m \in \M$ 
\begin{equation}
  \label{253.4}
 B\, (\cop\tp\id_\A)\big(\la(m)_\A\big) =
 (\cop\tp\id_\A)\big(\rho^{op}(m)_\A\big) \, B
\end{equation}
\end{lemma}
\begin{proof}
  This is straightforward from the intertwiner properties of 
  $\TT$ and the reassociators $\phi_\la, \phi_{\la\rho}$ and
   $\phi_\rho$, see \eqref{253.1}, \eqref{eq23}, \eqref{eq21} 
    and \eqref{eq25}.
\end{proof}
We now state the generalization of Lemma \ref{Lem 1.16}, where we
recall our notation for the left and right actions given in
Eqs.~\eqref{252.1} and \eqref{252.2}. 
\begin{proposition}$\,$
\label{prop2.15}   
\begin{enumerate}
\item 
Under the conditions of Definition \ref{def2.12} let $\delta_r : =
(\id\tp\rho)\circ \la$ and $\Psi_r \in \G\tp\G\tp\M\tp\G\tp\G$ as in
\eqref{eq215a} and let $p_\la,q_\la \in \G\tp\M$ be given by
Eqs. \eqref{2.5.14'},\eqref{2.5.15'}. Then the assignment
\begin{equation}
  \label{253.5}
  \TT \longmapsto \RR : = \TT \, (\id\tp\gam)(p_\la)
\end{equation}
provides a normality and coherence preserving isomorphism from the
space of $\la\rho$--inter\-twi\-ners onto the space of right
$\delta_r$--implementers with inverse given by
\begin{equation}
  \label{253.6}
  \RR \longmapsto \TT : = (\idG\tp\rho)(q_\la) \succ \RR 
\end{equation}
\item Similarly let $\delta_l : = (\la\tp\id)\circ \rho$ and $\Psi_l
  \in\G\tp\G\tp\M\tp\G\tp\G$ as in \eqref{eq215}, and let $p_\rho,
  q_\rho\in \M\tp\G$ be given by
  Eqs. \eqref{2.5.16'},\eqref{2.5.17'}. Then the assignment
  \begin{equation}
    \label{253.7}
    \TT \longmapsto \LL : = (\id \tp \gam)(q^{op}_\rho)\, \TT
  \end{equation}
provides a normality and coherence preserving isomorphism from the
space of $\la\rho$--inter\-twi\-ners onto the space of left
$\delta_l$--implementers with inverse given by
\begin{equation}
  \label{253.8}
  \LL \longmapsto \TT : = \LL \prec (\la\tp \idG)(p_\rho)
\end{equation}
\end{enumerate}
\end{proposition}
Again the reader may find it helpful to consult the representation
theoretic interpretation given at the end of this subsection, where
the statements of Proposition \ref{prop2.15} are expressed in terms of the
commuting diagrams \eqref{253.16} and \eqref{253.17}.
We also recall our previous results regarding the transformations
between left and right $\delta_r$--implementers 
or left and right $\delta_l$--implementers 
(Proposition \ref{prop2.9}), as well as the twist transformations
given in Eqs.~\eqref{6.45b} and \eqref{6.45c}. Using the identities
\eqref{5.57}--\eqref{5.60} together with Proposition \ref{prop2.15}
one may check that all these transformations also give rise to
commuting diagrams.
\begin{proof}
  Throughout, by a convenient abuse of notation, we are going to omit
  the symbol $\gam$. Again we only need to prove part 1, since part 2
  is functorially equivalent. 
If $\TT$ is a $\la\rho$--intertwiner then 
  \begin{align*}
    \delta_r(m) \succ \RR &\equiv \rho^{op}(\mo)\, \TT\, p_\la\,
    [\Si(\mme) \tp \eA] \\
    &= \RR \, [\eG\tp m]
  \end{align*}
by \eqref{253.1} and \eqref{2.5.18}, and therefore $\RR$ is a right
$\delta_r$--implementer. Moreover, \eqref{2.5.20} implies
\begin{align*}
  (\idG\tp\rho)(q_\la) \succ (\TT\, p_\la) &\equiv \rho^{op}(q^2_\la)
  \, \TT\, p_\la \,[\Si(q^1_\la) \tp \eM] \\
    &= \TT \, \la(q^2_\la) \, p_\la \, [\Si(q^1_\la)\tp \eM] \\
    &= \TT
\end{align*}
Conversely if $\RR$ is a right $\delta_r$--implementer, then 
\begin{align*}
  \rho^{op}(m)\, \TT &\equiv (\idG\tp\rho)\big((\eG\tp m)q_\la\big)
  \succ \RR \\
   & = \Big[ [(\idG\tp\rho)(q_\la)\, \delta_r(\mo)]\succ \RR \Big]\,
   [\mme\tp\eM] \\
  &= \TT \, \la(m)
\end{align*}
where in the second equation we have used \eqref{2.5.19} and in the
last line the right $\delta_r$--implementer property \eqref{252.4} of $\RR$.
Hence $\TT$ is a $\la\rho$--intertwiner. Moreover 
\begin{align*}
  \big[ (\idG\tp \rho)(q_\la) \succ \RR \big] \, p_\la &\equiv \
   \big[(\idG\tp\rho)\big((S(p^1_\la)\tp \eA)\,q_\la\big)\big] \succ 
   \big(\RR\,[\eG\tp p^2_\la]\big) \\
   & = (\idG\tp\rho)\Big((S(p^1_\la)\tp \eA)\, q_\la
   \,\la(p^2_\la)\Big) \succ \RR \\
   & = \RR
\end{align*}
Thus the correspondence $\TT \leftrightarrow \RR$ is one--to--one, and
since $q_\la$ and $p_\la$ are normal it is clearly normality
preserving. 

To prove that it is also coherence preserving assume now that the
$\la\rho$--intertwiner $\TT$ satisfies the coherence condition \eqref{253.2}
and let $\RR = \TT\, p_\la$. Then 
\begin{equation}
  \label{253.9}
(\cop\tp\id)(\RR)\, [h^{-1}\tp \eM] = \Fi[312]{\rho}\, \TT[13]\,
\Fii[132]{\la\rho}\, A
\end{equation}
where $A \in \G\tp\G\tp\M$ is given by
\begin{align}
  A &= \TT[23]\, \Fi{\la}\, (\cop\tp\idM)(p_\la)\,
  [h^{-1}\tp\eM]\notag\\
  &= \TT[23]\, (\idG\tp\la)\big(\la(\Fib[3]{\la}) \,p_\la\big)\,
  [\eG\tp p_\la]\,[\Si(\Fib[2]{\la})\tp \Si(\Fib[1]{\la}) \tp \eM]
  \notag\\
  &= (\idG\tp \rho^{op})\big(\la(\Fib[3]{\la})\,p_\la\big) \,
  \RR[23]\, [\Si(\Fib[2]{\la})\tp \Si(\Fib[1]{\la}) \tp \eM]
   \label{253.10}
\end{align}
Here we have used \eqref{2.5.23} in the second line and the
intertwining property
\eqref{253.1} of $\TT$ in the third line. Using the intertwiner
property \eqref{eq25} of $\phi_{\la\rho}$ and \eqref{2.5.38} we
further compute
\begin{align}
  & \TT[13] \, \Fii[132]{\la\rho} \,
  (\idG\tp\rho^{op})\big(\la(\Fib[3]{\la})\,p_\la\big)\notag \\
  &\quad\quad\quad= \label{253.11}
  \Big[(\rho^{op}\tp\id)\big(\rho(\Fib[3]{\la}) \, (\Fi[2]{\la\rho}\tp
  \Fi[3]{\la\rho})\big) \,(\RR \tp \eM)\Big]^{132} \, 
  [\Si(\Fi[1]{\la\rho}) \tp\eG\tp\eM] 
\end{align}
Putting \eqref{253.9} - \eqref{253.11} together we finally conclude 
\begin{equation}
  \label{253.12}
   (\cop\tp\id)(\RR) \, [h^{-1} \tp \eM] = 
  [\Pb_r^4 \tp \Pb^5_r \tp \Pb_r^3] \, \RR[13] \, \RR[23]\, 
  [\Si(\Pb_r^2) \tp \Si(\Pb_r^1)\tp \eM]
\end{equation}
where $\Pb_r \in \G\tp\G\tp\M\tp\G\tp\G$ is given by
\begin{equation*}
  \Pb_r : = (\eG\tp\eG\tp\phi_\rho)
  (\idG\tp\idG\tp\rho\tp\idG)\Big((\idG\tp\idG\tp\rho)(\Phi^{-1}_\la
  )\,(\eG\tp\phi_{\la\rho})\Big) 
\end{equation*}
By \eqref{eq215a} we have $\Pb_r = \Psi_r^{-1}$ and therefore
\eqref{253.12} is equivalent to the coherence condition \eqref{252.6}
for $\RR$ as a right $(\delta_r,\Psi_r)$--implementer.

Conversely, assume now that $\RR$ is a coherent right
$(\delta_r,\Psi_r)$--implementer and let \\
 $\TT : =
(\idG\tp\rho)(q_\la)\succ \RR$. To prove that $\TT$ is coherent we have
to show that
\begin{equation*}
  (\cop\tp\id)(\TT) = B
\end{equation*}
where $B$ is given by Eq. \eqref{253.3}. Now writing $\RR = \TT \,
p_\la$ and going backwards through the derivation
\eqref{253.9}$\leftarrow$\eqref{253.12} we conclude
\begin{equation}
  \label{253.13}
   (\cop\tp\id)(\TT\,p_\la) = B\, (\cop \tp\id)(p_\la)
\end{equation}
Thus, if $p_\la$ were invertible we could immediately conclude that
$\TT$ is coherent. It turns out that we may use Eq. \eqref{2.5.20} as
a substitute for the invertibility of $p_\la$, since it implies 
\begin{align*}
  (\cop\tp\id)(\TT) & \equiv (\cop
  \tp\id)\Big((\idG\tp\rho)(q_\la)\succ (\TT\, p_\la)\Big) \\
   & = (\cop\tp\id)\big(\rho^{op}(q^2_\la)\big) \, (\cop\tp\id)(\TT
   p_\la)\, \big[ \cop(\Si(q^1_\la))\tp \eM\big] \\
  & = (\cop\tp\id)\big(\rho^{op}(q^2_\la)\big) \, B\, 
    (\cop\tp\id)\Big(p_\la\,(\Si(q^1_\la)\tp \eM)\Big) \\
  &= B\, (\cop\tp\id)\Big(\la(q^2_\la)\, p_\la\, (\Si(q^1_\la)\tp
  \eM)\Big)\\
  &=B
\end{align*}
Here we have used the definition \eqref{252.1} in the second line,
\eqref{253.13} in the third line, \eqref{253.4} in the fourth line and
again \eqref{2.5.20} in the last line. Thus $\TT$ is a coherent
$\la\rho$--intertwiner, which concludes the proof of Proposition
\ref{prop2.15}. 
\end{proof}

We conclude this subsection with a representation categorical
interpretation of the notion of $\la\rho$--intertwiners. 
As before, given a fixed representation $(\gotH,\gamma)$ of $\M$ we consider
$\gamma: \M \pfeil \A \equiv \End_\Co (\gotH )$ as an algebra
map. Similarly as in Definition \ref{def2.11} we say
\begin{definition}{\rm
  \label{def2.11'}
Let $(\la,\phi_\la, \rho,\phi_\rho,\phi_{\la\rho})$ be a
quasi--commuting pair of $\G$--coactions on $\M$. A representation
$(\gotH,\gam_\gotH)$ of $\M$ is called {\it $\la\rho$--coherent}, if
there exists a normal coherent $\la\rho$--intertwiner $\TT \in \G\tp
\End_\Co (\gotH)$.
 }
\end{definition}
Proposition \ref{prop2.15} then says, that $\la\rho$--coherence is
 equivalent to $\delta_l$-coherence for $\delta_l = (\la\tp \id)\circ
 \rho$ (or to $\delta_r$--coherence for $\delta_r = (\id\tp \rho)\circ
 \la$).
 Associated with a $\la\rho$--intertwiner $\TT \in \G\tp \End_\Co (\gotH
)$ we now define a natural family of $\M$--linear morphisms 
\begin{equation}
  \label{253.14}
 t_V : V \so \gotH \pfeil \gotH \so V, \quad 
  v\tp \goth \mapsto \TT[21] (\goth \tp v),  
\end{equation}
where $(V,\pi_V) \in \Rep \G$ and $\TT[21]_V : = (\id\tp
\pi_V)(\TT[21])$. Eq. \eqref{253.1} guarantees that $t_V$ is a morphism
in $\Rep \M$ and naturality means that 
\begin{equation*}
  t_W \circ (g\tp \id_\gotH) = (\id_\gotH \tp g) \circ t_V
\end{equation*}
for any morphism $g : V \pfeil W$ in $\Rep \G$. The normality
condition  for $\TT$ implies $t_\Co = \id_\gotH$ and the coherence
condition for $\TT$ translates into the following coherence condition
for $t_V$
\begin{equation}
  \label{253.15}
  t_{V\bo W} = \phi_{V W \gotH} \circ (\id_V \tp t_W) \circ
  \phi^{-1}_{V\gotH W} \circ (t_V \tp \id_W) \circ \phi_{\gotH VW}
\end{equation}
Note that \eqref{253.15} looks precisely like one of the coherence conditions
for the braiding in a braided quasi--tensor category with
nontrivial associativity isomorphisms. Indeed, as will be shown in [HN1],
for the case $\M = \G$, the family of $t_V$'s  may be used to define a
braiding in the 
representation category of the quantum double $\D(\G) \equiv \G\reli
\dG$. 

Using the morphisms $P_{V\gotH},P_{\gotH V}$ and $Q_{V\gotH}, Q_{\gotH
  V}$ given in \eqref{2.5.31},\eqref{2.5.32} and
\eqref{2.5.36},\eqref{2.5.37}, the relation between  $\la\rho$--intertwiners
$\TT$, right $\delta_r$--implementers $\RR$ and left
$\delta_l$--implementers $\LL$ may now be described by the following
commuting diagrams connecting the intertwiner morphisms $t_V$ with the maps
$r_V$ (associated with $\RR$) and $l_V$ (associated with $\LL$) as given
in \eqref{252.12} and \eqref{252.13}.
\begin{equation}
  \label{253.16}
  \begin{picture}(25,8)
   \put(3,1){\makebox(0,0){$(\gotH \so V)\so V^*$}}
   \put(20,1){\makebox(0,0){$V^* \so(V \so\gotH )$}}
  \put(3,6){\vector(0,-1){4}}
   \put(2.5,4){\makebox(0,0)[r]{$t_V\so \id_{V^*}$}}    
\put(3,7){\makebox(0,0){$(V \so\gotH) \so V^*$}}
  \put(6,7){\vector(1,0){4}}
 \put(8,7.5){\makebox(0,0){$l_V$}}  
\put(11.5,7){\makebox(0,0){$\gotH$}}
    \put(13,7){\vector(1,0){4}}  
 \put(15,7.5){\makebox(0,0){$r_V$}}  
\put(20,7){\makebox(0,0){$V^* \so(\gotH \so V)$}}
    \put(20,2){\vector(0,1){4}}
 \put(20.5,4){\makebox(0,0)[l]{$\id_{V^*}\so t_V$}}      
\put(5,2){\vector(3,2){6}}
 \put(7.5,3){\makebox(0,0)[l]{$Q_{\gotH V}$}}    
\put(16,3){\makebox(0,0)[r]{$P_{V \gotH}$}}    
\put(12,6){\vector(3,-2){6}} 
  \end{picture}
\end{equation}
\begin{equation}
  \label{253.17}
  \begin{picture}(25,8)
   \put(11.5,1){\makebox(0,0){$\gotH \so V$}}
  \put(11.5,6){\vector(0,-1){4}}
   \put(12,4){\makebox(0,0)[l]{$t_V$}}    
\put(11.5,7){\makebox(0,0){$V\so\gotH$}} 
\put(1,7){\makebox(0,0){$\big((V \so\gotH) \so V^*\big)\so V$}}
\put(22,7){\makebox(0,0){$V \so\big( V^* \so(\gotH \so V)\big)$}}
  \put(10,7){\vector(-1,0){5}}
 \put(8,7.5){\makebox(0,0){$P_{(V\so\gotH)V^*}$}}
    \put(13,7){\vector(1,0){5}}  
 \put(15,7.5){\makebox(0,0){$\id_V \so r_V$}}    
\put(20,6){\vector(-3,-2){6}}
 \put(17,3){\makebox(0,0)[l]{$Q_{V^*(\gotH\so V)}$}}    
\put(3,6){\vector(3,-2){6}} 
 \put(6,3){\makebox(0,0)[r]{$l_V\so \id_V$}}    
  \end{picture}
\end{equation}

\subsection{Proof of the main theorem}
We are finally in the position to prove Theorem II as stated in the
introduction to Part II. 
  The uniqueness of $\M_1$ (up to equivalence) in part 2 of Theorem II
   follows by standard arguments. Indeed, suppose
  $\M\subset \tilde{\M}_1$ is another extension satisfying the same
  properties for $\tilde{\Gamma}: \dG\pfeil \tilde{\M}_1$. Choosing
  $\A = \tilde{\M}_1$ and $\gamma : \M \pfeil \A$ the inclusion map we
  have $\id : \tilde{\M}_1 \pfeil \tilde{\M}_1 \equiv \A$ as an extension of
  $\gamma$. Hence, putting 
  \begin{equation*}
    \GGt := e_\mu \tp \tilde{\Gamma}(e^\mu) \, \in \G\tp \tilde{\M}_1
  \end{equation*}
we conclude from the universality property of $\tilde{\M}_1$ that
$\GGt$ is a normal, coherent $\la\rho$--intertwiner in $\tilde{\M}_1$. Since
$\M_1$ also solves the universality property, there exists an algebra
map 
\begin{equation*}
  f \equiv \gamma_{\tilde{\Gamma}} : \M_1 \pfeil \tilde{\M}_1
\end{equation*}
restricting to the identity on $\M$ and satisfying $f \circ \Gamma =
\tilde{\Gamma}$. Since $\M_1$ is algebraically generated by $\M$ and
$\Gamma(\dG)$, $f$ is uniquely fixed by these
conditions. Interchanging the role of $\M_1$ and $\tilde{\M}_1$ we
also have a map $\tilde{f}: \tilde{\M}_1 \pfeil \M_1$ and clearly
$\tilde{f} = f^{-1}$.

We now prove the existence of $\M_1$ by choosing $\M_1 : = \M\relilr
\dG \equiv \M \relidr
\dG$ and $\Gamma : \dG \pfeil \M_1$, 
\begin{equation}
 \label{p1}
  \Gamma(\vi) : = \Big((q^2_\la)_{(0)} \relidr 
    \big( \Si(q^1_\la)\arr \vi \arl (q^2_\la)_{(1)}\big)\Big)
\end{equation}
where $q_\la \in \G\tp\M$ is given by \eqref{2.5.15'} and where for
$m\in\M$ we
write $\rho(m) \equiv \mo\tp\me$. Here $(\delta_r,\Psi_r)$ is the
two--sided $\G$--coaction constructed in Proposition \ref{prop20}. 
Putting $\RR_{\delta_r} : = e_\mu \tp (\eM
\reli e^\mu) \in \G\tp \M_1$ and comparing with the notation
\eqref{252.1} we conclude 
\begin{equation}
  \label{p2}
  \Gamma(\vi) = (\vi\tp \id_{\M_1})(\GG)
\end{equation}
where $\GG \in \G\tp\M_1$ is given by 
\begin{equation}
 \label{p3}
  \GG = (\idG\tp\rho)(q_\la) \succ \RR_{\delta_r} 
\end{equation}
Hence, by part 1 of Proposition \ref{prop2.15} $\GG$ is a normal
coherent $\la\rho$--intertwiner in $\M_1$. Moreover, putting $p \equiv p_\la
\in \G\tp\M$ given in Eq. \eqref{2.5.14'} we get for $m\in\M$ and
$\vi\in \dG$ 
\begin{align}
  \mu_R (m\tp \vi) :&= (m\reli \edG) \, \Gamma (\vi_{(1)}) \,
  (\vi_{(2)}\tp\id)(p) \notag \\
                   &=(m \reli \edG) \, (\vi\tp \id)(\GG \, p)\notag \\
                  &= (m\reli \edG) \, (\eM \reli \vi)\notag \\
                     & \equiv (m\reli \vi)
 \label{p4}
\end{align}
by Proposition \ref{prop2.15} 1.), and therefore $\mu_R: \M\tp \dG
  \pfeil \M_1$ becomes the identity map. This also shows that $\M_1$
  is algebraically generated by $\M \equiv (\M\reli \edG)$ and
  $\Gamma(\dG)$.

To prove that $\M_1$ solves the universality property we recall from
Corollary \ref{cor2.8} that there is a one--to--one correspondence
between extensions $\hat{\gamma} : \M_1 \pfeil \A$ of algebra maps
$\gamma : \M \pfeil \A$ and normal, coherent right
$\delta_r$--implementers $\RR_\A \in \G\tp \A$, such that 
\begin{equation*}
  (\idG\tp \hat{\gamma})(\RR_{\delta_r}) = \RR_\A.
\end{equation*}
Thus any such extension satisfies by \eqref{p3}
\begin{equation*}
  (\idG\tp\hat{\gamma})(\GG) = \TT_\A : = (\idG\tp\rho)(q_\la) \succ \RR_\A
\end{equation*}
By Proposition \ref{prop2.15} $\TT_\A \in \G\tp\A$ is a normal,
coherent $\la\rho$--intertwiner and the correspondence $\TT_\A
\leftrightarrow \RR_\A$ is one--to--one. Thus $\M_1$ solves the
universal property specified in part 1.\ of Theorem~II. 

We are left with showing that with $q_\rho \in \M\tp\G$  given by
\eqref{2.5.17'} also 
$\mu_L : \dG\tp\M \pfeil \M_1$ given in \eqref{2.03} provides a
linear isomorphism, which under the identification $\dG\tp\M
=\dG\relidl \M$ in fact becomes an algebra map. 
To this end we use the twist equivalence \\ 
$f : \M_1 \equiv \M\relidr
\dG \to \M\relidl \dG$ of Prop.~\ref{prop2.10} given by 
\begin{equation*}
  (\idG\tp f )(\RR_{\delta_r}) = \phi_{\la\rho} \succ \RR_{\delta_l},
\end{equation*}
see Eq.~\eqref{6.45c}. Putting $\LL_{\delta_l}: =e_\mu\tp (e^\mu \reli
\eM) \in \G\tp (\dG\relidl\M)$ this gives
\begin{align*}
  (\idG\tp f \circ\mu_L)(\LL_{\delta_l}) &\equiv (\idG\tp
  f )(q^{op}_\rho\,\GG) \equiv (\idG\tp f )\big((\eG\tp
  q_\rho)\succ \GG\big)\\
 & = \big[ (\eG\tp q_\rho)\, (\idG\tp\rho)(q_\la)\,\phi_{\la\rho}\big]
  \succ \RR_{\delta_l} \\
  &= q_{\delta_l} \succ \RR_{\delta_l}
\end{align*}
where we have used \eqref{p3} and \eqref{5.57}. Hence, by
Theorem~\ref{thm2.5}(iii) 
$f \circ \mu_L : \dG\relidl \M \to \M\relidl\dG$ provides an isomorphism.
 Since $f$ is invertible this concludes the proof of
part 3 of Theorem II. \qed

\vspace{0.5cm}
We think it to be gratifying that in view of this proof and the
preceding results all the different appearances of left and right
diagonal crossed products described by associated left/right
$\delta_{l/r}$--implementers may be replaced by the one universal
$\la\rho$--intertwiner $\GG$, which is independent of any left--right
conventions. This is also our main motivation for formulating Theorem
II in this way. In particular, from now on we may dispense with all
the nasty calculations involving $\delta_{l/r}$--implementers and just
work with the much more convenient generating relations \eqref{2.01},
\eqref{2.02}.


\section{Examples and applications}

\subsection{The quantum double $\D(\G)$}
In view of the identification of the quantum double $\D(\G)$ of an
ordinary Hopf algebra $\G$ with the diagonal crossed product $\G\reli
\dG$ in \eqref{1.32'} we propose the following 
\begin{definition}
{\rm
  \label{def2.61}
Let $(\G,\cop,\ep,\phi)$ be a quasi--Hopf algebra. 
  The diagonal crossed product $\D(\G): = \dG\relilr \G \cong
  \G\relilr \dG$ associated with the
  quasi--commuting pair $(\la=\rho=\cop,\,\, \phi_\la =
  \phi_\rho=\phi_{\la\rho} = \phi)$  of $\G$--coactions on $\M
  \equiv\G$ is called the {\it quantum 
  double of $\G$}.
}
\end{definition}
Following the notations of [N1], the universal 
$\la\rho$--intertwiner of the quantum double will be denoted by
$\DD \equiv \GG_{\D(\G)}\in \G\tp\D(\G)$. Hence it obeys the relations
$(\ep\tp\id)(\DD) = \e_{\D(\G)}$ and 
\begin{align}
  \label{2.61}
  \DD\, \cop(a) &= \cop^{op}(a)\, \DD, \quad \forall a \in \G \\
  \label{2.62}
   \Fi[312]{}\,\DD[13]\, \Fii[132]{}\, \DD[23]\, \phi &= (\cop\tp\id)(\DD)
\end{align}
where we have suppressed the embedding $\G \hookrightarrow \D(\G)$. 
Eq.~\eqref{2.61} motivates to call $\DD$ the {\it universal flip
  operator} for $\cop$.
Theorem II implies that the category of
representations of our quantum double $\D(\G)$ is nothing else but the
so--called ``double category of $\G$--modules'', which was used by
S. Majid [M2] to define the quantum double of a quasi--Hopf algebra using a
Tannaka--Krein like reconstruction procedure. 
Indeed, in our terminology this is precisely the subcategory of {\it
  $\la\rho$--coherent representations} $(\gotH,\gamma_\gotH,
\TT_\gotH)$ in $\Rep \G$ (where $\la =\rho =\cop$), see Definition 
\ref{def2.11'}, with morphisms given by the $\G$--linear maps $\tau :
\gotH \to \gotH'$ satisfying $(\id\tp\tau)(\TT_\gotH) = \TT_{\gotH'}$,
where $\TT_\gotH$ and $\TT_{\gotH'}$ are the normal coherent
$\la\rho$--intertwiners in $\G\tp \End_\Co (\gotH)$ and $\G\tp
\End_\Co (\gotH')$, respectively. A more detailed account of this will
be given in [HN1], where we will also 
show, that $\D(\G)$ is in fact a quasitriangular quasi--Hopf algebra. 

Here we restrict ourselves with proving the quasi--bialgebra
structure on $\D(\G)$. Analogously as in Proposition  \ref{Prop 1.9}
this will also guarantee that every diagonal crossed product\\ 
$\M_1 = \M\relilr \dG$ naturally admits a quasi--commuting pair
$(\la_D,\rho_D,\Fi{\la_D},\Fi{\rho_D},\Fi{\la_D \rho_D})$ of coactions
of $\D(\G)$ on $\M_1$. We begin with constructing $\la_D: \M_1 \to
\D(\G) \tp \M_1$ and $\rho_D: \M_1 \to \M_1 \tp \D(\G)$ as algebra
maps extending the left and right coactions $\la: \M_1 \supset \M \to
\G\tp \M \subset \D(\G)\tp \M_1$ and $\rho: \M_1 \supset \M \to
\M\tp\G \subset \M_1 \tp \D(\G)$, respectively (see
\eqref{1.32c},\eqref{1.32c'}).  
\begin{lemma}
  \label{lem2.62}
  Let $(\la,\rho,\phi_\la,\phi_\rho,\phi_{\la\rho})$ be a
  quasi--commuting pair of $\G$--coactions on $\M$ and let $\M_1
  \equiv \dG\relilr \M$ be the associated diagonal crossed product
  with universal $\la\rho$--intertwiner $\GG \in \G\tp \M_1$. Then
  there exist uniquely determined algebra maps $\la_D: \M_1 \to
  \D(\G) \tp \M_1$ and $\rho_D: \M_1 \to \M_1\tp \D(\G)$ satisfying  
  (we suppress all embeddings $\M \hookrightarrow \M_1$ and
  $\G\hookrightarrow \D(\G)$) 
  \begin{align}
    \label{2.63}
  \la_D(m) &= \la(m), \quad \forall m \in \M \subset \M_1 \\
    \label{2.64}
  (\id\tp \la_D)(\GG) & = \Fii[231]{\la\rho}\, \GG[13]\,
  \Fi[213]{\la}\, \DD[12] \, \phi^{-1}_\la \, \in \G\tp \D(\G)\tp \M_1
  \\
  \label{2.65}
\rho_D(m) & = \rho(m), \quad \forall m \in \M\subset \M_1 \\
 \label{2.66}
(\id\tp \rho_D)(\GG) & = \Fii[231]{\rho}\, \DD[13]\,
  \Fi[213]{\rho}\, \GG[12] \, \phi^{-1}_{\la\rho} \, \in \G\tp \D(\G)\tp \M_1
\end{align}
\end{lemma}
\begin{proof}
Viewing the left $\G$--coaction $\la : \M \pfeil \G\tp\M$ as a map
$\la: \M\pfeil \D(\G) \tp \M_1$, 
Theorem~II states that $\la_D$ is a unital algebra map extending
$\la$ if and only if $\TT_D : = (\id\tp\la_D)(\GG) \in \G \tp(\D(\G)
\tp \M_1)$ is a normal coherent $\la\rho$--intertwiner. Now normality of
$\TT_D$ follows from the normality of $\GG$. To 
prove that $\TT_D$ is a $\la\rho$--intertwiner we compute for all $m\in \M$ 
\begin{align*}
  \TT_D \, (\idG\tp \la_D)(\la(m)) & =
     \Fii[231]{\la\rho}\, \GG[13]\, \Fi[213]{\la} \,\DD[12]\,
     \phi^{-1}_\la \, (\idG\tp\la_D)(\la(m)) \\
    & = \big[(\la_D \tp \idG)(\rho(m))\big]^{231} \, 
             \Fii[231]{\la\rho}\, \GG[13]\, \Fi[213]{\la} \,\DD[12]\,
     \phi^{-1}_\la \\
    & = (\idG\tp \la_D)(\rho^{op}(m)) \, \TT_D
\end{align*}
where both sides are viewed as elements in $\G\tp \D(\G) \tp
\M_1$. Here we have used the intertwining properties of $\GG$ and $\DD$ and
of the three reassociators. 

To show that $\TT_D$ also satisfies the coherence condition,
i.e. Eq.~\eqref{2.02}, we compute in $\G\tp\G\tp\D(\G)\tp \M_1$ - again
suppressing all embeddings
\begin{align*}
  (\cop\tp\id)(\TT_D) &= \big[(\id\tp\id\tp\cop)(\phi_{\la\rho}^{-1})
  \,[\eG\tp\phi_\rho]\big]^{3412}\notag \\  
    &\quad\quad \GG[14] \, \Fii[142]{\la\rho}  \,
   \GG[24]  
   \big[[\eG\tp\phi_\la]\,(\id\tp\cop\tp\id)(\phi_\la)
    [\phi\tp\eM]\big]^{3124} \DD[13] \,\Fii[132]{}\,\DD[23] \notag \\
  &\quad\quad  [\phi\tp\eM]
  (\cop\tp\id\tp\id)(\phi^{-1}_\la) \notag \\
  &=  
  \big[ (\la\tp\id\tp\id)(\phi_\rho) [\phi^{-1}_{\la\rho}\tp\eG]
  (\id\tp\rho\tp\id)(\phi^{-1}_{\la\rho})\big]^{3412} \notag \\ 
     &\quad\quad \GG[14] \, \Fii[142]{\la\rho}  \,
   \GG[24]  
   \big[ (\id\tp\id\tp\la)(\phi_\la)(\cop\tp\id\tp\id)(\phi_\la)\big]^{3124}
    \DD[13] \,\Fii[132]{}\,\DD[23] \notag \\
   &\quad\quad (\id\tp\cop\tp\id)(\phi^{-1}_\la)
   [\eG\tp\phi^{-1}_\la](\id\tp\id\tp\la)(\phi_\la) \notag \\
   &= 
\big[ (\la\tp\id\tp\id)(\phi_\rho)
   [\phi^{-1}_{\la\rho}\tp\eG]\big]^{3412} \, \GG[14] \notag \\
    &\quad\quad
  \big[(\id\tp\la\tp\id)(\phi^{-1}_{\la\rho}) 
      \, [\eG \tp
      \phi^{-1}_{\la\rho}]\,(\id\tp\id\tp\rho)(\phi_\la)\big]^{3142} 
   \, \GG[24]\, \DD[13]\notag\\  
    &\quad\quad  
   \big[
 (\cop\tp\id\tp\id)(\phi_\la)[\phi^{-1}\tp\eM]\,
 (\id\tp\cop\tp\id)(\phi^{-1}_\la) \big]^{1324} 
   \notag \\
   &\quad\quad 
   \DD[23]\,[\eG\tp\phi^{-1}_\la](\id\tp\id\tp\la)(\phi_\la) \notag \\
   & = 
 \big[ (\la\tp\id\tp\id)(\phi_\rho)
   [\phi^{-1}_{\la\rho}\tp\eG]\big]^{3412} \, \GG[14]\,
    \Fi[314]{\la} \, \DD[13]\notag \\
    &\quad\quad \big[(\cop\tp\id\tp\id)\Fii{\la\rho}\,
    (\id\tp\id\tp\rho)(\phi^{-1}_\la)\big]^{1324}\, \GG[24] \\
    &\quad\quad \phi_\la^{324}\, \DD[23] 
\,[\eG\tp\phi^{-1}_\la](\id\tp\id\tp\la)(\phi_\la) \notag \\
  &= (\id\tp\id\tp\la_D)\Big( \Fi[312]{\rho} \, \GG[13]\, \Fii[132]{\la\rho}
  \, \GG[23] \, \Fi{\la} \Big)
\end{align*}
Here we have used several pentagon identities for the reassociators
involved and the intertwining and coherence properties of $\GG$ and
$\DD$. In the first equality we used \eqref{2.02} for $\GG$ and $\DD$,
and in the second the pentagons \eqref{eq27} and \eqref{eq22}. For the
third equality we used the intertwining properties of $\DD$ and $\GG$ to
move two more reassociators between $\DD[13]$ and $\DD[23]$ and to more
between $\GG[14]$ and $\GG[24]$. To arrive at the fourth equality  we
commuted $\DD[13]$ and 
$\GG[24]$ and used the pentagons \eqref{eq26} and \eqref{eq22} and
then again the intertwining properties of $\DD$ and $\GG$ to bring two
reassociators back between $\DD[13]$ and $\GG[24]$. The last equality
holds by  \eqref{2.63}, \eqref{2.64}.
Thus we have
shown that $\TT_D$ is coherent and therefore the definitions
\eqref{2.63}, \eqref{2.64}  uniquely define a
unital algebra map $\la_D$ extending $\la$. Similarly one shows that $\rho_D$
defines a unital algebra map $\rho_D : \M_1 \pfeil \M_1 \tp \D(\G)$
extending $\rho$.
\end{proof}
Choosing in Lemma \ref{lem2.62} also $\M = \G$ (i.e. $\M_1 = \D(\G)$)
we arrive at the following
\begin{theorem}
  \label{thm2.63}
Let $(\G,\cop,\ep,\phi)$ be a quasi--Hopf algebra,
denote $i_D: \G \hookrightarrow \D(\G)$ the canonical embedding and 
    let $\DD \in \G\tp\D(\G)$ be the universal
    flip operator.
  \begin{itemize}
  \item[(i)]  Then
    $(\D(\G),\cop_D,\ep_D, \phi_D)$ is a quasi--bialgebra, where
    \begin{align}
      \label{2.67}
\phi_D &: = ( i_D\tp i_D\tp i_D)(\phi) \\
      \label{2.68}
 \ep_D( i_D(a)) &: = \ep(a), \quad (\id\tp\ep_D)(\DD) := \e_{\D(\G)} \\
   \label{2.69}
 \cop_D (i_D(a)) &: = (i_D\tp i_D)(\cop(a)), \quad \forall a \in \G
 \\
    \label{2.610}
  (i_D\tp\cop_D)(\DD) &:= \Fii[231]{D}\, \DD[13]\, \Fi[213]{D}\, \DD[12]
  \, \phi_D^{-1}
    \end{align}
\item[(ii)]
Under the setting of Lemma \ref{lem2.62} denote $i_{\M_1}: \M
\hookrightarrow \M_1$ the embedding and define 
 \begin{align*}
   \phi_{\la_D} &: = (i_D \tp i_D \tp i_{\M_1})(\phi_\la) \in
 \D(\G)\tp\D(\G)\tp\M_1 \\
  \phi_{\rho_D}& := ( i_{\M_1} \tp i_D \tp i_D )(\phi_\rho)
   \in  \M_1\tp\D(\G)\tp \D(\G) \\
  \phi_{\la_D \rho_D} &:=(i_D\tp i_{\M_1}\tp i_D)(\phi_{\la\rho})
  \in  \D(\G)\tp\M_1 \tp \D(\G) 
 \end{align*}
 Then
 $(\la_D,\rho_D,\phi_{\la_D},\phi_{\rho_D},\phi_{\la_{D} \rho_D})$
 provides a  quasi--commuting pair of $\D(\G)$--coactions on $\M_1
 \equiv \dG\relilr \M$.
  \end{itemize}
\end{theorem}
\begin{proof}
  Setting $\M := \G$ and $\la = \cop$ in Lemma \ref{lem2.62} implies
  that $\cop_D$ is a unital algebra morphism. The property of $\ep_D$
  being a counit for $\cop_D$ follows directly from the fact that
  $(\id\tp\ep\tp\id)(\phi) = \eG\tp\eG$. 
To show that $\cop_D$ is quasi--coassociative one computes that
\begin{equation*}
  [\eG\tp \phi_D] \cdot (\id\tp\cop_D \tp \id)\Big((\id\tp\cop_D)(\DD)
  \Big)
  = (\id\tp\id\tp\cop_D)\Big((\id\tp\cop_D)(\DD)\Big) \cdot [\eG\tp\phi_D],
\end{equation*}
where one has to use \eqref{2.610}, the pentagon equation for $\phi$
and the intertwiner property \eqref{2.61} of $\DD$ similarly as in the
proof of Lemma 
\ref{lem2.62}. Thus $\cop_D$ is quasi--coassociative and this concludes
the proof of part (i).

Part (ii) is shown by direct calculation using the intertwiner properties of
$\GG$ and $\DD$ and several pentagon identities for the reassociators
involved. The details are left to the reader.
\end{proof}
Note that viewed in $\D(\G)\tp \D(\G)$ and $\D(\G)^{\tp^3}$,
respectively, the relations \eqref{2.61}, \eqref{2.62} and
\eqref{2.610} 
are precisely the defining properties of a quasitriangular $R$--matrix
[Dr2]. Thus $R_D : = (i_D \tp \id)(\DD)$ is an $R$--matrix for
$\D(\G)$. It is shown in [HN1] that there also exists an antipode $S_D$
for $\D(\G)$ extending the antipode of $\G$, thus making the quantum
double $\D(\G)$ into a quasitriangular quasi--Hopf algebra.


\subsection{Two--sided crossed products}
As in the associative case, a simple recipee to produce two--sided
$\G$--comodule algebras $(\M,\delta)$ is by tensoring a right
$\G$--module algebra $(\A,\rho_\A)$ and a left $\G$--comodule algebra
$(\B,\la_\B)$, i.e. by setting $\M = \A\tp\B$ and 
\begin{equation*}
  \delta(A\tp B) : = B_{(-1)} \tp (A_{(0)}\tp B_{(0)})\tp A_{(1)}
\end{equation*}
as in Eq.~\eqref{1.36a}. 
Obviously $\delta = (\la\tp\id)\circ \rho = (\id\tp\rho)\circ \la$,
where $(\la,\phi_\la)$ and 
$(\rho,\phi_\rho)$ are
the trivially extended left and right coactions $(\la_\B,\phi_{\la_\B})$ and
$(\rho_\A,\phi_{\rho_\A})$, respectively. Hence
$(\la,\rho,\phi_\la,\phi_\rho,\phi_{\la\rho}= \eG\tp\e_{\M}\tp \eG)$
is a {\it strictly commuting} pair of coactions. 
In the terminology of Proposition \ref{prop20} we have $\delta =
\delta_r = \delta_l$, whereas $\Psi = \Psi_r = \Psi_l$ is given by
\begin{equation*}
  \Psi = [\eG\tp\eG\tp\phi^{-1}_\rho] \, [\phi_\la \tp
  \eG\tp\eG]. 
\end{equation*}
According to Theorem II the diagonal crossed
  product $\M_1 = (\A\tp \B) \relilr \dG$ is generated by $\{A,B,
  \Gamma(\vi)\mid A\in \A,\, B\in \B,\,\vi\in\dG\}$ satisfying the defining
  relations 
  \begin{align}
   \label{2.612'}
   A\,B &= B\,A \\
    \label{2.612}
   [\eG\tp B]\, \GG &= \GG \, \la(B) \\
    \label{2.613}
   \rho^{op} (A)\, \GG &= \GG \, [\eG\tp A]\\
    \label{2.614}
   (\cop\tp\id)(\GG) &= \Fi[312]{\rho}\, \GG[13]\,\GG[23]\, \phi_\la
  \end{align}
where $\GG = e_\mu \tp \Gamma(e^{\mu})$ is the universal
$\la\rho$--intertwiner. 

The next Proposition is an analogue of Proposition \ref{Prop 1.12}
saying that the diagonal crossed product $(\A\tp\B)\relilr \dG$ may be
realized as a {\it two--sided crossed product} $\A\crosr\dG\crosl\B$.
Note that the isomorphism $\mu$ in Eq.~\eqref{2.616} below
 is different from the isomorphisms $\mu_R$ and $\mu_L$
constructed in Theorem II.

\begin{proposition}
  \label{prop2.64}
Let $\re:\hat\G\o\A\to\A$ and $\li:\B\o\hat\G\to\B$ be the left and 
the right action corresponding to a right $\G$--coaction
$(\rho,\phi_\rho)$ on $\A$ 
and a left $\G$--coaction $(\la,\phi_\la)$ on $\B$,
respectively. Extend $\la$ and $\rho$ trivially to $\A\tp\B$ and let
$\M_1 : = (\A\tp\B) \relilr \dG \equiv (\A\tp\B)\relidr \dG$,
$\delta_r : = (\id\tp\rho)\circ \la$, be the diagonal crossed product
with universal $\la\rho$--intertwiner $\GG \in \G\tp \M_1$.
\begin{itemize}
\item[(i)] There is a linear bijection $\mu: \A\tp\dG\tp \B \pfeil
  \M_1$ given by 
\begin{equation}
\label{2.616}
\mu(A\tp \vi\tp B) =  A\, \Gamma(\vi)\, B 
\end{equation}
where we have suppressed the embeddings $\A\hookrightarrow \M_1$ and
$\B\hookrightarrow \M_1$.
\item[(ii)]
Denote the induced algebra structure on $\A\tp\dG\tp\B$ by 
$\A\crosr\hat\G \crosl \B \equiv \mu^{-1}(\M_1)$. Then we get the
following multiplication structure with unit
$\one_\A\>cros\hat\one\<cros\one_\B$ on $\A\crosr\hat\G \crosl \B$ 
\begin{align}
&\big(A\>cros\vi\<cros B\big)\,\big(A'\>cros\psi\<cros B'\big) \notag\\
&\quad\quad\quad \label{2.615}  =
A(\vi\1\re A')\Fib[1]{\rho} \>cros
[\Fib[1]{\la}\arr\vi\2\arl\Fib[2]{\rho}]\,[\Fib[2]{\la}\arr
\psi\1\arl \Fib[3]{\rho}]  
\<cros\Fib[3]{\la}(B\li\psi\2)B'
\end{align}
\end{itemize}
\end{proposition}
\begin{proof}
  Let $p_\la \in \G\tp\B \equiv \G\tp (\e_\A \tp \B)$ be given by
  Eq.~\eqref{2.5.14'}. Then, using \eqref{2.612}
  \begin{align*}
    \mu\Big(A\tp\vi_{(1)} \tp (\vi_{(2)} \tp \id_\B)\big(\la(B)
    p_\la\big)\Big) &=
   A \,(\vi\tp \id_{\M_1})\big(\GG\, \la(B) \, p_\la \big) \\
  & = A\,B\, (\vi\tp \id_{\M_1})(\GG\, p_\la) \\
  & = \mu_R (A\tp B\tp \vi) 
  \end{align*}
where $\mu_R : (\A\tp\B)\tp \dG \to \M_1$ is the linear bijection
constructed in part 3.~of Theorem II, see also \eqref{p4}. Hence $\mu$
is surjective. Conversely, let $\RR : = \GG\, p_\la \in \G\tp \M_1$
then by Proposition~\ref{prop2.15} $\RR$ is a right
$\delta_r$--implementer and 
\begin{equation*}
  \GG = (\idG\tp\rho)(q_\la) \succ \RR  \equiv
  [\eG\tp q^2_\la]\, \RR \, [\Si(q^1_\la) \tp \e_{\M_1}]
\end{equation*}
where $q_\la \in \G \tp \B$ if given by \eqref{2.5.15'} and where we
have used that $\rho$ is trivial on $\B$. Hence we get
for all $B\in \B$ 
\begin{align}
  \notag
  \GG \, [\eG\tp B] &= [\eG \tp q^2_\la] \, [\delta_r (B) \succ
  \RR]\,[\Si(q^1_\la)\tp \e_{\M_1}] \\
   \label{2.617}
  & = [\eG \tp q^2_\la B_{(0)}]\, \RR \,
    [\Si(q^1_\la B_{(-1)}) \tp \e_{\M_1}]
\end{align}
where $B_{(-1)}\tp B_{(0)} = \la(B)$ and where we have used $\delta_r
(B) = B_{(-1)}\tp (\eA \tp B_{(0)}) \tp \eG$. Eq.~\eqref{2.617}
implies for all $A\in\A,\vi \in \dG$ and $B \in \B$ 
\begin{equation*}
  A\, \Gamma(\vi)\, B = \mu_R \Big( A \tp
  \big(\hat{S}^{-1}(\vi_{(2)})\tp\id_\B\big)\big(q_\la \la(B)\big)\tp
  \vi_{(1)}\Big)  
\end{equation*}
and therefore the injectivity of $\mu_R$ implies the injectivity of
$\mu$.

This proves part (i). Part (ii) follows since one straightforward
checks that the multiplication rule \eqref{2.615} is equivalent to the
defining relations \eqref{2.612}-\eqref{2.614}.
\end{proof}

\subsection{Hopf spin chains and lattice current algebras}
In this subsection we describe how the Hopf algebraic quantum chains
considered in [NSz] generalize to quasi--Hopf algebras $\G$ (the case
of weak quasi--Hopf algebras like semisimple quotients of quantum
groups at roots of unity requires only minor changes as sketched in
Part~III). To this end we use the two--sided crossed product theory of
Section~11.2 to generalize the constructions \eqref{1.43},\eqref{1.44}
and \eqref{1.46'}.

First we show, that analogously as in \eqref{1.46'} the two--sided
crossed product construction given in Proposition~\ref{prop2.64} may
be iterated if 
one of the two algebras $\A$ and $\B$ 
 admits a quasi--commuting pair of coactions. 
\begin{proposition}
\label{prop2.65}
  Let $(\A,\rho_\A,\phi_{\rho_\A})$, $(\C,\la_\C,\phi_{\la_\C})$ 
and $(\B,\rho_\B,\la_\B,\phi_{\la_\B},\phi_{\rho_\B},\phi_{\la_\B,\rho_\B})$
 be a right, a left and  a two--sided comodule
  algebra, respectively, and denote the universal
  $\la\rho$--intertwiners
  \begin{align}
   \label{GAB}
    \GG_{\A\B} &:= e_\mu \tp (\eA\>cros e^\mu \<cros \eB) \in
    \A\crosrA \dG \croslB \B \\
    \label{GBC}
\GG_{\B\C} &:= e_\mu \tp (\eB\>cros e^\mu \<cros \eC) \in
    \B\crosrB \dG \croslC \C
  \end{align}
Then
\begin{itemize}
\item[(i)]
 $\A\crosrA \dG \croslB \B$ admits a right $\G$--coaction 
$(\tilde{\rho},\phi_{\tilde{\rho}})$ given by 
  $\phi_{\rhot} : = \phi_{\rho_\B}$, $\rhot |_{(\A\tp \B)} : = \idA
 \tp \rho_\B$ and
 \begin{equation}
   \label{2.622}
   (\idG\tp \rhot)(\GG_{\A\B}) : = (\GG_{\A\B} \tp \eG) \,
   \phi^{-1}_{\la_B\rho_B} 
 \end{equation}
\item[(ii)]
$\B\crosrB \dG \croslC \C$ admits a left $\G$--coaction
$(\lat,\phi_\lat)$ given by $\phi_\lat : = \phi_{\la_\B}$, $\lat
|_{(\B\tp\C)} : = \la_B \tp \id_\C$ and 
\begin{equation}
  \label{2.623}
  (\idG\tp \lat)(\GG_{\B\C}) : = \Fii[231]{\la_\B\rho_\B} \, \GG[13]_{\B\C}
\end{equation}
\item[(iii)]
We have an algebra isomorphism 
\begin{equation}
  \label{2.624}
 (\A\crosrA \dG \croslB \B) \crosrt \dG \croslC \C \equiv
  \A \crosrA\dG \croslt (\B \crosrB\dG \croslC \C)
\end{equation}
given by the trivial identification.
\end{itemize}
\end{proposition}
\begin{proof}
  (i) To show that \eqref{2.622} provides a well defined algebra map 
$ \rhot:   \A\crosrA \dG \croslB \B \to (\A\crosrA \dG \croslB \B)\tp
\G$ extending $\id_\A \tp \rho_\B$ we have to check that the relations
\eqref{2.612}-\eqref{2.614} are respected. To this end we put 
$\TT_{\A\B} : = (\GG_{\A\B} \tp \eG)\, \phi^{-1}_{\la_\B \rho_\B}$ and
compute for $B\in \B$
 \begin{align*}
 [\eG\tp \rhot(B)]\, \TT_{\A\B} 
     &= 
      (\GG_{\A\B} \tp \eG)\, (\la_\B\tp\id)(\rho_\B(B))\,
    \phi^{-1}_{\la_B\rho_B} \\
   & = \TT_{\A\B}\,  (\id_G\tp \rhot)(\la_\B(B))  
  \end{align*}
which is the relation \eqref{2.612}. Trivially one also has (since
$\phi_{\la_B\rho_B} \in \G\tp (\eA\tp\B)\tp \G$)
\begin{equation*}
  \TT_{\A\B} \, [\eG\tp A\tp \eG] = [\rho^{op}(A) \tp \eG] \, \TT_{\A\B}
\end{equation*}
verifying \eqref{2.613}. Finally, the coherence condition
\eqref{2.614} is respected, since in 
\\
$\G\tp\G\tp (\A\>cros \dG \<cros
\B)\tp \G$ we have 
\begin{align*}
&(\cop\tp\id\tp\id)(\TT_{\A\B})  =  (\cop\tp \id)(\GG_{\A\B}) ^{123} \,
  (\cop\tp\id\tp\id)\Fii{\la_\B\rho_\B} \\
  &\quad\quad=(\cop \tp \id)(\GG_{\A\B}) ^{123} \, 
    [\phi^{-1}_{\la_\B} \tp
  \eG]\,(\id\tp\la_\B\tp\id)\Fii{\la_\B\rho_\B} \, [\eG\tp 
  \phi^{-1}_{\la_\B\rho_\B}]\, (\id\tp\id\tp\rho_\B)(\phi_{\la_\B}) \\
&\quad\quad= \Fi[312]{\rho_\A} \, \GG[13]_{\A\B} \,\GG[23]_{\A\B} \,
  (\id\tp\la_\B\tp\id)\Fii{\la_\B\rho_\B} \, [\eG\tp
  \phi^{-1}_{\la_\B\rho_\B}]\, (\id\tp\id\tp\rho_\B)(\phi_{\la_\B})\\
&\quad\quad=
  \Fi[312]{\rho_\A} \, \GG[13]_{\A\B} \, \Fii[134]{\la_\B\rho_\B} \, \GG[23]_{\A\B}
  \Fii[234]{\la\rho} \, (\id\tp\id\tp \rho_B)(\phi_{\la_\B}) \\
&\quad\quad=
  (\id\tp\id\tp\rhot)(\Fi[312]{\rho_\A})\,\TT[13]_{\A\B}\, \TT[23]_{\A\B}\,
  (\id\tp\id\tp\rhot)(\Fi{\la_\B})  
\end{align*}
Here we have used the pentagon
identity \eqref{eq26} in the second line, the coherence property
(\ref{2.614}) of $\GG_{\A\B}$ in the third line and finally the intertwining
property \eqref{2.612}.
Thus $\rhot$ provides a well defined algebra map, which is also unit
preserving since $(\ep\tp\id\tp\id)(\TT_{\A\B}) = (\eA\>cros\edG\<cros
\eB) \tp \eG$. Similarly one shows by a straight
forward calculation that the pair $(\rhot,\phi_\rhot)$ satisfies
\eqref{eq23}. Since  $\phi_\rhot = \phi_{\rho_\B}$, the pentagon
equation \eqref{eq24} and the counit equations \eqref{eq23b} and
\eqref{eq23c} are clearly
satisfied. This proves part (i). Part (ii) follows analogously. To
prove part (iii) we have to check that we may consistently identify
\begin{align}
 \label{*}
  \G\tp \Big[ (\A\>cros \dG \<cros \B) \>cros \dG \<cros \C \Big]\ni
   \GG_{\A\B} &
\stackrel{!}{\equiv}  \GG_{\A(\B\>cros \dG \<cros\C)} \in 
  \G\tp \Big[ \A\>cros \dG\<cros (\B\>cros\dG \<cros \C)\Big] \\
  \label{**}
  \G\tp \Big[ (\A\>cros \dG \<cros \B) \>cros \dG \<cros \C \Big]\ni
   \GG_{(\A\>cros\dG\<cros\B)\C} &\stackrel{!}{\equiv} \GG_{\B\C} \in 
  \G\tp \Big[ \A\>cros \dG\<cros (\B\>cros\dG \<cros \C)\Big]
\end{align}
For this the nontrivial commutation relations to be looked at are
according to \eqref{2.612},\eqref{2.613}
\begin{align*}
  \GG_{(\A\>cros\dG\<cros\B)\C} \, [\eG\tp X] &= \rhot^{op}(X)\,
    \GG_{(\A\>cros\dG\<cros\B)\C},\quad X \in \A\>cros\dG\<cros\B \\
 [\eG\tp Y]\, \GG_{\A(\B\>cros\dG\<cros\C)} &=
 \GG_{\A(\B\>cros\dG\<cros\C)}\, \lat(Y), \quad Y\in \B\>cros\dG\<cros\C
\end{align*}
The remaining cases being trivial, it is enough to consider $X\in
\Gamma_{\A\B} (\dG)$ and $Y \in \Gamma_{\B\C}(\dG)$ for which we get
\begin{align*}
  \GG[23]_{(\A\>cros\dG\<cros\B)\C}\, \GG[13]_{\A\B} &= 
   (\idG \tp \rhot)(\GG_{\A\B})^{132} \, \GG[23]_{\A(\B\>cros \dG
   \<cros\C)}
   = \GG[13]_{\A\B} \, \Fii[132]{\la_B\rho_B} \, 
   \GG[23]_{(\A\>cros\dG\<cros\B)\C}
    \\
  \GG[23]_{\B\C} \, \GG[13]_{\A(\B\>cros \dG \<cros\C)} &=
   \GG[13]_{\A(\B\>cros \dG \<cros\C)}\, (\idG\tp \lat)(\GG_{\B\C})^{213} 
   = \GG[13]_{\A(\B\>cros \dG \<cros\C)}\, \Fii[132]{\la_B\rho_B}\, 
\GG[23]_{\B\C}
\end{align*}
by the definitions \eqref{2.622} and \eqref{2.623}. This shows that
the identifications \eqref{*} and \eqref{**} are indeed consistent and
therefore proves part (iii)
\end{proof}

Note that due to part (iii) of the above proposition the notations
$\GG_{\A\B}$ and $\GG_{\B\C}$ as in 
\eqref{GAB},\eqref{GBC} are still well defined in iterated two--sided
crossed products and the commutation relations of ``neighboring''
$\la\rho$--intertwiners are given by    
\begin{equation}
  \label{2.626}
  \GG[13]_{\A\B} \, \phi^{-1}_{\la_\B\rho_\B} \, \GG[23]_{\B\C} = 
  \GG[23]_{\B\C} \, \GG[13]_{\A\B}
\end{equation}

Due to Proposition \ref{prop2.65} the definition of Hopf spin chains
as reviewed in Section 4.3 immediately generalizes to the
quasi--coassociative case. As in Section 4.3 we interpret even
integers as 
sites and odd integers as links of a one dimensional lattice and we
set $\A_{2i} \cong\G$, 
$\A_{2i+1} \cong \dG$, the latter just being a linear space. A local
net of associative algebras $\A_{n,m}$ is then constructed inductively
for all 
$n,m\in 2\mathbb{Z}, n \leq m$, by first putting 
\begin{equation*}
  \A_{2i,2i+2}:= \A_{2i}\>cros \A_{2i+1} \<cros \A_{2i+2} \cong
  \G\>cros \dG \<cros \G 
\end{equation*}
where $\G$ is equipped with its canonical two--sided comodule structure 
$(\la = \rho =\cop$, $\phi_\la = \phi_\rho = \phi_{\la\rho} = \phi)$.
Due to Proposition \ref{prop2.65} this procedure may be iterated as in
\eqref{1.43}, by setting 
\begin{equation*}
  \A_{2i,2j+2} := \A_{2i,2j} \>cros \dG \<cros \G \equiv 
                 \G\>cros \dG\<cros  \A_{2i-2,2j+2} 
\end{equation*}
where the last equality follows from \eqref{2.624} by iteration. 
More generally one has as in \eqref{1.44} for all $i\leq \mu \leq j-1$
\begin{equation}
  \label{2.627'}
  \A_{2i,2j} = \A_{2i,2\mu}\>cros \dG \<cros \A_{2\mu+2,2j}
\end{equation}
Defining the generating ``link operators''
  $\LL_{2i +1} : = \GG_{\A_{2i}\A_{2i+2}}$ as in \eqref{GAB},
 the commutation relation \eqref{2.626} implies
 \begin{equation*}
   \LL[13]_{2i -1}\, \phi^{-1}_{2i} \, \LL[23]_{2i +1} =
   \LL[23]_{2i+1}\, \LL[13]_{2i-1}, \quad \phi_{2i} : = (\id\tp\id\tp
   A_{2i})(\phi) .
 \end{equation*}
Writing $A_{2i +1} (\vi) =
(\vi\tp\id)(\LL_{2i+1})$ this is equivalent to 
\begin{equation}
  \label{2.627}
  A_{2i-1}(\vi) \, A_{2i+1}(\psi) = A_{2i+1}(\psi \arl
  \Fib[3]{}) \, A_{2i}(\Fib[2]{}) \, A_{2i-1}(\Fib[1]{}
  \arr \vi)   
\end{equation}
Thus link operators on neighboring links do not commute any more
in  contrast to the coassociative setting! But the algebras $\A_{2i
  -4,2i-2}$ and $\A_{2i+2,2i+4}$ still commute, which means that the
above construction still yields a local net of algebras.

Next we remark that Lemma \ref{2.62} applied to the special case of
two--sided crossed products (Prop.~\ref{prop2.64}) provides us with
localized left and right coactions of the 
quantum double $\D(\G)$ on the above quantum chain. This generalizes the
$\D(\G)$--cosymmetry discovered by [NSz]. 

Also the construction of the periodic chain by closing a finite open
chain $\A_{0,2i}$ may again be described by a diagonal crossed
product
  $\A_{0,2i}\, \relilr \dG$, 
where $\la$ and $\rho$ are nontrivial only on $\A_0 \>cros \A_1 \<cros
\A_2$ and $\A_{2i-2} \>cros \A_{2i-1}\<cros \A_{2i}$, respectively,
where they are defined as in Proposition \ref{prop2.65}. \\

Let us conclude this section by indicating, how the equivalence of the
 Hopf spin chains of [NSz] and the lattice current algebras of [AFFS] as
shown in [N1] generalizes to the quasi--Hopf setting. 
Following [N1] we define the generating current operators by 
\begin{equation}
  \label{2.628}
  \JJ_{2i+1} : =  (\id\tp A_{2i})(R^{op}) \, \LL_{2i+1},
\end{equation}
where $R\in \G\tp\G$ is supposed to be quasi--triangular.
This yields the following commutation relations
\begin{align*}
  [\eG\tp A_{2i}(a)]\, \JJ_{2i-1} &=
  \JJ_{2i-1} \, [a_{(1)}\tp A_{2i}(a_{(2)})], \quad \forall a \in \G
  \\
  [a_{(1)}\tp A_{2i}(a_{(2)})]\, \JJ_{2i+1} & = \JJ_{2i+1} \, [\eG \tp
  A_{2i}(a)] \\
  \JJ[13]_{2i+1}\, \JJ[23]_{2i+1} &= \hat{R}_{2i}\, \phi_{2i}\, (\cop\tp
  \id)(\JJ)\, \phi^{-1}_{2i+2} \\
   \JJ[13]_{2i-1} \, \hat{R}_{2i} \JJ[23]_{2i+1} & = \JJ[23]_{2i+1}\,
  \JJ[13]_{2i-1} 
\end{align*}
where $\hat{R}_{2i} : = (\id\tp\id\tp A_{2i})(\phi^{213}\, R^{12}\, \phi^{-1})$
and $\phi_{2i} : = (\id\tp\id\tp A_{2i})(\phi)$. These relations
generalize the defining relations of lattice current algebras as given
in [AFFS] to quasi--Hopf algebras. They have also appeared in 
[AGS] as lattice Chern--Simons algebras (restricted to the boundary of
a disk) in the weak quasi--Hopf algebra setting, where the copies
of $\G$ sitting at the sites are interpreted as gauge transformations. \\ 

A more detailed account of these constructions will be given elsewhere [H].

\subsection{Field algebra construction with quasi--Hopf symmetry}
In this section we sketch 
how the field algebra construction of Mack
and Schomerus [MS,S] may be described by two--sided crossed
products.
Here we anticipate the result of Part III, where it is shown how to
generalize our constructions to {\it weak} quasi--Hopf algebras $\G$
(i.e.\ satisfying $\cop(\e) \neq \e\tp\e$).
Choosing $\B \equiv \G$ and $\la_B \equiv \cop$, one gets a
unital algebra $\M_1 :=(\A\tp\G)\reli \dG \cong \A\>cros\dG\<cros \G$
associated with every right $\G$--comodule algebra $(\A,\rho,\phi_\rho)$. 
Here  $\A \subset \M_1$ is to be interpreted as
the algebra of observables, the universal $\la\rho$--intertwiner $\GG$
is a ``master'' 
field operator, 
 and $\G \subset \A\>cros\dG\<cros \G$ represents
the global ``quantum symmetry''. The fields are said to {\it transform
covariantly}, which means that 
\begin{equation}
  \label{2.620}
  [\eG\tp a]\, \GG = \GG \, \la(a) \equiv \GG\, \cop(a),\quad a \in \G,
\end{equation}
whereas the observables $A\in \A$ are {\it $\G$--invariant}
\begin{equation*}
  a \, A = A\, a, \quad \forall a \in \G, A \in \A
\end{equation*}
The linear subspace $\F \equiv \A \>cros \dG : = \A\>cros \dG \<cros \eG$
is called the {\it field algebra}. Note that in the
quasi--coassociative setting $\F$ is {\it not} a subalgebra of
$\A\>cros \dG \<cros \B$. But similarly as in [MS] one may define a new
non--associative 
``product''  $\times$ on $\F$ by setting $A\times A' = A A'$ for $A,A'
\in \A$ and 
\begin{equation}
  \label{2.621}
  \GG[13]\times \GG[23] : = (\cop\tp\id)(\GG)
\end{equation}
This product is quasi--coassociative in the sense that
\begin{equation*}
  [\phi\tp\e] \, \big(\GG[14]\times \GG[24]\big) \times \GG[34] = 
  \GG[14] \times \big(\GG[24] \times \GG[34]\big) \, [\phi \tp \e]
\end{equation*}
and it satisfies
\begin{equation*}
  [\eG\tp\eG\tp a]\, (\GG[13]\times \GG[23]) = (\GG[13]\times
  \GG[23])\, (\cop\tp\id)(\cop(a))
\end{equation*}
which is the reason why it is called {\it covariant product} in
[MS]. Moreover, if $\G$ is quasitriangular,  the field operators
satisfy the braiding relations 
\begin{equation*}
  \GG[13]\times \GG[23] = R^{12} \, (\GG[23]\times \GG[13])\, (R^{-1})^{12}
\end{equation*}
The difference of this approach with the setting of [MS,S] lies in
the fact that here the field operators appear in the form of
irreducible matrix--multiplets 
\begin{equation*}
  F^{ij}_I : =(\pi^{ij}_I \tp \id)(\GG), \quad \pi_I \in \text{Irrep}\, \G
\end{equation*}
Correspondingly, the superselection sectors of the observable algebra
$\A$ are given by the {\it amplimorphisms} 
\begin{equation*}
  \rho^{ij}_I(A) \equiv (\pi^{ij}\tp \id)(\rho(A)) = F^{ik}_I \, A\,
  (F^{kj}_I)^* 
\end{equation*}
which is equivalent to the defining relation \eqref{2.613}. In this
sense the above construction fits to the formulation of DHR--sector
theory as proposed for lattice theories in [SzV,NSz]. In the
terminology of [NSz] the
$\G$--coactions $\rho: \A \to \A\tp\G$ would then be a ``universal
cosymmetry''.

We will show elsewhere [HN2], that such a $\rho$ may indeed be
constructed also for continuum theories (more precisely for
``rational'' theories, where the number of sectors is finite). This
will be done by providing suitable multiplets $W^i_I \in \A$
satisfying
\begin{equation*}
  \sum_i W_I^i \,W_I^{i*} = \eA, \quad \forall I
\end{equation*}
and putting
\begin{equation*}
  \rho^{ij}_I (A) : = W^{i*}_I\,\sigma_I (A)\,W^j_I
\end{equation*}
where $\sigma_I$ is a DHR--endomorphism representing the irreducible
sector $I$. The field operators of [MS,S] are then recovered as 
\begin{equation*}
  \Psi^i_I = W_I^j\,F^{ji}_I
\end{equation*}
In general the $W_I^i$'s will not generate a Cuntz algebra (i.e.\ they will
not be orthogonal isometries)  and therefore the amplimorphisms
$\rho_I$ will be non--unital, implying $\G$ to be a {\it weak}
quasi--Hopf algebra. The reassociator $\phi_\rho \in \A\tp\G\tp\G$
will be obtained by suitably ``blowing up'' the quantum field
theoretic ``$6j$--symbols'' with the help of the $W^i_I$'s. In this
way we will finally be able to show [HN2], that the inclusion
$\A\subset \M_1 \equiv \A\>cros \dG\<cros \G$ provides a finite index
and depth--2 inclusion of von--Neumann factors satisfying $\A' \cap
\M_1 = \G$.


%
%
\part{Weak quasi--Hopf algebras}
In this last part we will generalize the definitions and constructions
of Part II to {\it weak quasi--Hopf algebras} as introduced in
[MS]. This contains the physical relevant examples such as truncated
quantum groups at roots of unity.  
As will be shown, it is nearly straightforward to extend all results
obtained so far to the case of weak quasi--Hopf algebras. The new
feature of {\it weak} quasi--Hopf algebras is to allow the coproduct to
be non unital, i.e. $\cop(\eG)\neq \eG\tp\eG$. This results in a
{\it truncation} of tensor products of
representations, i.e. the representation $(\pi_V\tp
\pi_W)\circ \cop $ operates only on the subspace $V\bo W : = (\pi_V\tp
\pi_W)(\cop(\eG))(V\tp W)$. Also the invertibility requirement on certain
universal 
elements - such as the reassociator or the $R$--matrix - 
is weakened by only postulating the existence of so--called {\it
  quasi--inverses}. Correspondingly coactions and two--sided coactions
of a weak quasi-Hopf 
algebra are no longer supposed to
be unital and the associated reassociators are only required to possess
quasi--inverses. The diagonal crossed products $\M_1 \equiv\dG \reli \M$ may
then again be defined by the same relations as before, with the
additional requirement that the universal $\la\rho$--intertwiner $\GG
\in \G\tp \M_1$ has to satisfy $\GG \, \la(\eM) \equiv
\rho^{op}(\eM)\, \GG = \GG$. 
This will also imply that now as a linear space $\M_1$ is only
isomorphic to a certain subspace $\dG \relilr \M \subset \dG \tp \M$ (or
$\M\relilr \dG \subset \M\tp \dG$). More specifically Theorem II
now reads
\vspace{0.25cm}
  \begin{sloppypar} \noindent{\bf Theorem III}
{\it 
Let $\G$ be a finite dimensional weak quasi--Hopf algebra and let
$(\la,\Fi{\la},\rho,\Fi{\rho},\Fi{\la\rho})$ be a quasi--commuting pair
of (left and right) $\G$--coactions on an associative algebra $\M$.
\begin{enumerate}
\item Then their exists a unital associative algebra extension $\M_1
  \supset \M$ together with a linear map $\Gamma : \dG\pfeil \M_1$
  satisfying the following universal property: \\
  $\M_1$ is algebraically generated by $\M$ and $\Gamma(\dG)$ and for
  any algebra map $\gamma :\M\pfeil\A$ into some target algebra $\A$
  the relation 
  \begin{equation*}
   \gam_T (\Gamma(\vi)) = (\vi\tp\id)(\TT)
  \end{equation*}
provides a one--to--one correspondence between algebra maps
$\gam_T: \M_1 \pfeil \A$ extending $\gam$ and elements $\TT \in
\G\tp\A$ satisfying $(\ep\tp\id)(\TT) = \gam(\eM)$ and 
\begin{align*}
\TT \, \la_\A (m) &= \rho^{op}_\A (m)\,\TT, \quad \forall m \in \M \\
 \TT \, \la_\A (\eM)&\equiv\rho_\A^{op}(\eM)\, \TT = \TT \\
 (\Fi[312]{\rho})_\A \, \TT[13] \, (\Fi[132]{\la\rho})^{-1}_\A \, \TT[23] \,
 (\Fi{\la})_\A &= (\cop \tp\idA)(\TT),
\end{align*}
where $\la_\A (m) := (\id\tp\gam)(\la(m))$, $(\Fi{\la})_\A : =
(\idG\tp\idG\tp\gam)(\Fi{\la})$, etc., and where $\phi^{-1}_{\la\rho}$
is the quasi--inverse of $\phi_{\la\rho}$. 
\item If $\M\subset \tilde{M}_1$ and $\tilde{\Gamma}:
  \dG\pfeil\tilde{\M}_1$ satisfy 
  the same universality property as in part 1.), then there exists a
  unique algebra isomorphism $f: \M_1\pfeil \tilde{\M}_1$ restricting
  to the identity on $\M$, such that $\tilde{\Gamma} = f\circ \Gamma$ 
\item
There exist elements $p_\la \in \G\tp\M$ and $q_\rho \in \M\tp \G$
such that the linear maps 
  \begin{align*}
   \mu_L : \dG\tp\M \ni (\vi\tp m) &\mapsto (\id \tp
   \vi_{(1)})(q_\rho)   
   \Gamma (\vi_{(2)}) m \, \in \M_1\\
   \mu_R : \M\tp\dG \ni (m\tp \vi)&\mapsto m \Gamma (\vi_{(1)})
   (\vi_{(2)}\tp \id)(p_\la) \, \in \M_1
  \end{align*}
are surjective.
\item
Let $P_L : \dG \tp \M \to \dG \tp \M$ and $P_R: \M\tp\dG \to \M\tp\dG$
be the linear maps given by 
\begin{align*}
  P_L(\vi\tp m) \equiv \vi\reli m : = &
    \vi_{(2)} \tp
    (\vi_{(3)}\tp\id\tp\hat{S}^{-1}(\vi_{(1)}))\big(\delta_l(\eM)\big) 
      m \\
   P_R( m \tp \vi) \equiv m\reli \vi : =&
   m ( \hat{S}^{-1}(\vi_{(3)})\tp\idM\tp
    \vi_{(1)})\big(\delta_r(\eM)\big) \tp \vi_{(2)} 
\end{align*}
where $\delta_l = (\la\tp\id)\circ \rho$ and $\delta_r =
(\rho\tp\id)\circ \la$. 
Then $P_L$ and $P_R$ are projections and $\Ker \mu_{L/R} = \Ker P_{L/R}$.
\end{enumerate}
  }
\end{sloppypar}\vspace{0.25cm}
Part 3.~and 4.~of Theorem III imply that we may put $\dG \relilr 
M : = P_L (\dG \tp \M)$ and $\M\relilr\dG : = P_R (\M\tp \dG)$ to
conclude that analogously as in Eqs.~\eqref{2.05} and \eqref{2.06} the
restrictions 
\begin{align*}
  \mu_L : \dG\relilr \M &\pfeil \M_1 \\
   \mu_R : \M\relilr \dG & \pfeil \M_1
\end{align*}
are linear isomorphisms inducing a concrete realization of the
abstract algebra $\M_1$ on the subspaces $\dG\relilr\M\subset
\dG\tp\M$ and $\M\relilr\dG \subset \M\tp\dG$, respectively. As
before, we call these concrete realizations the {\it left and right
  diagonal crossed products}, respectively, associated with the
quasi--commuting pair of $\G$--coactions
$(\la,\rho,\phi_\la,\phi_\rho,\phi_{\la\rho})$ on $\M$\footnote{
Actually we have again four versions of diagonal crossed products as
in Part II, given by \\
$\dG \relidl \M \equiv \dG \relilr \M$,
$\dG\relidr \M$, $\M\relidl \dG$ and $\M\relidr \dG \equiv
\M\relilr\dG$.}. 
To actually prove Theorem~III we follow the same strategy as in
Part~II., i.e. we first construct these diagonal crossed products explicitly
and then show that they solve the universal properties defining
$\M_1$.


\section{Weak quasi--bialgebras}
We start with a little digression on the notion of {\it
  quasi--inverses}. Let $\A$ be an associative algebra and let $x,p,q
\in \A$ satisfy
\begin{align}
  \label{12.1}
   & px = x =xq \\
  \label{12.2}
   & p^2 = p, \quad\quad q^2 = q
\end{align}
Then we say that $y \in \A$ is a {\it quasi--inverse} of $x$ with
respect to $(p,q)$, if 
\begin{equation}
  \label{12.3}
  yx =q,  \quad\quad  xy =p, \quad\quad   yxy = y 
\end{equation}
Clearly, given $(p,q)$, a quasi--inverse of $x$ is uniquely determined,
provided it exists. This is why we also write $y =x^{-1}$, if the
idempotents $(p,q)$ are understood. Also note that we have $qy=y=yp$
and $xyx =x$ and therefore $x$ is the quasi--inverse of $y$ with
respect to $(q,p)$. All this generalizes in the obvious way to
$\A$--module morphisms $x\in \Hom_\A(V,W)$, $p \in \End_\A (W)$ and $q
\in \End_\A 
(V)$, in which case the quasi--inverse would be an element  $x^{-1}
\in \Hom_\A(W,V)$.

Note that in place of \eqref{12.1} we could equivalently add to
\eqref{12.3} the requirement 
\begin{equation}
  \label{12.5}
  xyx =x
\end{equation}

In our setting of weak quasi--Hopf algebras the idempotents $p,q$
always appear as images of $\eG$ of non--unital algebra maps defined
on $\G$, like $\cop(\e)$, $\cop^{op}(\e)$,
$(\cop\tp\id)(\cop(\e)),\dots$ etc., whereas the element $x$ will be
an intertwiner between two such maps, like a reassociator $\phi$
etc. Hence, throughout we will adopt the convention that if 
$\alpha : \G \to \A$ and $\beta : \G\to \A$ are two algebra maps and
$x \in \A$ satisfies 
\begin{equation*}
  x \, \alpha(g) = \beta(g)\, x ,\quad \forall g \in \G
\end{equation*}
then the quasi--inverse
$y=x^{-1}\in \A$ is defined to be the unique (if existing) element
satisfying 
\begin{align*}
  yx = \alpha(\e),&  \quad xy = \beta(\e)\\
 xyx=x,& \quad yxy =y  
\end{align*}
Clearly this implies conversely 
\begin{equation*}
  \alpha(g) \, y = y \, \beta(g), \quad \forall g \in \G
\end{equation*}
and therefore $x = y^{-1}$. We also note the obvious identities
\begin{align*}
  \beta(g) = x \, \alpha(g) \, x^{-1},& \quad \alpha(g) = x^{-1}\,
  \beta(g) \, x \\
  x \, \alpha(\e) = \beta(\e)\, x = x ,& \quad \alpha(\e) \, x^{-1} =
  x^{-1} \, \beta(\e) = x^{-1}
\end{align*}

After this digression we now define, following [MS]
  a {\it weak quasi--bialgebra} $(\G,\e,\cop,\ep,\phi)$ to be an
  associative algebra $\G$ with unit $\e$, a non--unital algebra
  map $\cop : \G \to \G\tp\G$, an algebra map $\ep: \G \to \Co$ and an
  element
  $\phi \in \G\tp\G\tp\G$ satisfying \eqref{eq11}-\eqref{eq13},
  whereas \eqref{eq14} is replaced by 
  \begin{equation}
    \label{3.1}
    (\id\tp\ep\tp\id)(\phi) = \cop(\e) 
  \end{equation}
and where in place of invertibility $\phi$ is supposed to have a
quasi--inverse $\Fib{} \equiv \phi^{-1}$ with respect to the
intertwining property \eqref{eq11}. Hence we have $\phi\Fib{}\phi =
\phi$, $\Fib{}\phi\Fib{} = \Fib{}$ as well as 
\begin{equation}
   \label{3.2}
    \phi \, \Fib{} = (\id\tp\cop)(\cop(\e)),\quad 
   \Fib{}\,\phi = (\cop\tp\id)(\cop(\e)) 
\end{equation}
implying the further identities
\begin{align}
  \label{3.4}
  (\id\tp\cop)(\cop(a)) &= \phi \, (\cop\tp\id)(\cop(a)) \, \Fib{},
  \quad \forall a \in \G \\
   \label{3.4'}
  \phi = \phi \, (\cop\tp\id)(\cop(\e)), &\quad \Fib{} = \Fib{}\,
  (\id\tp\cop)(\cop(\e))\\
  \label{3.3}
  (\id\tp\ep\tp\id)(\Fib{}) & = \cop(\e)
\end{align}
A weak quasi--bialgebra is called {\it weak quasi--Hopf algebra}, if there
exists  a unital algebra antimorphism $S$ and elements $\alpha,\beta\in
\G$ satisfying \eqref{eq16} and \eqref{eq17}. We will also always
suppose that $S$ is invertible. The remarks about the quasi--Hopf
algebras $\G_{op}$, $\G^{cop}$ and $\G_{op}^{cop}$ remain valid as in
Section 6.

A quasi--invertible element $F \in \G\tp\G$ satisfying 
$(\ep\tp\id)(F) = (\id\tp\ep)(F)  = \e$ induces a {\it twist transformation}
 from $(\G,\cop,\ep,\phi)$ to $(\G,\cop_F,\ep,\phi_F)$ as in \eqref{eq113}
 and \eqref{eq114}, where 
   $\G_F : = (\G,\e,\cop_F,\ep,\phi_F)$ is again a weak
   quasi--bialgebra. The two bialgebras $\G_F$ and $\G$ are called
   {\it twist--equivalent}.
 
Finally the properties of the twists $f,h$ defined as in \eqref{eq118} and
\eqref{defh} are still valid. In particular \eqref{eq119} defines the
quasi--inverse of $f$ with respect to the intertwining property
\eqref{eq120}.


\section{Weak coactions}
The notion of coactions may easily be generalized as well.  By a left
$\G$--coaction $(\la,\phi_\la)$ of a weak quasi--bialgebra $\G$ on
a unital algebra $\M$  we mean a (not necessarily unital) algebra map
$\la: \M \to \G\tp\M$ and a quasi--invertible element
$\phi_\la  \in \G\tp\G\tp\M$ satisfying
\eqref{eq21}-\eqref{eq22b} as in Definition \ref{def2.3} and 
\begin{equation}
 \label{3.6}
  (\id\tp\ep\tp\id)(\phi_{\la}) = (\ep\tp\id\tp\id)(\phi_\la) = \la(\eM) 
\end{equation}
The definition of right coactions is generalized analogously. 
Lemma \ref{lemtwist} about twist equivalences of coactions stays
valid, where one has to make the adjustments that a twist $U \in
\M\tp\G$ only is  supposed to be quasi--invertible.

By now it should become clear how one has to proceed: Definition
\ref{def2.4} of two--sided coactions is generalized by allowing
$\delta$ to be non--unital and $\Psi$ to be non--invertible but with
quasi--inverse $\bar{\Psi}\equiv \Psi^{-1}$ 
and by replacing \eqref{eq24f} by 
\begin{equation}
  \label{3.9}
(\idG\tp\ep\tp\idM\tp\ep\tp\idG)(\Psi) =
(\ep\tp\idG\tp\idM\tp\idG\tp\ep)(\Psi) = \delta(\eM)  
\end{equation}
The definitions of quasi--commuting pairs of coactions, twist
equivalences of two--sided coactions etc. are generalized
similarly. With these adjustments all results of Section 8 stay
valid for weak quasi--Hopf algebras and are proven analogously. \\

Let us now shortly review the representation theoretical
interpretation as given in Section 9 in the light of the present
setting. Due to the coproduct being non--unital the definition of the
tensor product functor in $\Rep \G$ has to be slightly modified. First note
that the element $\cop(\e)$ (as well as higher coproducts of $\e$)
is idempotent and commutes with all elements in $\cop(\G)$. Thus,
given two representations $(V,\pi_V), (W,\pi_W)$, the operator
$(\pi_V\tp\pi_W)(\cop(\e))$ is a  projector,
whose image is precisely the $\G$--invariant subspace of $V\tp W$ on
which the tensor product representation operates non trivial. Thus one
is led to define the tensor product $\bo$ of two representations of $\G$
by setting
\begin{equation*}
  V\bo W : = (\pi_V \tp \pi_W)(\cop(\e))\, V\tp W, \quad
  \pi_V\bo \pi_W : = (\pi_V \tp \pi_W)\circ \cop |_{ V\bo W}
\end{equation*}
One readily verifies that with these definition $\phi_{UVW}$ -
restricted to the subspace $(U\bo V)\bo W$ - furnish a natural family of
isomorphisms defining an associativity constraint for the tensor
product functor $\bo$. Moreover, $\Rep \G$ becomes a rigid
monoidal category with rigidity structure defined as before by
\eqref{2.5.3}-\eqref{2.5.8.}. 

The {\it left action} of $\Rep \G$ on $\Rep \M$ induced by a left
$\G$--coaction $(\la,\phi_\la)$ on $\M$ has to be modified analogously
by defining
\begin{equation*}
  V\so \gotH : = (\pi\tp \gamma)(\la(\eM))\, V \tp \gotH, \quad 
  \pi \so \gamma : = (\pi\tp\gamma )\circ \cop |_{V \so \gotH}
\end{equation*}
The modifications of right actions and of two--sided actions of $\Rep
\G$ on $\Rep \M$ (induced 
by $(\rho,\phi_\rho)$ and by $(\delta,\Psi)$, respectively) should by now
be obvious and are left to the reader. 

With these adjustments all categorical identities such as the
definition of natural families and commuting diagrams given in
Section\ 9
 stay valid. Translating these into algebraic
identities one has to take some care with identities in higher tensor
products of $\G$. The only equations which have to be modified are
\eqref{2.5.20},\eqref{2.5.21}, where the r.h.s. becomes $\la(\eM)$
instead of $\eG\tp\eM$ and similarly (\ref{2.5.27}/\ref{2.5.28}) and 
(\ref{2.5.47}/\ref{2.5.48}) where the r.h.s. has to be replaced by
$\rho(\eM)$ and by $\delta({\eM})$, respectively. This is rather obvious
from the categorical point of view, since for example \eqref{2.5.20}
is directly connected with \eqref{2.5.33} where the r.h.s. is given by
$\id_{V\so \gotH} \equiv (\pi_V\tp \pi_\gotH)(\la(\eM))$.


\section{Diagonal crossed products}
The definition of diagonal crossed products $\dG \relid \M$ and
$\M\relid \dG$ as equivalent algebra extensions of $\M$, given in
Definition \ref{def2.9} and Theorem \ref{thm2.5} need some more care
in the present context. We will proceed in two steps. First we
define an associative algebra structure on $\dG\tp\M$ (or $\M \tp
\dG$) exactly as in Definition \ref{def2.9}. Unfortunately in general
this algebra 
is not unital unless the two--sided coaction $\delta$ is
unital. But the element $\edG\tp\eM$ is still a right unit ($\eM\tp
\edG$ is still a left unit)
and in particular idempotent. The second step then consists in
defining the subalgebra 
$\dG\relid \M \subset \dG \tp \M$ as the right ideal generated by
$\edG\tp\eM$, i.e. $\dG \relid \M : =
(\edG\tp\eM) \cdot (\dG \tp \M)$ (the left ideal generated by $\eM\tp\edG$,
i.e. $\M\relid \dG := (\M\tp \dG)\cdot(\eM\tp\edG)$). These  algebras are
then unital algebra extensions of $\M \equiv \edG \reli \M$ and of $\M
\equiv \M\reli\edG$, respectively. As in Section 10 one may proceed
to a description by left and right $\delta$--implementers and
equivalently by $\la\rho$--intertwiners, thus getting a proof of
Theorem~III.
\begin{definition}{\rm
  \label{def3.1}
 Let $(\delta,\Psi)$ be a two--sided coaction
 of a weak quasi--Hopf algebra $\G$ on an algebra $\M$. We define
 $\dG\tp_\delta \M$ to be the vector space $\dG \tp \M$ with
 multiplication structure given as in \eqref{def261} and the
 left {\it diagonal crossed product} $\dG \relid \M$ to be the subspace 
 \begin{equation}
   \label{3.11}
   \dG\relid \M : = (\edG\tp \eM) \cdot (\dG\tp_\delta \M)
 \end{equation}
Analogously $\M\tp_\delta \dG$ is defined to be the vector space
$\M\tp \dG$ with multiplication structure \eqref{def262} and the
right {\it diagonal crossed product} $\M\relid \dG$ to be the subspace
\begin{equation}
  \label{3.11'}
   \M\relid\dG : = (\M\tp_\delta \dG) \cdot (\eM\tp\edG)
\end{equation}
The elements spanning $\dG\relid \M$ and $\M\relid \dG$ will be denoted by,
respectively 
\begin{align}
  \label{3.10a}
  \vi\reli m &: = (\edG \tp \eM)(\vi\tp m) \equiv (\vi\tp\eM)(\eM\tp
  m) \equiv \vi_{(2)} \tp \big(\hat{S}^{-1}(\vi_{(1)})\re \eM \li
  \vi_{(3)}\big) m \\
   \label{3.10b}
  m\reli \vi & : = (m\tp\vi)(\eM\tp\edG) \equiv (m\tp \edG)(\eM\tp
  \vi) \equiv  m \big( \vi_{(1)} \re \eM \li
  \hat{S}^{-1}(\vi_{(3)})\big) \tp \vi_{(2)}
\end{align}}
\end{definition}
Note that $\dG\tp_\delta \M = \dG\relid\M$, if $\delta(\eM) =
\eG\tp\eM\tp\eG$, which means that the above definition generalizes
Definition \ref{def2.9}. We now state the analogue of
Theorem \ref{thm2.5}.
\begin{theorem}$\,$
  \label{thm3.2}
  \begin{itemize}
  \item[(i)] 
  $\dG\tp_\delta \M$ and $\M\tp_\delta\dG$ are associative
    algebras with left unit $\edG\tp \eM$ and right unit $\eM\tp\edG$,
  respectively. Consequently, 
   the diagonal crossed products $\dG\relid \M$ and
    $\M\relid\dG$ are subalgebras of $\dG\tp_\delta\M$ and
    $\M\tp_\delta\dG$, respectively, with unit given by 
    $\edG\reli \eM \equiv \edG \tp \eM$ and $\eM\reli \edG \equiv
    \eM\tp \edG$, respectively.
  \item[(ii)] $\M \equiv  \edG \tp \M = \edG \reli \M \subset
    \dG\relid \M$ and
    $\M  \equiv \M\tp \edG = \M\reli \edG \subset \M\relid \dG$  are
    unital algebra inclusions.
  \item[(iii)] The algebras $\dG\reli \M$ and $\M\reli\dG$ provide
    equivalent extensions of $\M$ with isomorphism given as in
    \eqref{2.5.1.1},\eqref{2.5.1.2}. 
  \end{itemize}
\end{theorem}
\begin{proof}
We will sketch the proof of part (i) for $\dG\tp\M$, the case $\M\tp\dG$
  being analogous. From \eqref{def261} one computes that
  \begin{align}
    \notag
    (\vi\tp m)(\edG\tp\eM)& = (\vi \tp m)\\
   \notag
   (\edG\tp\eM)(\vi\tp m) &= \vi_{(2)} \tp \big(
    \hat{S}^{-1}(\vi_{(1)}) \re \eM\li \vi_{(3)}\big) \, m \\
     \notag
      & = \big( \e_{(-1)} \arr \vi \arl \Si(\e_{(1)})\big) \tp \e_{(0)} \, m
  \end{align}
 where $\delta(\eM) = \e_{(-1)} \tp \e_{(0)} \tp \e_{(1)}$. This shows
 that $\edG \tp\eM$ is a right unit in $\dG \tp_\delta \M$ but in
 general not a left unit.

To proof the associativity of the product one proceeds as in the proof
of Theorem \ref{thm2.5}. Here one has to take some notational care
when translating 
(\ref{24.21}/\ref{24.22}) into relations of a generating matrix
$\LL$. First, it is necessary to distinguish $\LLt \equiv \Lt^1 \tp
\Lt^2 := e_\mu \tp (e^\mu \tp \eM) \in \G\tp (\dG\tp \M)$ and 
$\LL : = e_\mu \tp (e^\mu \reli \eM) \in \G\tp (\dG\reli \M)$.
Eq.~(\ref{216a}) must then be replaced by 
\begin{equation}
  \label{Lt}
  [\eG \tp \eM]\, \LLt = \LL
\end{equation}
and Eqs.~\eqref{eq216b}, \eqref{eq216c} are at first only valid for
$\LL$ but not for $\LLt$, since $\eM \equiv \edG \tp \eM$ is not a
left unit in $\dG \tp_\delta \M$. This can be cured by rewriting for
example \eqref{eq216b} for $\LLt$ more carefully as 
\begin{equation*}
  (\Lt^1 \tp m \Lt^2) = \Si(\me) \, \Lt^1\, \mme \tp \Lt^2\, \mo,
  \quad \forall m \in \M
\end{equation*}
in which form it would still be valid.
Taking this into account and
using that $\bar{\Psi} \, (\eG\tp \delta(\eM) \tp \eG) = \bar{\Psi}$,
the proof proceeds as the one of Theorem \ref{thm2.5} (i). 
The proof of the remaining parts of Theorem \ref{thm2.5} is
straightforwardly  adjusted in the same spirit.
\end{proof}
 From now on we will disregard the ``unphysical'' non--unital algebras
 $\M\tp_\delta\dG$ and $\dG\tp_\delta\M$, and stay with $\dG \reli \M$
 and $\M\reli\dG$ as our objects of interest. With the appropriate
 notations \eqref{3.10a},\eqref{3.10b} all relations of Section 10.1
 remain valid for these algebras. We also remark that the left
 multiplication by $\edG\tp \eM$ (right multiplication by
 $\eM\tp\edG$) precisely gives the projections $P_{L/R}$ mentioned in
 part 4. of Theorem III.

Defining left and right $\delta$--implementers as in Definition~\ref{def2.7}
 all results in Section 10.2 also carry
over to the present setting. 
A $\la\rho$--intertwiner is then supposed to have the additional property
that 
\begin{equation}
  \label{3.15}
  \TT \, \la_\A(\eM) \equiv \rho_\A^{op}(\eM) \, \TT = \TT
\end{equation}
With this the one--to--one correspondence $\TT \leftrightarrow \RR$  and $\TT
\leftrightarrow \LL$  of Proposition \ref{prop2.15} is still valid,
where \eqref{3.15} becomes equivalent to
$\LL\prec \delta(\eM) = \LL$ and $\delta(\eM) \succ \RR =\RR$,
respectively, which follow from \eqref{252.3} and \eqref{252.4}.
One may now prove Theorem III analogously as Theorem II
where the modifications in part 3.\ and 4.\ have their origin in
\eqref{3.15}.


\section{Examples and applications}
Finally, all examples given in Section 11 generalize to weak
quasi--Hopf algebras.
The definition of the quantum double $\D(\G)$ for weak quasi--Hopf
algebras as given in Theorem \ref{2.63} yields a weak
quasi--bialgebra. 

The definition of {\it two--sided crossed products} as in
 Proposition \ref{prop2.64} (i.e.  an alternative
description of the algebra $\M\relilr\dG$ for the case that $\M = \A\tp
\B$), has to be slightly modified. The
unital algebra $\A\crosr\dG \crosl \B$ is now defined on the subspace
 of $\A\tp\dG\tp\B$ given by 
\begin{equation}
    \label{15.1}
  \A\>cros\dG\<cros\B : = \text{span}\big\{
   A\>cros\vi\<cros B \equiv A\eA^{(0)} \tp
   \eA^{(1)}\re\vi\li\eB^{(-1)}\tp \eB^{(0)}B \mid
   A\in\A,B\in\B,\vi\in\dG \big\},  
\end{equation}
where $\e^{(0)}_\A \tp \e^{(1)}_\A \equiv \rho(\e_\A)$ and 
$\e_\B^{(-1)} \tp \e^{(0)}_\B \equiv \la(\e_\B)$
Again we have a linear bijection \\
$\mu: \A\>cros\dG\<cros\B \to \M_1 = (\A\tp\B) \relilr \dG$
\begin{equation*}
  \mu (A\>cros\dG\<cros B) = A\,\Gamma(\vi)\, B
\end{equation*}
inducing the multiplication rule described by Eq.~\eqref{2.615}.
With these notations all results of Section~11.3 are still valid for
weak quasi--Hopf algebras. As for part (iii) of Proposition \ref{prop2.65}
we note that \eqref{15.1} still allows the identification
\begin{equation*}
  (A\>cros \vi \<cros B)\>cros\psi\<cros C = A\>cros\vi \<cros (B
  \>cros\psi \<cros C)
\end{equation*}
Moreover, putting $(\B,\la,\phi_\la) = (\G,\cop,\phi)$ this is also
the setting underlying the field algebra constructions proposed in
Section 11.4.


\begin{appendix}
\part*{Appendix}
\addcontentsline{toc}{part}{\protect\numberline{}{Appendix}}
\setcounter{section}{1}
\setcounter{equation}{0}
In this appendix we will give the proofs which have been omitted in
the previous sections. Let us first
collect additional identities of the reassociator $\Psi \in
\G\tp\G\tp\M\tp\G\tp\G$ given in Definition \ref{def2.4}, which follow
by applying $\ep_i
\tp\idM \tp \ep_j, \, 1\leq i,j \leq 3,$ to both sides of
\eqref{eq24e}, where $\ep_i : \G^{\tp^3} \rightarrow \G^{\tp^2}$ is
given by acting with $\ep$ on the $i$--th tensor factor. 
\begin{align}
 & [(\idG^{\tp^2} \tp\idM\tp \ep \tp\idG)(\Psi)\tp\eG]\,
  [(\ep\tp\cop\tp\idM\tp\idG^{\tp^2})(\Psi)]  \notag\\
  & \quad\quad \quad\quad    \label{ep1} =
  [(\ep\tp\idG\tp \la\tp\idG^{\tp^2})(\Psi)]\,
  [(\idG^{\tp^2} \tp\idM\tp\ep\tp\cop)(\Psi)] \\
 & [(\idG^{\tp^2} \tp\idM\tp\idG\tp\ep)(\Psi)\tp\eG]\,
  [(\ep\tp\cop\tp\idM\tp\idG^{\tp^2})(\Psi)]  \notag\\
  & \quad\quad \quad\quad   \label{ep2} =
  [(\ep\tp\idG\tp\delta\tp\ep\tp\idG)(\Psi)]\, \Psi \\
&  [\eG\tp(\ep\tp\idG\tp\idM\tp\ep\tp\idG)(\Psi)\tp\eG]\, \Psi \notag \\
  & \quad\quad \quad\quad    \label{ep3} =
   [(\idG\tp\ep\tp \la\tp\idG^{\tp^2})(\Psi)]\,
  [(\idG^{\tp^2} \tp\idM\tp\ep\tp\cop)(\Psi)] \\
&  [\eG\tp (\idG\tp\ep\tp\idM\tp\idG^{\tp^2})(\Psi)]\,
  [(\idG^{\tp^2} \tp\idM\tp\cop\tp\ep)(\Psi)]  \notag\\
  & \quad\quad\quad\quad     \label{ep4} =
   [(\idG^{\tp^2} \tp\rho\tp\idG\tp\ep)(\Psi)] \,
   [(\cop\tp\ep\tp\idM\tp\idG^{\tp^2})(\Psi)] \\
&  [\eG\tp (\ep\tp\idG\tp\idM\tp\idG^{\tp^2})(\Psi)]\,
  [(\idG^{\tp^2} \tp\idM\tp\cop\tp\ep)(\Psi)]  \notag\\
  & \quad\quad\quad\quad     \label{ep5} =
   [(\idG \tp\ep\tp\delta\tp\idG\tp\ep)(\Psi)] \, \Psi \\
&  [\eG\tp(\idG\tp\ep\tp\idM\tp\idG\tp\ep)(\Psi)\tp\eG]\, \Psi \notag \\
  & \quad\quad\quad\quad     \label{ep6} =
   [(\idG^{\tp^2}\tp\rho\tp \ep\tp\idG)(\Psi)]\,
  [(\cop\tp\ep\tp\idM\tp\idG^{\tp^2})(\Psi)].
\end{align}
Here we have used $\ep_i (\phi) = \eG\tp\eG$ and we have introduced
the notation 
  $\la : = (\idG\tp\idM\tp\ep)\circ \delta$, 
  $\rho : = (\ep\tp\idM\tp\idG)\circ \delta$.
Note that Eqs. \eqref{ep4} - \eqref{ep6} are the left--right mirror
images of \eqref{ep1} - \eqref{ep3}, and that the axioms \eqref{eq24c}
- \eqref{eq24f} are symmetric under this reflection.

\begin{proof}[{\bf Proof of Proposition \ref{prop20}:}] 
We prove the implication   $(ii)\Rightarrow (i)$:   
So given that $(\delta_l,\Psi_l)$ is a two--sided coaction we use the
counit axioms for $\Fi{\la}$ and $\Fi{\rho}$ to conclude
\begin{align}
  \label{2.36a}
 (\idG^{\tp^2}\tp\idM\tp\idG\tp\ep)(\Psi_l) &= 
     (\idG\tp\la\tp\idG)(\Fi{\la\rho})\,(\Fi{\la}\tp\eG) \\
  \label{2.36b}
 (\idG\tp\ep\tp\idM\tp\idG\tp\ep)(\Psi_l) &= \Fi{\la\rho}
\end{align}
In particular this implies $(\ep\tp\idM\tp\idG)(\Fi{\la\rho}) =
\eM\tp\eG$ and $(\idG\tp\idM\tp\ep)(\Fi{\la\rho}) = \eG\tp\eM$ and
therefore also
\begin{align}
  \label{2.36c}
  (\ep\tp\idG\tp\idM\tp\ep\tp\idG)(\Psi_l) &= \eG\tp\eM\tp\eG \\
    \label{2.36d}
  (\id^{\tp^2}\tp\idM\tp\ep\tp\ep)(\Psi_l) &= \Fi{\la} \\
  \label{2.36f}
  (\ep\tp\ep\tp\idM\tp\id^{\tp^2})(\Psi_l) &= \phi^{-1}_{\rho} 
\end{align}
Now we use $\la = (\idG\tp\idM\tp\ep)\circ \delta_l$ and 
$\rho = (\ep\tp\idM\tp\idG) \circ \delta_l$ implying
\begin{equation*}
  \delta_r = (\idG\tp\ep\tp\idM\tp\ep)\circ \delta_l^{(2)} =
     Ad \,\Fi{\la\rho} \circ \delta_l
\end{equation*}
by \eqref{2.36b} and the assumption \eqref{eq24e} for
$(\delta_l,\Psi_l)$. This proves \eqref{eq25}.
Eq. \eqref{eq26} now follows by acting with $(\idG^{\tp^2}
\tp\idM\tp\idG\tp\ep)$ on \eqref{ep4}. 

To prove Eq. \eqref{eq27} we first use that by \eqref{2.36c},
\eqref{2.36f} and \eqref{ep3} 
\begin{equation*}
  (\la\tp\idG^{\tp^2})(\phi^{-1}_{\rho}) =
  (\ep\tp\idG\tp\idM\tp\idG^{\tp^2})(\Psi_l).
\end{equation*}
Hence Eq. \eqref{eq27} follows from
\begin{align*}
  & (\id\tp\rho\tp\id)(\Fi{\la\rho}) \,(\Fi{\la\rho}\tp\eG)
    (\la\tp\id\tp\id)\Fii{\rho} \\
       &\quad\quad =
[(\idG\tp\ep\tp\rho\tp\idG\tp\ep)(\Psi_l)]
  [(\idG\tp\ep\tp\idM\tp\idG\tp\ep)(\Psi_l)\tp\eG]
[(\ep\tp\idG\tp\idM\tp\idG^{\tp^2})(\Psi_l)] \\
       &\quad\quad =
    [(\idG\tp\ep\tp\rho\tp\idG\tp\ep)(\Psi_l)]
[(\idG\tp\ep\tp\idM\tp\idG^{\tp^2})(\Psi_l)] \\
       &\quad\quad =
 [\eG\tp(\ep\tp\ep\tp\idM\tp\idG^{\tp^2})(\Psi_l)] 
 [(\idG\tp\ep\tp\idM\tp\cop\tp\ep)(\Psi_l)] \\
      &\quad\quad =
 (\eG\tp \phi_\rho^{-1})\, (\idG\tp\idM\tp\cop)(\phi_{\la\rho})
\end{align*}
Here we have acted with $(\idG\tp\ep\tp\idM\tp\idG^{\tp^2})$ on
\eqref{ep2} in the second equation  and on \eqref{ep4} in the third
equation.
Hence we have shown $(ii) \Rightarrow (i)$. The implication
$(iii)\Rightarrow (i) $ is proven analogously. \\

We are left with showing that under the conditions $(i) - (iii)$ $U
\equiv \Fi{\la\rho}$ provides a twist from $(\delta_l,\Psi_l)$
to $(\delta_r,\Psi_r)$. Now the intertwiner property \eqref{eq2.31}
holds by definition and \eqref{eq2.33} follows from \eqref{eq27a}. To
prove \eqref{eq2.32} for $\Psi' = \Psi_r$ and $\Psi = \Psi_l$ we note
that \eqref{eq27} and \eqref{eq215} imply
\begin{equation*}
\Psi_l = (\idG\tp\delta\tp\id)(\phi^{-1}_{\la\rho})
[\eG\tp(\la\tp\idG\tp\idM)(\phi^{-1}_\rho)]\,
(\idG\tp\la\tp\cop)(\Fi{\la\rho})\,[\phi_\la \tp\eG\tp\eG] 
\end{equation*}
and therefore, putting $U=\phi_{\la\rho}$ in \eqref{eq2.32} 
\begin{align*}
  &(\eG\tp U\tp\eG)\,(\idG\tp\delta\tp\idG)(U)\,\Psi_l \\
     &\quad =  
     [\eG\tp(\Fi{\la\rho}\tp\eG)(\la\tp\idG\tp\idG)\Fii{\rho}]
     (\idG^{\tp^2} \tp\idM\tp\cop)\Big(
     (\idG\tp\la\tp\idG)(\Fi{\la\rho})(\Fi{\la}\tp\eG)\Big) \\
  &\quad =  
     [\eG\tp(\idG\tp\rho\tp\idG) \Fii{\la\rho}]
     [\eG\tp\eG\tp\phi^{-1}_\rho]\, \\
   &\quad\quad\quad
     (\idG^{\tp^2} \tp\idM\tp\cop)\Big( (\eG\tp\Fi{\la\rho})
     (\idG\tp\la\tp\idG)(\Fi{\la\rho})(\Fi{\la}\tp\eG)\Big) \\
  &\quad = 
 [\eG\tp(\idG\tp\rho\tp\idG) \Fii{\la\rho}]
     [\eG\tp\eG\tp\phi^{-1}_\rho]\,
     (\idG^{\tp^2} \tp(\idM\tp\cop)\circ\rho)(\Fi{\la}) \,
     (\cop\tp\idM\tp\cop)(\Fi{\la\rho}) \\
  &\quad =  
\Psi_r \, (\cop\tp\idM\tp\cop)(U)
\end{align*}
Here we have used \eqref{eq27} in the second equality, \eqref{eq26} in
the third equality and \eqref{eq23} in the last equality. This proves
Eq. \eqref{eq2.32} for $\Psi = \Psi_l$ and $\Psi'=\Psi_r$. Hence
$\Fi{\la\rho}$ provides a twist from $(\delta_l,\Psi_l)$ to
$(\delta_r,\Psi_r)$. This concludes the proof of Proposition
\ref{prop20}.
\end{proof}
\begin{proof}[{\bf Proof of Proposition \ref{prop21} (ii):}]
  The identities $\delta_{l/r} = Ad U_{l/r} \circ \delta$ follow
  immediately from \eqref{eq24c}, \eqref{la2} and \eqref{rho2}. We are
  left to show that
  \begin{equation}
    \label{2.36.4}
   (\eG\tp U_l\tp\eG)\,(\idG\tp\delta\tp\idG)(U_l)\, \Psi = 
   \Psi_l\,(\cop\tp\idM\tp\cop)(U_l)
  \end{equation}
thus proving that $(\delta_l,\Psi_l)$ is a two--sided coaction twist
equivalent to $(\delta,\Psi)$ via $U_l$\footnote{
Recall that \eqref{2.36.4} together with $\delta_{l/r} = Ad U_{l/r}
\circ \delta$ already guarantees that $(\delta_l , \Psi_l)$ satisfies
all axioms of Definition \ref{def2.4}.}.
 Using \eqref{2.36.1} and
\eqref{ep2} the l.h.s. of \eqref{2.36.4} gives  
\begin{align}
  &   (\eG\tp U_l\tp\eG)\,(\idG\tp\delta\tp\idG)(U_l)\, \Psi \notag \\
    &\quad = 
  \big[ \eG\tp (\ep\tp\idG\tp\idM\tp\ep\tp\idG)(\Psi) \tp\eG\big]
   \big[ (\idG^{\tp^2}\tp\idM\tp\idG\tp\ep)(\Psi)\tp \eG\big] \notag \\
   & \quad\quad\quad
  \big[ (\ep\tp\cop\tp\idM\tp\idG^{\tp^2})(\Psi)\big] \notag \\
    &\quad = \label{2.36.5}
  \big[(\idG\tp\ep\tp\la\tp\idG\tp\ep)(\Psi) \tp\eG\big]
   \big[ (\idG^{\tp^2}\tp\idM\tp\ep\tp\idG)(\Psi) \tp\eG\big] \\
  &\quad \quad \quad \notag
  \big[ (\ep\tp\cop\tp\idM\tp\idG^{\tp^2})(\Psi)\big],
\end{align}
where in the second equation we have used Eq. \eqref{ep3} acted upon
by $(\idG^{\tp^2} \tp\idM\tp\idG\tp\ep)$. On the other hand, the
r.h.s.\ of \eqref{2.36.4} gives, using 
\eqref{eq215}, \eqref{2.36.1}-\eqref{2.36.3}, \eqref{la1}, \eqref{rho1}
 and \eqref{eq21} 
\begin{align}
&   \Psi_l\,(\cop\tp\idM\tp\cop)(U_l) \notag \\
&\quad = 
\big[ (\idG\tp\ep\tp\la\tp\idG\tp\ep)(\Psi) \tp \eG \big]
       \big[(\ep\tp\idG\tp\la\tp\ep\tp\idG)(\Psi^{-1})\tp\eG\big]
       \notag \\
   &\quad\quad\quad
       \big[(\idG^{\tp^2}\tp\idM\tp\ep\tp\ep)(\Psi) \tp\eG\tp\eG\big]
     \big[ (\ep\tp\ep\tp (\cop\tp\idM)\circ \la
     \tp\idG^{\tp^2})(\Psi)\big]
        \notag \\
 &\quad\quad\quad
     \big[ (\ep\tp\cop\tp\idM\tp\ep\tp\cop)(\Psi)\big] \notag \\
&\quad =
 \big[ (\idG\tp\ep\tp\la\tp\idG\tp\ep)(\Psi) \tp \eG \big]
       \big[(\ep\tp\idG\tp\la\tp\ep\tp\idG)(\Psi^{-1})\tp\eG\big]
       \notag \\
 &\quad\quad\quad
       \big[(\idG^{\tp^2}\tp\idM\tp\ep\tp\ep)(\Psi) \tp\eG\tp\eG\big]
     \big[ (\ep\tp\cop\tp \idM\tp\ep\tp\idG)(\Psi)\tp\eG\big]
       \notag \\
 &\quad\quad\quad
     \big[ (\ep\tp\cop\tp\idM\tp\idG^{\tp^2})(\Psi) \big]\notag \\
&\quad =  \label{2.36.6}
  \big[ (\idG\tp\ep\tp\la\tp\idG\tp\ep)(\Psi) \tp \eG \big]
    \big[(\idG^{\tp^2}\tp\idM\tp\ep\tp\idG)(\Psi) \tp\eG\big]
       \notag \\
 &\quad\quad\quad 
        \big[ (\ep\tp\cop\tp\idM\tp\idG^{\tp^2})(\Psi) \big]
\end{align}
Here we have used Eq. \eqref{ep1} acted upon by
$(\ep\tp\cop\tp\idM\tp\idG^{\tp^2})$ in the second equation and
Eq. \eqref{ep1} acted upon by $(\idG^{\tp^2}\tp\idM\tp\ep\tp\idG)$ in
the last equation. Comparing \eqref{2.36.6} and \eqref{2.36.5} we have
proven \eqref{2.36.4}. Hence $U_l$ provides a twist equivalence from
$(\delta,\Psi)$ to $(\delta_l,\Psi_l)$. Similarly, taking a
left--right mirror image of the above proof, $U_r$ provides a twist
equivalence from $(\delta,\Psi)$ to $(\delta_r,\Psi_r)$.
\end{proof}
\begin{proof}[{\bf Proof of Lemma \ref{lem22}:}]
As has been remarked after Lemma \ref{lem22} it suffices to prove
Eqs.\ \eqref{2.5.26},\eqref{2.5.28} and \eqref{2.5.30}. 
Let us begin with \eqref{2.5.26}. Denoting the multiplication in $\G^{op}$
by $\mu^{op}$ one computes  
\begin{align*}
  &[\eM\tp \Si(\me)]\, q_\rho\, \rho(\mo)
    =  
      [\Fi[1]{\rho} \tp \Si(\alpha \Fi[3]{\rho}\me)\Fi[2]{\rho}]\,
   \rho(\mo) \\
   &\quad \quad =   
   (\idM\tp \mu^{op})\circ (\idM\tp \idG\tp \Si)\Big( 
       [\eM \tp\eG\tp\alpha]\,\phi_\rho \,(\rho\tp\idG)(\rho(m)) \, 
        \Big) \\
     &\quad\quad =(\idM\tp \mu^{op})\circ (\idM\tp \idG\tp \Si)\Big( 
      [\eM \tp\eG\tp\alpha]\,(\idM\tp\cop)(\rho(m)) \,
        \phi_\rho \Big) \\
    &\quad\quad = [m\tp \eG]\, q_\rho,
\end{align*}
where we have plugged in the definition \eqref{2.5.17'} of
$q_\rho$ and used the intertwiner property \eqref{eq23} of $\phi_\rho$
and the antipode property \eqref{eq16}. This proves 
\eqref{2.5.26}.

To prove \eqref{2.5.28} we introduce for $a,b,c \in \G$ the notation
 $\sigma(a\tp b\tp c) : = c \,\Si(\alpha b \beta)\,a$, to compute for
 the l.h.s. 
\begin{align*}
  [\eM\tp \Si(p^2_\rho)]\, q_\rho \, \rho(p^1_\rho) &\equiv 
   \big[\Fi[1]{\rho} \tp \Fib[3]{\rho}\, \Si(\alpha
   \Fi[3]{\rho}\Fib[2]{\rho}\beta)\Fi[2]{\rho}\big]\, 
   \rho(\Fib[1]{\rho}) \\
   &= (\idM\tp\sigma)\Big(
     [\Fi{\rho}\tp\eG]\,
   (\rho\tp\idG\tp\idG)\Fii{\rho} \Big) \\
    & = (\idM\tp\sigma)\Big(
     (\idM\tp\cop\tp\idG)\Fii{\rho} \, [\eM\tp\phi^{-1}]\,
     (\idM\tp\idG\tp\cop)(\Fi{\rho})\Big)  
     \\ 
    & = \eM\tp\Fib[3]{}\Si(\alpha\Fib[2]{}\beta)\Fib[1]{}
                   = \eM\tp\eG,
  \end{align*}
  where we have used the  pentagon identity (\ref{eq24}), then the
 two  antipode properties (\ref{eq16}) together with 
  $(\id\tp\id\tp\ep)(\Fi{\rho}) = \eM\tp\eG$ to drop  the 
  reassociators $\phi_\rho$ and $\phi^{-1}_\rho$ and finally \eqref{eq17}.  
  
The proof of \eqref{2.5.30} is more complicated.
First we  rewrite 
\begin{equation*}
  \text{l.h.s. \eqref{2.5.30}} = \om(X),
\end{equation*}
where 
\begin{equation}
  \label{a2}
 X = 
  [\Fi[1]{\rho} \tp \Fi[2]{\rho} \tp \eG\tp\eG\tp \Fi[3]{\rho}]\,
  [(\rho \tp \idG\tp\idG)(\phi_\rho) \tp \eG] \, 
   [\phi^{-1}_\rho \tp \eG\tp\eG]
\end{equation}
and where the map $\om : \M\tp \G^{\tp^4} \to \M\tp\G^{\tp^2}$ is given
by
\begin{equation*}
  \om(m\tp a \tp b \tp c \tp d) : = m \tp \Si(\alpha d) \,a \tp
  \Si(\alpha c)\,b
\end{equation*}
To rewrite the r.h.s. of \eqref{2.5.30} in the same fashion we first
use the identities \eqref{eq122} and \eqref{eq122a} and the definition
\eqref{2.5.17'} of $q_\rho$ to compute
\begin{align*}
  [\eM\tp h]\, (\idM\tp \cop)(q_\rho) & = 
   [\eM\tp h] \, [\Fi[1]{\rho} \tp \cop\Big(\Si(\alpha \Fi[3]{\rho})
   \Fi[2]{\rho}\Big)] \\
    & = [ \Fi[1]{\rho} \tp (\Si \tp \Si)(\cop^{op}(\Fi[3]{\rho}))] \,
   [\eM\tp h] \, [\eM\tp \cop(\Si(\alpha)\Fi[2]{\rho})] \\ 
  &=
   \big[\eM\tp (\Si\tp\Si)\big(\gamma^{op}\cop^{op}(\Fi[3]{\rho}) \,
   \big)\big] \, [\Fi[1]{\rho}\tp \cop(\Fi[2]{\rho})]. 
\end{align*}
Now we use the formula \eqref{eq116} for $\gamma$ implying
\begin{equation*}
  (\Si\tp\Si)(\gamma^{op}) = \Si(\alpha \Fib[3]{} \Fi[3]{(2)})\Fi[1]{}
  \tp \Si(\alpha \Fib[2]{}\Fi[3]{(1)})\Fib[1]{} \Fi[2]{}
\end{equation*}
to obtain 
\begin{equation*}
  \text{r.h.s. \eqref{2.5.30}} = \om(Y),
\end{equation*}
where
\begin{align}
  \label{a1}
Y &= 
    [\eM\tp\eG\tp \phi^{-1}] \,
    (\idM\tp\idG\tp\idG\tp\cop)\big([\eM\tp\phi] \,
    (\idM\tp\cop\tp\idG)(\phi_\rho)\big)\notag\\
   &\quad\quad\quad\quad\quad 
     ((\id\tp\cop)\circ \rho \tp
    \idG\tp\idG)(\phi_\rho)  
\end{align}
Using the pentagon eq. \eqref{eq24} to replace the second and third
reassociator in \eqref{a1} yields
\begin{align}
    Y& = \notag
    [\eM\tp\eG\tp \phi^{-1}] \,
    (\idM\tp\idG^{\tp^2}\tp\cop)\Big((\id^{\tp^2}\tp\cop)(\phi_\rho) \,
    (\rho\tp\id^{\tp^2})(\phi_\rho)\,[\phi^{-1}_\rho \tp \eG]\Big)  \\ 
   \notag & \quad\quad\quad\quad\quad\quad\quad 
    ((\id\tp\cop)\circ \rho \tp
    \id^{\tp^2})(\phi_\rho)  \\
   & \notag = 
   (\id^{\tp^2} \tp (\cop\tp\id)\circ\cop)(\phi_\rho) \,
   (\rho\tp\id\tp\id)\Big( [\eM\tp \phi^{-1}] \,
   (\id^{\tp^2}\tp\cop)(\phi_\rho)\,(\rho\tp\id^{\tp^2})(\phi_\rho) \Big)
   \\
   \notag&\quad\quad\quad\quad\quad\quad\quad
      [\phi^{-1}_\rho \tp\eG\tp\eG] \\
    &\label{a3} = 
   (\id^{\tp^2} \tp (\cop\tp\id)\circ \cop )(\phi_\rho) \, (\rho \tp
   \cop\tp\id)(\phi_\rho) \, [(\rho\tp\id^{\tp^2})(\phi_\rho) \tp \eG]\,
   [\phi^{-1}_\rho \tp \eG\tp\eG]
\end{align}
where in the second equation we have used \eqref{eq11} and
\eqref{eq23} to shift the reassociators $\phi^{-1}$ and
$\phi^{-1}_\rho$ by one step to the right, and in the third equation
again the pentagon identity \eqref{eq24}.
Hence, when computing $\om(Y)$
the second factor in \eqref{a3}
may be dropped due to the antipode property \eqref{eq16} and the two
coproducts in the first factor disappear by the same reason. Comparing
with \eqref{a2} proves, that 
$\om(X) = \om(Y)$ and therefore
both sides of \eqref{2.5.30} are equal.
\end{proof}
\begin{proof}[{\bf Proof of Lemma \ref{lem22a}}]
  We prove the identity \eqref{2.5.41}. Introducing for $a,b \in \G$
  the map $\nu (a\tp b) : = \Si(\alpha b)\, a$ and using the formula
  \eqref{2.5.17'} for $q_\rho$ we compute  
 \begin{align*}
   (\la\tp\id)(q_{\rho}) \,\phi^{-1}_{\la\rho} &\equiv 
    [\la(\Fi[1]{\rho})\tp\Si(\alpha\Fi[3]{\rho})\Fi[2]{\rho}]\,
    \phi^{-1}_{\la\rho} \\
   &= (\idG\tp\idM\tp\nu)\Big( 
   (\la\tp\id\tp\id)(\Fi{\rho})\, [\phi^{-1}_{\la\rho}\tp\eG]
   \Big) \\
   &= (\idG\tp\idM\tp\nu)\Big( 
      (\id^{\tp^2}\tp\cop)(\phi^{-1}_{\la\rho}) 
     [\eG\tp\Fi{\rho}] \,
      (\id\tp\rho\tp\id)(\Fi{\la\rho}) \Big) \\
   &=[\eG \tp\eG\tp\Si(\Fi[3]{\la\rho})] \, [\eG\tp
      q_\rho ]\, [\Fi[1]{\la\rho}\tp\rho(\Fi[2]{\la\rho})].
 \end{align*}
Here we have plugged in the pentagon equation (\ref{eq27}) and used the fact
that $\phi_{\la\rho}$ may be dropped due to \eqref{27a} and the antipode
property (\ref{eq16}). This proves eq. \eqref{2.5.41} and therefore
Lemma \ref{lem22a} (see remark after Lemma
  \ref{lem22a}).
\end{proof}
\end{appendix}

\end{document}